\documentclass[12pt, a4paper]{article}
\usepackage[left=1in,right=1in,top=1.0in,bottom=1in]{geometry}

\usepackage{amsmath,amssymb,amsthm,amsfonts,afterpage,hyperref}
\usepackage{amscd,subeqnarray}

\usepackage{graphicx}
\usepackage{caption}
\usepackage{subcaption}
\usepackage{authblk}
\usepackage{tikz}

\usepackage{array}
\usepackage{wrapfig}
\usepackage{color}
\usepackage{ragged2e}
\usepackage{stmaryrd}
\usepackage{siunitx}
\usepackage{multirow}
\usepackage{xcolor,cancel}
\usepackage{boldfonts}
\usepackage{symbols}
\usepackage{footmisc}
\usepackage{float}
\usepackage{setspace}
\usepackage{bm}
\usepackage{booktabs}
\usepackage{tikz, xcolor}
\usepackage{soul}
\usepackage{subcaption}

\usepackage[T1]{fontenc}
\usepackage[utf8]{inputenc}
\usepackage{mathrsfs}

\usepackage[linesnumbered,ruled,vlined]{algorithm2e}
\SetKwInput{KwInput}{Input}                
\SetKwInput{KwOutput}{Output} 

\usepackage[sort, numbers]{natbib}
\setcitestyle{square}

\usetikzlibrary{shapes,arrows, positioning,calc}	
\usetikzlibrary{arrows,snakes,backgrounds}

\newcommand{\bfa}[1]{\boldsymbol{#1}} 			
			
\newcommand{\bfeps}{\boldsymbol{\epsilon}}

\newcommand{\Sym}{\text{Sym}}   			%
\newcommand{\tr}{\text{tr}}       				%

\DeclareMathAlphabet{\mathpzc}{OT1}{pzc}{m}{it}

\newcommand{\bfu}{\boldsymbol{u}}

\newcommand{\bfT}{\boldsymbol{T}}		
\newcommand{\bfzero}{\boldsymbol{0}}

\newcommand{\bff}{\boldsymbol{f}}	
	
\newtheorem{theorem}{Theorem}[section]

\newtheorem{lemma}[theorem]{Lemma}

\providecommand{\keywords}[1]
{
  \small	
  \textbf{\textit{Keywords---}} #1
}

\title{Computational Insights into Orthotropic Fracture: Crack-Tip Fields in Strain-Limiting Materials under Non-Uniform Loads\footnote{This work is dedicated to the cherished memory of Professor K. R. Rajagopal (Professor of Mechanical Engineering, Texas A\&M University, College Station, Texas, USA). His seminal contributions to solid mechanics were crucial to this research and profoundly influenced the authors.}}

\author[1]{Saugata Ghosh}
\author[2]{Dambaru Bhatta}
\author[3]{S. M. Mallikarjunaiah} 

\affil[1,2]{School of Mathematical \& Statistical Sciences,\\
The University of Texas - Rio Grande Valley, \\
Edinburg, Texas 78539, USA \\
Email: saugata.ghosh01@utrgv.edu, dambaru.bhatta@utrgv.edu}

\affil[3]{Department of Mathematics \& Statistics,\\
Texas A\&M University-Corpus Christi, \\
Corpus Christi, Texas 78412-5825, USA \\
Email: m.muddamallappa@tamucc.edu (Corresponding Author)}
\date{}

\begin{document}

\maketitle
	    
\begin{abstract}
A finite element framework is presented for analyzing crack-tip phenomena in transversely isotropic, strain-limiting elastic materials. Mechanical response is characterized by an algebraically nonlinear constitutive model, relating stress to linearized strain. Non-physical strain singularities at the crack apex are mitigated, ensuring bounded strain magnitudes. This methodology significantly advances boundary value problem (BVP) formulation, especially for first-order approximate theories. For a transversely isotropic elastic solid with a crack, the governing equilibrium equation, derived from linear momentum balance and the nonlinear constitutive model, is reduced to a second-order, vector-valued, quasilinear elliptic BVP. This BVP is solved using a robust numerical scheme combining Picard-type linearization with a continuous Galerkin finite element method for spatial discretization. Numerical results are presented for various loading conditions, including uniform tension, non-uniform slope, and parabolic loading, with two distinct material fiber orientations. It is demonstrated that crack-tip strain growth is substantially lower than stress growth. Nevertheless, strain-energy density is found to be concentrated at the crack tip, consistent with linear elastic fracture mechanics principles. The proposed framework provides a robust basis for formulating physically meaningful, rigorous BVPs, critical for investigating fundamental processes like crack propagation, damage, and nucleation in anisotropic, strain-limiting elastic solids under diverse loading conditions.
 \end{abstract}
	 
\noindent \keywords{Finite element method; Plane-strain fracture; Transversely Isotropic; Nonlinear constitutive relations}

\section{Introduction}\label{Introduction}

The study of crack-tip fields in orthotropic elastic materials holds immense importance across various advanced engineering disciplines, including aerospace, biomedical applications, civil infrastructure, and composite material design. A key difference with orthotropic materials is their anisotropic nature; unlike isotropic materials, their mechanical properties depend significantly on the direction of applied forces. This inherent directional variability profoundly influences how stress and strain intensify at the tip of a propagating crack. From a fracture mechanics standpoint, this anisotropy introduces considerable complexities, often leading to highly intricate and challenging-to-predict crack propagation paths and failure modes. The exact magnitudes and distributions of these localized crack-tip stresses and strains determine not only a crack's \textit{initiation} but also its subsequent \textit{growth trajectory}, \textit{stability}, and ultimately, the material's \textit{catastrophic failure point}. Without a thorough and precise characterization of these complex crack-tip fields, engineers face significant hurdles in accurately predicting fatigue life, preventing brittle fracture, designing for damage tolerance, and optimizing material selection and orientation \cite{broek2012practical,anderson2017fracture,barsoum1977triangular,barsoum1976degenerate}. In essence, a lack of precise knowledge about these intricate stress and strain distributions at the crack tip poses formidable challenges to ensuring the structural integrity, long-term reliability, and overall safety of components and structures made from such directionally sensitive materials. This highlights the critical need for advanced analytical and computational tools to model these complex phenomena comprehensively \cite{sih1973handbook,broberg1999}.

The complexities of real-world applications frequently subject orthotropic elastic materials to non-uniform loading conditions, creating a critical need for detailed investigations into the resulting crack-tip stresses, strains, and strain energy density. While uniform loading scenarios provide foundational insights, they often fail to capture the nuanced and often more severe stress states induced by concentrated forces, bending moments, or thermal gradients encountered in service. Orthotropic materials, with their directional mechanical properties, respond uniquely to such non-uniform loads, leading to intricate stress and strain distributions around crack tips that can deviate significantly from those predicted by simpler models \cite{wang2013elasticity}. A thorough understanding of crack-tip stresses is essential for identifying potential failure initiation sites and predicting the magnitude of local forces. Simultaneously, analyzing crack-tip strains provides crucial information about the material's local deformation capacity and the onset of plasticity or micro-damage \cite{broberg1999}. Furthermore, the strain energy density criterion, which accounts for both stress and strain, offers a more comprehensive and robust approach to characterizing the energy stored and dissipated at the crack tip, providing a powerful predictor for crack initiation and propagation under complex loading spectra \cite{sih1973handbook, yoon2021quasi}. Neglecting the effects of non-uniform loading on these critical crack-tip parameters in orthotropic materials would lead to substantial limitations in predictive capabilities, potentially compromising the structural integrity and safety of components in high-performance applications.

Accurately characterizing the stress and strain distributions around geometric irregularities like notches, slits, or holes is a cornerstone challenge in both engineering design and theoretical solid mechanics. Historically, the analytical foundation for understanding these critical stress concentrations has been rooted in linearized elasticity theory \cite{Inglis1913, lin1980singular, love2013treatise, murakami1993stress}. However, a key drawback of this classical framework is its inherent prediction of unbounded strain singularities at the tips of such discontinuities---a physically unrealistic outcome stemming from its first-order approximation of finite deformations. This discrepancy has spurred extensive efforts to develop more physically representative constitutive models \cite{gurtin1975, sendova2010, MalliPhD2015, ferguson2015, zemlyanova2012, WaltonMalli2016, rajagopal2011modeling, gou2015modeling}, often integrated with sophisticated numerical techniques like collocation methods \cite{yadav2024gegenbauer, gazonas2023numerical}. Yet, a significant hurdle persists: striking a balance between achieving higher model fidelity and maintaining computational efficiency and experimental validation \cite{broberg1999}. Many proposed theoretical enhancements, while offering greater realism, often incur substantial computational costs or prove difficult to verify empirically.

Furthermore, linear elastic fracture mechanics (LEFM), despite its widespread use in modeling crack initiation and propagation, faces its own inherent limitations. Beyond the well-known strain singularity, LEFM also predicts a physically implausible blunt crack-opening profile and the problematic interpenetration of crack faces, especially in bimaterial interfaces. Crucially, the issue of crack-tip singularity isn't entirely resolved even within various nonlinear elasticity frameworks, as highlighted by works such as \cite{knowles1983} and models incorporating specific constraints, like the bell constraint model \cite{tarantino1997}. These ongoing challenges raise a fundamental question: can algebraic nonlinear constitutive models effectively mitigate the crack-tip strain singularity, even if some stress singularity remains? This inquiry serves as a significant driving force behind the present research.

Traditional elasticity theories, such as those by Cauchy and Green, often fall short when describing material behavior near extreme deformations, particularly at crack tips where unphysical singularities can arise. To address these limitations, Rajagopal and his collaborators have developed a generalized framework for elasticity \cite{rajagopal2003implicit, rajagopal2007elasticity, rajagopal2007response, rajagopal2009class, rajagopal2011non, rajagopal2011conspectus, rajagopal2014nonlinear, rajagopal2016novel, rajagopal2018note}. This innovative approach, often referred to as Rajagopal's theory of elasticity, utilizes implicit constitutive models rooted in robust thermodynamic principles. Within this framework, an elastic body, by definition non-dissipative, is characterized by implicit relationships linking the Cauchy stress and deformation gradient tensors \cite{bustamante2018nonlinear,bustamante2021new,bustamante2015implicit,rodriguez2021stretch}. A particularly powerful aspect of Rajagopal's theory is its ability to yield a hierarchy of 'explicit' nonlinear relationships, allowing linearized strain to be expressed as a nonlinear function of stress. Crucially, a specific subclass of these implicit models possesses a unique 'strain-limiting' property: they can represent linearized strain as a uniformly bounded function across the entire material domain, even under conditions of significant stress. This characteristic makes these models exceptionally well-suited for investigating crack and fracture behavior in brittle materials \cite{rajagopal2011modeling, gou2015modeling, Mallikarjunaiah2015, MalliPhD2015}, offering a path toward analyzing both quasi-static and dynamic crack evolution. The utility of these strain-limiting models has been demonstrated through numerous studies that have revisited and provided new insights into classical elasticity problems \cite{kulvait2013,rajagopal2018bodies,bustamante2009some,bulivcek2014elastic, erbay2015traveling, zhu2016nonlinear, csengul2018viscoelasticity, itou2018states, itou2017contacting, yoon2022CNSNS, yoon2022MMS,rodriguez2022longitudinal}. Their versatility in elucidating the mechanical behavior of a broad spectrum of materials, especially concerning crack and fracture phenomena, is a significant advantage. Recent research, for instance, has shown that formulating quasi-static crack evolution problems within this strain-limiting framework can predict complex crack patterns and even increased crack-tip propagation velocities \cite{lee2022finite, yoon2021quasi}.

This present study applies this advanced theoretical framework to investigate the behavior of a singular crack in a transversely isotropic solid under various forms of the non-uniform tensile load. We develop a specialized constitutive relationship specifically designed to capture the stress-strain response of orthotropic materials accurately. Combining this algebraically nonlinear constitutive equation with the principle of linear momentum balance yields a vector-valued, quasi-linear elliptic boundary value problem. Given the inherent intractability of analytical solutions for such nonlinear partial differential equations, we employ a finite element-based numerical methodology to approximate the solution. The finite element method, widely recognized for its ability to accurately capture crack-tip fields in elastic materials, provides a flexible and robust framework for discretizing the domain and solving the governing partial differential equations. To effectively handle the system's inherent nonlinearities, Picard's iterative algorithm is implemented, with numerical convergence being systematically demonstrated through the progressive reduction of the residual in each iteration. Our findings reveal several intriguing results concerning stress concentration, the controlled growth of crack-tip strains, and a notable decrease in strain-energy density for a single crack subjected to tensile loading. This foundational study opens several promising avenues for future research, including the exploration of thermo-elastic static and quasi-static cracks, as well as the intricate dynamics of crack propagation in transversely isotropic materials.

The structure of this paper is designed to guide the reader through the theoretical development, problem formulation, numerical implementation, and results of our investigation. We begin in Section~\ref{math_formulation}, where the foundational implicit theory is introduced, followed by a detailed derivation of the specific nonlinear constitutive relation employed in this study. Subsequently, Section~\ref{sec:BVP} formally establishes the mathematical model for a static crack situated within a transversely isotropic solid subjected to non-uniform tensile loading, concurrently providing a demonstration of the existence of a unique solution to its weak formulation. Our numerical approach is then elucidated in Section~\ref{sec:fem}, outlining the methodology which combines continuous Galerkin-type finite elements with Picard's iterative algorithm. A comprehensive discussion of the numerical solutions obtained and an analysis of the influence of various model parameters are presented in Section~\ref{sec:rd}. Finally, the paper concludes with a summary of the key findings and their implications in the concluding section.

\section{Mathematical formulation}\label{math_formulation}

The physical setting for our analysis is a bounded, two-dimensional domain, denoted by $\mathcal{D} \subset \mathbb{R}^2$, which is occupied by the material body of a transversely isotropic elastic solid. The boundary of this domain, $\partial \mathcal{D}$, is assumed to be Lipschitz continuous, a condition that ensures the outward unit normal vector, $\bfa{n}$, is well-defined almost everywhere. This boundary is partitioned into two disjoint regions: $\Gamma_D$, on which Dirichlet (displacement) boundary conditions are prescribed, and $\Gamma_N$, on which Neumann (traction) boundary conditions are applied, such that $\partial \mathcal{D} = \overline{\Gamma_N} \cup \overline{\Gamma_D}$. We require that the Dirichlet boundary has a non-zero measure, i.e., $\Gamma_D \neq \emptyset$, to prevent rigid body motion. Within this domain, we introduce a one-dimensional manifold, $\Gamma_c \subset \mathcal{D}$, which represents a pre-existing crack that bifurcates the domain $\mathcal{D}$. The motion of the body is described by the displacement vector field $\bfa{u}: \mathcal{D} \to \mathbb{R}^2$, which maps points from the reference configuration $\bfa{X} \in \mathcal{D}$ to the deformed (current) configuration $\bfa{x} \in \mathcal{D}$. This relationship is formally expressed as $\bfa{u} := \bfa{x} - \bfa{X}$. The mathematical framework is established within the vector space of second-order symmetric tensors, $\Sym(\mathbb{R}^{2 \times 2})$. This space is equipped with the standard inner product $\bfa{A} \colon \bfa{B} = \tr(\bfa{A}^T\bfa{B}) = \sum_{i,j=1}^2 \bfa{A}_{ij} \bfa{B}_{ij}$ and its associated induced Frobenius norm $\| \bfa{A} \| = \sqrt{\bfa{A} \colon \bfa{A}}$. Under the assumption of small deformations, the infinitesimal strain tensor $\bfeps$ is defined as the symmetric part of the displacement gradient:
\begin{equation}
\bfeps(\bfa{u}) := \frac{1}{2} \left( \nabla \bfa{u} + (\nabla \bfa{u})^T \right).
\end{equation}

To ensure a rigorous formulation of the boundary value problem, we introduce the standard Lebesgue and Sobolev function spaces. Let $L^{p}(\mathcal{D})$ for $p\in[1, \infty)$ be the space of $p$-integrable functions, with $L^2(\mathcal{D})$ being the Hilbert space of square-integrable functions equipped with the usual inner product $(\cdot, \cdot)$ and norm $\| \cdot \|$. The Sobolev space $W^{k, p}(\mathcal{D})$ contains functions in $L^p(\mathcal{D})$ whose weak derivatives up to order $k$ are also in $L^p(\mathcal{D})$. For our purposes, the key function space is the classical Sobolev space $H^1(\mathcal{D}) \equiv W^{1, 2}(\mathcal{D})$, defined as:
\begin{equation}
H^{1}(\mathcal{D}) := \left\{ v \in L^{2}(\mathcal{D}) \; \colon \; \nabla v \in (L^{2}(\mathcal{D}))^2 \right\}. \label{eq:H1_def}
\end{equation}
Associated with this is the subspace $H^1_0(\mathcal{D})$, which incorporates homogeneous boundary conditions:
\begin{equation}
H^1_0(\mathcal{D}) := \left\{ v \in H^1(\mathcal{D}) \; \colon \; \left. v \right|_{\partial \mathcal{D}} = 0 \right\}, \label{eq:H10_def}
\end{equation}
where the trace operator implies the boundary condition is satisfied in a generalized sense. Based on these, we define the appropriate spaces for the displacement fields. The space of admissible displacements (trial space) $\bfa{V}$ contains all vector-valued functions that satisfy the essential boundary conditions. The corresponding space of variations (test space) $\bfa{V}_{\bfzero}$ contains functions that vanish on the Dirichlet boundary $\Gamma_D$. These are defined as subspaces of $(H^1(\mathcal{D}))^2$:
\begin{subequations}
\begin{align}
\bfa{V} &:= \left\{ \bfa{u} \in (H^{1}(\mathcal{D}))^2 \colon \; \bfa{u} = \bfa{u}_0 \quad \mbox{on} \;\; \Gamma_D \right\}, \label{eq:trial_space} \\
\bfa{V}_{\bfzero} &:= \left\{ \bfa{u} \in (H^{1}(\mathcal{D}))^2 \colon \; \bfa{u} = \bfa{0} \quad \mbox{on} \;\; \Gamma_D \right\}. \label{eq:test_space}
\end{align}
\end{subequations}

\subsection{A strain-limiting constitutive framework}

The mechanical behavior of the material is described using a nonlinear constitutive framework derived from Rajagopal's theory of implicit elasticity \cite{rajagopal2003implicit,rajagopal2007elasticity}. This theory provides a powerful generalization of classical elasticity, where the relationship between a stress measure and a strain measure is given implicitly. In its general form, the theory posits a functional relationship between the Cauchy stress tensor $\bfa{T}$ and the left Cauchy-Green stretch tensor $\bfa{B}$ as:
\begin{equation}
\widetilde{\bfa{F}}(\bfa{B}, \bfa{T}) = \bfa{0}. \label{eq:implicit_general}
\end{equation}
A prominent subclass of this theory, which forms the basis of our model, expresses the stretch tensor as an explicit function of the stress tensor:
\begin{equation}
\bfa{B} = \widehat{\bfa{F}}(\bfa{T}). \label{eq:explicit_B_of_T}
\end{equation}
By adopting the standard assumption of geometrically linear kinematics (i.e., small displacement gradients), the relationship in \eqref{eq:explicit_B_of_T} simplifies to a nonlinear constitutive law relating the infinitesimal strain $\bfeps$ to the Cauchy stress $\bfa{T}$ via a response function $\bfa{F}$:
\begin{equation}
\bfeps = \bfa{F}(\bfa{T}). \label{eq:strain_from_stress}
\end{equation}
A critical feature of the chosen response function is that it is inherently \textbf{strain-limiting}, meaning the attainable strain is bounded by a material constant $M > 0$, such that $\max_{\bfa{T} \in \Sym} \| \bfa{F}(\bfa{T}) \| \leq M$. This property is physically significant as it precludes the possibility of infinite strains, a non-physical artifact predicted by linear elasticity at stress concentrators like crack tips.

Following recent developments in the field \cite{yoon2021quasi,lee2022finite,itou2017nonlinear}, we adopt a specific algebraic form for the nonlinear response function $\bfa{F}$:
\begin{equation}
\bfa{F}(\bfa{T}) = \frac{\mathbb{K}[\bfa{T}]}{\left(1 + \beta^{\alpha} \| \mathbb{K}^{1/2}[\bfa{T}] \|^{\alpha} \right)^{1/\alpha}}, \quad \text{which implies} \quad \sup_{\bfa{T} \in \Sym} \| \bfa{F}(\bfa{T}) \| \leq \frac{1}{\beta}. \label{eq:specific_F}
\end{equation}
Here, $\alpha > 0$ and $\beta > 0$ are scalar modeling parameters. The parameter $\beta$ governs the degree of nonlinearity and sets the ultimate strain limit; setting $\beta=0$ recovers classical linear elasticity. The parameter $\alpha$ controls the sharpness of the transition to the limiting-strain regime. The fourth-order tensor $\mathbb{K}$ is the compliance tensor, defined as the inverse of the elasticity tensor, $\mathbb{E}$. For a transversely isotropic material, $\mathbb{E}$ is given by:
\begin{equation}
\mathbb{E}[\bfeps] := 2\mu\bfeps + \lambda\,\tr(\bfeps)\,\bfa{I} + \gamma(\bfeps \colon \bfa{M})\,\bfa{M},
\end{equation}
where $\mu > 0$ and $\lambda > 0$ are the Lam\'{e} parameters, $\gamma$ is an additional modulus characterizing the anisotropic response, and the structural tensor $\bfa{M}$ defines the preferred material direction (e.g., fiber orientation) \cite{Mallikarjunaiah2015,MalliPhD2015}.

The tensor-valued function $\bfa{F}$ in \eqref{eq:specific_F} possesses several key mathematical properties \cite{itou2018states} that ensure the well-posedness of the resulting BVP:
\begin{itemize}
    \item[(i)] \textbf{Boundedness:} The function is uniformly bounded, $\| \bfa{F}(\bfa{T}) \| \leq 1/\beta$ for all $\bfa{T} \in \Sym(\mathbb{R}^{2 \times 2})$. This directly enforces the strain-limiting nature of the model.
    \item[(ii)] \textbf{Strict Monotonicity:} The function is strictly monotone, meaning $(\bfa{F}(\bfa{T}_1) - \bfa{F}(\bfa{T}_2)) \colon (\bfa{T}_1 - \bfa{T}_2) > 0$ for all distinct $\bfa{T}_1, \bfa{T}_2 \in \Sym(\mathbb{R}^{2 \times 2})$. This property is fundamental for guaranteeing the uniqueness of the solution to the BVP.
    \item[(iii)] \textbf{Lipschitz Continuity:} The function is Lipschitz continuous, satisfying $\| \bfa{F}(\bfa{T}_1) - \bfa{F}(\bfa{T}_2) \| \leq \widehat{c}_1 \| \bfa{T}_1 - \bfa{T}_2 \|$ for some constant $\widehat{c}_1 > 0$. This ensures the mapping from stress to strain is stable and well-behaved.
    \item[(iv)] \textbf{Coercivity:} A key property for ensuring the existence of a solution within the framework of variational calculus is \textbf{coercivity}. This condition requires that there exists a positive constant, $\widehat{c}_2$, such that:
$$
\left| \bfa{v} \cdot \bfa{F}(\bfa{\Pi}) \bfa{v} \right| \geq \widehat{c}_2 \| \bfa{v} \|^2
$$
This inequality must hold for any symmetric tensor $\bfa{\Pi} \in \mathrm{Sym}(\mathbb{R}^{2 \times 2})$ and any non-zero vector $\bfa{v} \in \mathbb{R}^{2}$. The coercivity constant $\widehat{c}_2$ is not universal; its value depends on the specific material parameters, model choices, and the dimension of the problem. The tensor-valued function $\bfa{F}(\cdot)$ is Coercive. 
\end{itemize}

For sufficiently small values of the nonlinearity parameter $\beta$, the constitutive relation \eqref{eq:strain_from_stress} is invertible \cite{mai2015strong,mai2015monotonicity}. Furthermore, the inverted relationship can be derived from a scalar potential, meaning it is \textit{hyperelastic}. This inverted form, which expresses stress as a function of strain, is computationally advantageous for displacement-based finite element methods and is given by:
\begin{equation}
\bfa{T}(\bfeps) := \Psi\left( \| \mathbb{E}^{1/2}[\bfeps] \| \right) \, \mathbb{E}[\bfeps], \quad \text{where} \quad \Psi(s) = \frac{1}{\left(1 - (\beta s)^{\alpha} \right)^{1/\alpha}}. \label{eq:inverted_hyperelastic}
\end{equation}
In the subsequent analysis, this hyperelastic formulation \eqref{eq:inverted_hyperelastic} will be employed to construct the governing BVP. A primary objective of this work is to systematically compare and contrast the mechanical fields (stress, strain, and strain energy density) at the crack tip as predicted by this nonlinear, strain-limiting model with the predictions derived from its classical linear elastic counterpart (i.e., by setting $\beta=0$).

\section{Governing equations and well-posedness}
\label{sec:BVP}

The study of fracture mechanics in transversely isotropic materials is a cornerstone of modern materials science and structural engineering. The significance of this research area stems from the widespread use of such materials in high-performance and safety-critical applications. This class of materials, characterized by directionally dependent mechanical properties, includes advanced composites, natural substances such as wood and rock, various geological formations, and numerous biological tissues. The initiation and propagation of cracks within these materials can severely compromise structural integrity, potentially leading to catastrophic failure. Consequently, a thorough understanding of crack behavior is essential for designing safe, durable, and reliable structures. In the aerospace sector, for instance, where fiber-reinforced composites are ubiquitous, the ability to predict and arrest crack growth is fundamental to aircraft safety. Likewise, in civil engineering, assessing the long-term durability of structures requires a precise comprehension of fracture processes.

Accordingly, this section formally establishes the governing boundary value problem (BVP) for a cracked, transversely isotropic solid exhibiting strain-limiting behavior. We then outline the key assumptions required to ensure the problem is well-posed and conclude by presenting the weak formulation, which provides the foundation for both the theoretical existence proof and the subsequent numerical approximation. The mechanical system is described by the balance of linear momentum, which, in the absence of inertial effects, dictates that the divergence of the Cauchy stress tensor $\bfa{T}$ is balanced by the body force vector $\bff$. This equilibrium equation is coupled with the inverted, hyperelastic constitutive relationship previously defined, which expresses stress as a nonlinear function of strain, $\bfa{T} = \bfa{T}(\bfeps)$. The system is closed by imposing Dirichlet and Neumann boundary conditions on complementary parts of the domain's boundary. The resulting strong form of the BVP is to find the displacement field $\bfa{u}$ that satisfies:
\begin{subequations}
\label{eq:strong_form_bvp}
\begin{align}
-\nabla \cdot \bfa{T}(\bfeps(\bfa{u})) &= \bfa{f},  && \text{in } \mathcal{D}, \label{eq:equilibrium} \\
\bfa{u} &= \bfa{u}_0, && \text{on } \Gamma_D, \label{eq:dirichlet_bc} \\
\bfa{T}(\bfeps(\bfa{u})) \bfa{n} &= \bfa{g}, && \text{on } \Gamma_N, \label{eq:neumann_bc}
\end{align}
\end{subequations}
where $\bfa{f} \in (L^2(\mathcal{D}))^2$ represents the body force per unit volume and $\bfa{g}$ is the prescribed traction vector on the Neumann boundary $\Gamma_N$.

For the BVP described by \eqref{eq:strong_form_bvp} to be mathematically well-posed, we introduce the following physically motivated assumptions:

\begin{itemize}
    \item[\textbf{A1:}] \textbf{Model and material parameters.} The scalar modeling parameters, $\alpha$ and $\beta$, which control the nonlinear constitutive response, are assumed to be positive constants throughout the domain. Similarly, the Lamé parameters $\mu$ and $\lambda$ are taken as constants, consistent with a homogeneous material. A crucial validation of the model involves verifying numerically that in the limit $\beta \to 0^+$, the predictions of the nonlinear model converge to those of classical linear elasticity.

    \item[\textbf{A2:}] \textbf{Static equilibrium condition.} For the case of a pure traction problem where $\Gamma_D = \emptyset$, the external loads must be self-balanced. This physical necessity imposes the following integral compatibility condition on the source terms:
    \begin{equation}
    \int_{\mathcal{D}} \bfa{f} \, d\bfa{x} + \int_{\partial \mathcal{D}} \bfa{g} \, ds = \bfa{0}.
    \end{equation}
    This ensures that the net force acting on the body is zero, precluding rigid body motion.

    \item[\textbf{A3:}] \textbf{Regularity of boundary data.} The prescribed displacement on the Dirichlet boundary, $\bfa{u}_0$, is assumed to possess sufficient regularity, typically $\bfa{u}_0 \in (W^{1,1}(\mathcal{D}))^2$. This ensures that the prescribed boundary deformation is physically reasonable and mathematically tractable.
\end{itemize}

The existence of a solution to the nonlinear BVP in \eqref{eq:strong_form_bvp} is established by analyzing its corresponding weak (or variational) formulation by following \cite{beck2017existence}. This is derived by multiplying the equilibrium equation \eqref{eq:equilibrium} by an arbitrary test function $\bfa{w}$ from the space of kinematically admissible variations $\bfa{V_0}$ and integrating over the domain $\mathcal{D}$. Applying the divergence theorem and incorporating the Neumann boundary condition \eqref{eq:neumann_bc} yields the following integral statement:\\

\noindent\textbf{Weak formulation.} Find the displacement $\bfa{u} \in \bfa{V}$ such that for all test functions $\bfa{w} \in \bfa{V_0}$:
\begin{equation}
\int_{\mathcal{D}} \bfa{T}(\bfeps(\bfa{u})) \colon \bfeps(\bfa{w}) \, d\bfa{x} = \int_{\mathcal{D}} \bfa{f} \cdot \bfa{w} \, d\bfa{x} + \int_{\Gamma_N} \bfa{g} \cdot \bfa{w} \, ds. \label{eq:weak_form}
\end{equation}
The left-hand side represents the internal virtual work, while the right-hand side represents the external virtual work done by the body forces and surface tractions. The following theorem, stated without proof, formally asserts the existence of a solution.

\begin{theorem}[Existence of a weak solution]
\label{thm:existence}
Let $\mathcal{D} \subset \mathbb{R}^2$ be a bounded Lipschitz domain, given the body force $\bfa{f} \in (L^2(\mathcal{D}))^2$ and $\bfa{g} \in (L^2(\Gamma_N))^2$, and assuming the constitutive law $\bfa{T}(\bfeps)$ satisfies the properties of continuity, monotonicity, and coercivity, and that assumptions A1-A3 hold, then there exists at least one solution $\bfa{u} \in \bfa{V}$ to the weak formulation \eqref{eq:weak_form}.
\end{theorem}

This theorem provides the essential theoretical guarantee that the formulated BVP is solvable, paving the way for the development of stable and convergent numerical methods, such as the finite element method, for its approximation. The theoretical underpinning for the existence of a solution to our boundary value problem is provided by the framework established in the work of Beck et al.~\cite{beck2017existence}. The existence theorem presented therein is directly applicable to our current formulation, with only a minor modification in the algebraic form of the Cauchy stress tensor, as specified in equation~\eqref{eq:inverted_hyperelastic}. This congruence holds because our chosen constitutive model preserves the essential mathematical properties, such as strict monotonicity and coercivity, that are central to the theorem's proof. We therefore posit that the material properties ($\mu, \lambda$), the constitutive modeling parameters ($\alpha, \beta$), and the problem data ($\bfa{f}, \bfa{g}, \bfa{u}_0$) all satisfy the necessary conditions stipulated in~\cite{beck2017existence}. With these requirements met, we can confidently assert the existence of a unique solution pair $(\bfu, \bfT)$ residing in the function space $(W^{1,1}(\mathcal{D}))^2 \times \Sym({L}^1(\mathcal{D})^{2 \times 2})$ for the weak formulation previously delineated. The formal proof for our specific strain-limiting model is a direct adaptation of the arguments presented in~\cite{beck2017existence}, providing a rigorous mathematical foundation for the numerical investigations that follow.

\section{Finite element discretization and numerical scheme}
\label{sec:fem}

This section details the numerical approximation of the previously established boundary value problem using a conforming Galerkin finite element method. We begin by presenting the continuous variational formulation, discuss the properties of the resulting weak solution, and then construct the discrete problem that leads to a solvable algebraic system.

\subsection{Variational formulation and solution properties}

The foundation of the finite element method lies in recasting the strong form of the BVP, given by the partial differential equation in \eqref{eq:equilibrium}, into an equivalent integral form, known as the weak or variational formulation. This is achieved through the standard method of weighted residuals. We multiply the governing equilibrium equation by an arbitrary test function $\bfa{v}$ from the space of admissible variations $\bfa{V}_{\bfzero}$ (defined in \eqref{eq:test_space}) and integrate over the domain $\mathcal{D}$. Applying Green's formula (integration by parts) and incorporating the Neumann boundary condition \eqref{eq:neumann_bc} allows us to arrive at the following continuous variational problem.

\vspace{\baselineskip}
\noindent\textbf{Problem (continuous weak formulation):}
Given the prescribed data and material parameters, find the displacement field $\bfa{u} \in \bfa{V}$ such that:
\begin{equation}
\label{eq:weak_formulation}
a(\bfa{u}; \bfa{v}) = L(\bfa{v}), \quad \forall\, \bfa{v} \in \bfa{V}_{\bfzero},
\end{equation}
where the semilinear form $a(\cdot; \cdot)$ and the linear functional $L(\cdot)$ are defined as follows:
\begin{subequations}
\label{def:A-L}
\begin{align}
a(\bfa{u}; \bfa{v}) &:= \int_{\mathcal{D}} \Psi\left( { \|\mathbb{E}^{1/2} [\bfeps(\bfa{u})]\|}\right) \, \mathbb{E}[\bfeps(\bfa{u})] \colon \bfeps(\bfa{v}) \, d\bfa{x}, \\
L(\bfa{v}) &:= \int_{\mathcal{D}} \bfa{f} \cdot \bfa{v} \, d\bfa{x} + \int_{\Gamma_N} \bfa{g} \cdot \bfa{v} \, ds.
\end{align}
\end{subequations}
Note that the form $a(\cdot; \cdot)$ is nonlinear in its first argument due to the dependence of the function $\Psi$ on the solution $\bfa{u}$.

Under assumptions of sufficient smoothness on the problem data—specifically, for $\bfa{f} \in (L^{2}(\mathcal{D}))^2$, $\bfa{u}_0 \in (H^{1/2}(\Gamma_D))^2$, and $\bfa{g} \in (L^{2}(\Gamma_N))^2$—it can be shown that the variational problem \eqref{eq:weak_formulation} admits a unique weak solution $\bfa{u}$ with enhanced regularity. This solution belongs to the space $\mathcal{U}_s := \{ \bfa{v} \in (H^2(\mathcal{D}) \cap W^{1, \infty}(\mathcal{D}))^2 \colon \left. \bfa{v} \right|_{\Gamma_D} = \bfa{u}_0 \}$. Furthermore, the solution satisfies the following \textit{a priori} stability estimate:
\begin{equation}
\| \bfa{u} \|_{H^2(\mathcal{D})} \leq \widehat{c} \left( \| \bfa{f} \|_{L^2(\mathcal{D})} + \| \bfa{u}_0 \|_{H^{1/2}(\Gamma_D)} + \| \bfa{g} \|_{L^2(\Gamma_N)} \right), \label{eq:stability_estimate}
\end{equation}
where $\widehat{c} > 0$ is a regularity constant independent of the solution and the data. This estimate is crucial as it confirms that the continuous problem is well-posed.

\vspace{\baselineskip}
\noindent\textbf{Remark (numerical treatment of nonlinearity):} The constitutive nonlinearity introduced by the function $\Psi(\cdot)$ requires an iterative approach for its numerical solution. In this work, we employ Picard's iterative method (a fixed-point iteration), which is a natural choice for this class of semilinear problems. The convergence of Picard's method is, however, only guaranteed for a sufficiently close initial guess. To ensure a robust solution strategy, we first solve the corresponding linear elastic problem (i.e., by setting $\beta = 0$, which makes $\Psi \equiv 1$) and use its solution as the initial estimate for the subsequent nonlinear iterations.

\subsection{Galerkin discretization and properties}

To obtain a numerical approximation, we discretize the domain $\overline{\mathcal{D}}$ using a shape-regular family of partitions, $\{\mathcal{T}_h\}_{h > 0}$, consisting of non-overlapping quadrilateral elements $\mathcal{K}$. The mesh parameter $h$ is defined as $h := \max_{\mathcal{K} \in \mathcal{T}_h} \operatorname{diam}(\mathcal{K})$. The collection of all boundary edges is denoted by $\mathscr{E}_{bd,h}$ and is partitioned into Dirichlet and Neumann sets, $\mathscr{E}_{D,h}$ and $\mathscr{E}_{N,h}$, respectively.

We then define a finite-dimensional subspace $\bfa{V}_h \subset (H^1(\mathcal{D}))^2$ and $\bfa{V}_{h, \, 0} \subset \bfa{V}_{\bfzero}$ for the trial and test functions. For this study, we use the space of continuous vector-valued functions that are piecewise bilinear on each element:
\begin{equation}
\bfa{V}_{h} := \left\{ \bfa{u}_h \in (C^{0}(\overline{\mathcal{D}}))^2 \colon \left. \bfa{u}_h \right|_{\mathcal{K}} \in (\mathbb{Q}_{1}(\mathcal{K}))^2, \forall \mathcal{K} \in \mathcal{T}_h \right\}.
\end{equation}
Here, $\mathbb{Q}_{1}(\mathcal{K})$ is the space of bilinear polynomials on element $\mathcal{K}$. The continuous Galerkin method seeks an approximate solution within this discrete space by restricting the weak formulation \eqref{eq:weak_formulation} to $\bfa{V}_h$.

\vspace{\baselineskip}
\noindent\textbf{Problem (discrete weak formulation):}
Find the discrete displacement $\bfa{u}_h \in \bfa{V}_{h}$ (satisfying the discrete boundary conditions) such that for all test functions $\bfa{v}_h \in \bfa{V}_{h,0}$:
\begin{equation}\label{eq:discrete_wf}
a(\bfa{u}_h; \bfa{v}_h) = L(\bfa{v}_h).
\end{equation}
The discrete forms are computed by summing the contributions from each element:
\begin{align}\label{}
a(\bfa{u}_h; \bfa{v}_h) &= \sum_{\mathcal{K} \in \mathcal{T}_h} \int_{\mathcal{K}} \Psi\left( \|\mathbb{E}^{1/2}[\bfeps(\bfa{u}_h)]\| \right) \, \mathbb{E}[\bfeps(\bfa{u}_h)] \colon \bfeps(\bfa{v}_h) \, d\bfa{x}, \label{eq:discrete_a} \\
L(\bfa{v}_h) &= \sum_{\mathcal{K} \in \mathcal{T}_h} \int_{\mathcal{K}} \bfa{f} \cdot \bfa{v}_h \, d\bfa{x} + \sum_{e \in \mathscr{E}_{N,h}} \int_{e} \bfa{g} \cdot \bfa{v}_h \, ds. \label{eq:discrete_L}
\end{align}
The well-posedness of this discrete nonlinear system relies on the mathematical properties of the semilinear form $a(\cdot; \cdot)$. The following two lemmas state its Lipschitz continuity and strong monotonicity, which are essential for proving the existence and uniqueness of the discrete solution and for guaranteeing the convergence of iterative solvers.

\begin{lemma}[Lipschitz continuity]
\label{lem:lipschitz}
For any $\bfa{u}_1, \bfa{u}_2, \bfa{w} \in \bfa{V}_h$, there exists a constant $k_1 > 0$, independent of the mesh size $h$, such that the semilinear form $a(\cdot; \cdot)$ is Lipschitz continuous in its first argument:
\begin{equation}
|a(\bfa{u}_1; \bfa{w}) - a(\bfa{u}_2; \bfa{w})| \leq k_1 \| \bfa{u}_1 - \bfa{u}_2 \|_{H^1} \| \bfa{w} \|_{H^1}.
\end{equation}
\end{lemma}

\begin{lemma}[Strong monotonicity]
\label{lem:monotonicity}
For any $\bfa{u}_1, \bfa{u}_2 \in \bfa{V}_h$, there exists a constant $k_2 > 0$, independent of the mesh size $h$, such that the semilinear form $a(\cdot; \cdot)$ satisfies the strong monotonicity condition:
\begin{equation}
a(\bfa{u}_1; \bfa{u}_1 - \bfa{u}_2) - a(\bfa{u}_2; \bfa{u}_1 - \bfa{u}_2) \geq k_2 \| \bfa{u}_1 - \bfa{u}_2 \|_{H^1}^2.
\end{equation}
\end{lemma}
With the fundamental properties of the discrete operator established, we can now address the well-posedness of the numerical scheme. The conditions of Lipschitz continuity and strong monotonicity, proven in the preceding lemmas, are precisely the hypotheses required to apply foundational results from nonlinear functional analysis. Specifically, these two properties together ensure that the operator associated with the discrete problem is a contraction mapping. Therefore, the existence and uniqueness of a solution are a direct consequence of the \textit{Banach Fixed-Point Theorem}. This pivotal result provides a rigorous guarantee that our discrete formulation is sound. We consolidate this conclusion in the following theorem.

\begin{theorem}
\label{thm:discrete_existence}
Assuming the semilinear form $a(\cdot; \cdot)$ satisfies the Lipschitz continuity and strong monotonicity conditions detailed in Lemma~\ref{lem:lipschitz} and Lemma~\ref{lem:monotonicity}, the discrete variational problem~\eqref{eq:discrete_wf} admits a unique solution $\bfa{u}_h \in \bfa{V}_h$.
\end{theorem}

The proof of Theorem~\ref{thm:discrete_existence} is a standard application of the arguments underlying the fixed-point theorem (see, for instance, \cite{manohar2024hp}). This result confirms that the final algebraic system derived from our finite element discretization is non-singular and has a unique solution, thereby justifying the computational approach.

\section{Computational study and discussion}
\label{sec:rd}

This section presents a series of numerical experiments designed to validate the proposed strain-limiting framework and to quantify its impact on crack-tip mechanics. We investigate the behavior of a single crack embedded in a transversely isotropic solid, contrasting the predictions of our nonlinear model against those of classical linear elasticity. The primary objective is to demonstrate how the inherent strain-limiting property of the model mitigates the non-physical singularities predicted by linear theory, resulting in a more realistic depiction of the near-tip fields.

For the numerical approximation, we employ a conventional continuous Galerkin finite element method using bilinear quadrilateral elements. While more advanced discretization techniques exist, this standard approach is sufficient and well-suited for the present goal of demonstrating the fundamental physical differences between the two constitutive models. A comprehensive \textit{a priori} error analysis and the application of more sophisticated discretization schemes are deferred to future work, for which the analysis developed in \cite{manohar2024hp} may serve as a valuable starting point. All simulations were implemented using the open-source, object-oriented finite element library \texttt{deal.II}~\cite{2023dealii, dealii2019design} and were performed on structured computational meshes.

The constitutive nonlinearity is resolved using the iterative procedure detailed in Algorithm~\ref{alg:picard}. To assess convergence at each Picard iteration $n$, we compute the norm of the residual functional, which measures the imbalance in the weak formulation:
\begin{equation}\label{eq:residual}
\mathcal{R}(\bfa{u}^n_h; \bfa{\varphi}_h) := a(\bfa{u}^n_h; \bfa{\varphi}_h).
\end{equation}
The iteration is terminated when the norm of this residual falls below a prescribed tolerance. The complete computational procedure is outlined below.

\begin{algorithm}[H]
\SetAlgoLined
\KwInput{Finite element mesh $\mathcal{T}_h$; material and model parameters ($\mu, \lambda, \gamma, \alpha, \beta$); max iterations $M_{max}$; tolerance $TOL$; boundary data $\bfa{g}, \bfa{u}_0$.}
\KwOutput{Converged discrete solution $\bfa{u}_h$.}

\tcp{Step 1: Initialization}
Set iteration counter $n \leftarrow 0$. \\
Compute the initial guess $\bfa{u}^0_h$ by solving the linear elastic problem (i.e., setting $\beta=0$ in the weak form \eqref{eq:discrete_wf}). 

\tcp{Step 2: Picard Iteration Loop}
\While{$n < M_{max}$}{
  $n \leftarrow n + 1$. \\
  Assemble the linear system arising from the weak form, using the previous iterate $\bfa{u}^{n-1}_h$ to evaluate the nonlinearity. \\
  Solve the resulting algebraic system for the new iterate $\bfa{u}^n_h$.

  \tcp{Step 3: Convergence Check}
  Compute the norm of the residual, e.g., $\| \mathcal{R}(\bfa{u}^n_h) \|$. \\
  \If{$\| \mathcal{R}(\bfa{u}^n_h) \| \leq TOL$}{
    Break loop.
  }
}

\tcp{Step 4: Output}
Return the converged solution $\bfa{u}_h \leftarrow \bfa{u}^n_h$. \\
Perform postprocessing (e.g., calculate stress, strain, strain energy).

\caption{Iterative Picard scheme for the nonlinear fracture problem.}
\label{alg:picard}
\end{algorithm}

\subsection{Benchmark problem: Mode-I static crack in a plate}

\begin{figure}[H]
  \centering
  \includegraphics[width=0.5\linewidth]{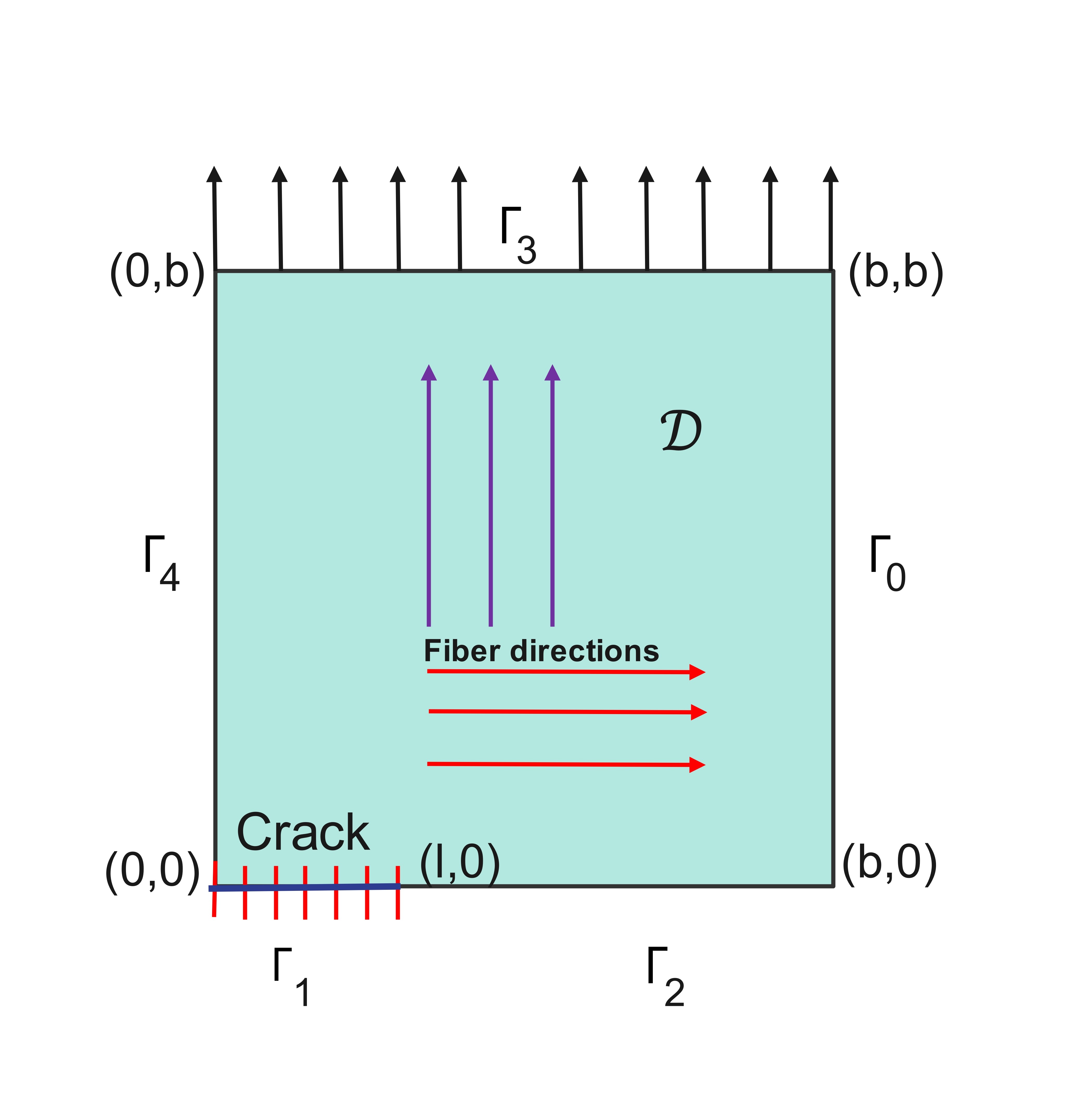}
  \caption{Schematic of the computational domain and boundary conditions for the Mode-I crack problem.}
  \label{fig:cd}
\end{figure}

The central numerical experiment involves assessing the model's performance on a benchmark problem: a plate with a single edge crack subjected to Mode I tensile loading. Our goal is to compare the mechanical fields predicted by our nonlinear model with the standard linear elastic solution, thereby highlighting the physical realism introduced by the strain-limiting constitutive model. For all simulations, the Picard iteration was executed with a convergence tolerance of $TOL = 10^{-6}$ and a maximum of $M_{max} = 10$ iterations, which was sufficient to achieve convergence in all cases presented. The geometry and boundary conditions for this problem are illustrated in Figure~\ref{fig:cd}. The computational domain is a rectangular plate with a horizontal crack extending from the left edge along the $x$-axis ($0 \leq x \leq 1, y=0$). The loading and constraints are as follows:
\begin{itemize}
    \item The top boundary, denoted by $\Gamma_{3}$, is subjected to a combination of a linearly varying slope load and a non-uniform load.
    \item The right boundary ($\Gamma_{0}$) is traction-free.
    \item The left boundary ($\Gamma_{4}$) is constrained from horizontal movement ($u_1 = 0$).
    \item The bottom boundary ($\Gamma_{2}$) is constrained from vertical movement ($u_2 = 0$), enforcing symmetry.
    \item The crack boundary ($\Gamma_{1}$) is kept traction-free. 
\end{itemize}

To explore the influence of material anisotropy, we analyze the system's response for two distinct fiber orientations. In transversely isotropic materials like fiber-reinforced composites or bone, the axis of symmetry typically aligns with the fiber direction. This orientation is critical as it governs the directional stiffness and profoundly influences the stress distribution and strain mitigation mechanisms in the vicinity of the crack tip.

A defining characteristic of transversely isotropic materials is their directionally dependent mechanical response, which is fundamentally governed by their internal microstructure. In our mathematical framework, this anisotropy is explicitly captured through the orientation of the structural tensor, $\bfa{M}$, within the fourth-order elasticity tensor, $\mathbb{E}$. The alignment of this structural tensor dictates the material's principal axis of stiffness. To systematically investigate the influence of this anisotropy on crack-tip phenomena, our study explores two distinct scenarios for the material's fiber orientation. This approach is physically motivated, as it is well-established that the axis of material symmetry in materials such as fiber-reinforced composites, wood, and bone is determined by the alignment of their constituent fibers or other microstructural features. The orientation of these internal structures plays a crucial role in determining the material's overall mechanical behavior. It fundamentally governs how stresses are distributed and concentrated under applied loads, particularly in the critical region surrounding a crack tip. By analyzing these different configurations, we aim to gain critical insight into how material anisotropy can influence stress shielding, strain localization, and ultimately, the fracture resistance of the solid.

\subsection{Case 1(a)-Slope load: Fibers aligned with the crack  plane}
In our first numerical investigation, we consider a cracked material with reinforcing fibers oriented parallel to the $x$-axis, subjected to a linearly varying slope-type load on its top surface. This loading condition is critical as it realistically simulates many practical scenarios, such as the bending of beams or thermal gradients, which induce linearly distributed stresses. Applying this load enables a more accurate simulation of the stress and strain fields at the crack tip, which is crucial for calculating critical parameters such as a structure's fracture toughness and predicting crack propagation.

For this simulation, the principal axis of material stiffness is aligned with the crack plane. This configuration is significant because it directly influences the crack opening behavior and stress distribution. To implement this physical setup in our computational model, we define the structural tensor as $\mathbf{M} = \mathbf{e}_1 \otimes \mathbf{e}_1$, where $\mathbf{e}_1$ is the unit vector in the $x$-direction. For a particular value of the material's Poisson's ratio constant, we solved the resulting nonlinear system of equations using the Picard iterative scheme. The numerical solver was robust and efficient, demonstrating rapid and monotonic convergence as documented in Table~\ref{table1}. This confirms the stability of our iterative method for this specific fiber orientation.

\begin{table}[H]
\centering
\begin{tabular}{|c|c|}
\hline
Iteration No & Residual \tabularnewline
\hline 
1 & 3.00312e-05 \tabularnewline
\hline 
2 & 2.14996e-06 \tabularnewline
\hline 
3 & 2.66311e-07 \tabularnewline
\hline 
4 & 2.09593e-07 \tabularnewline
\hline 
5 & 2.10718e-07 \tabularnewline
\hline 
6 & 2.10626e-07 \tabularnewline
\hline 
7 & 2.10632e-07 \tabularnewline
\hline 
\end{tabular}
 \caption{Residual computed at each iteration for the case of fiber's orientation is along the plane of the crack and a slope-type top load.}
   \label{table1}
\end{table}
Since the residual value stagnated after the 7th iteration and remained unchanged up to the 100th iteration, a different solver with an appropriate preconditioner is preferred to reduce the residual even further. Investigating this combination will be a topic for future study.

\begin{figure}[H]
    \centering
    \begin{subfigure}{0.3\linewidth}
        \centering
        \includegraphics[width=\linewidth]{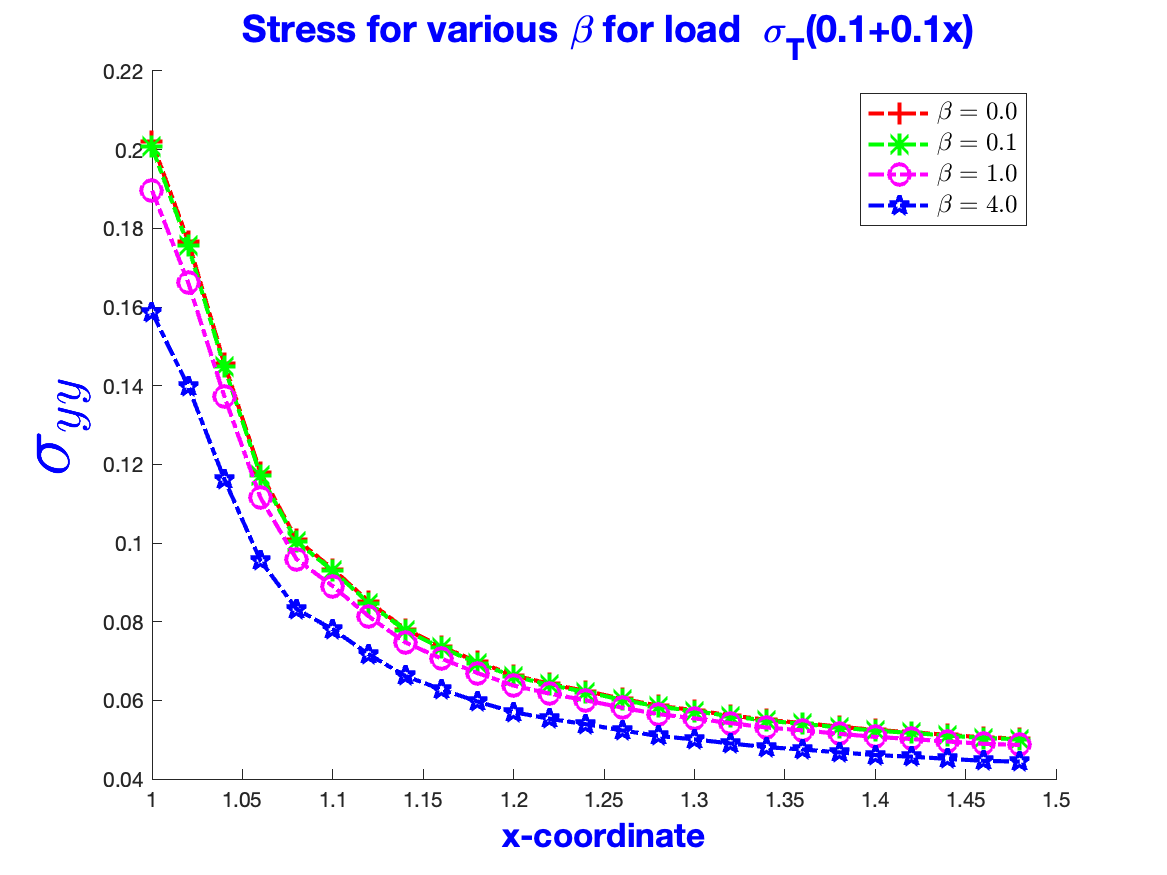}
        \caption{Stress for various $\beta$ for $\alpha = 1.0$ and $\sigma_{T} = 0.1$}
        \label{fig:stress_beta}
    \end{subfigure}
    \hfill
    \begin{subfigure}{0.3\linewidth}
        \centering
        \includegraphics[width=\linewidth]{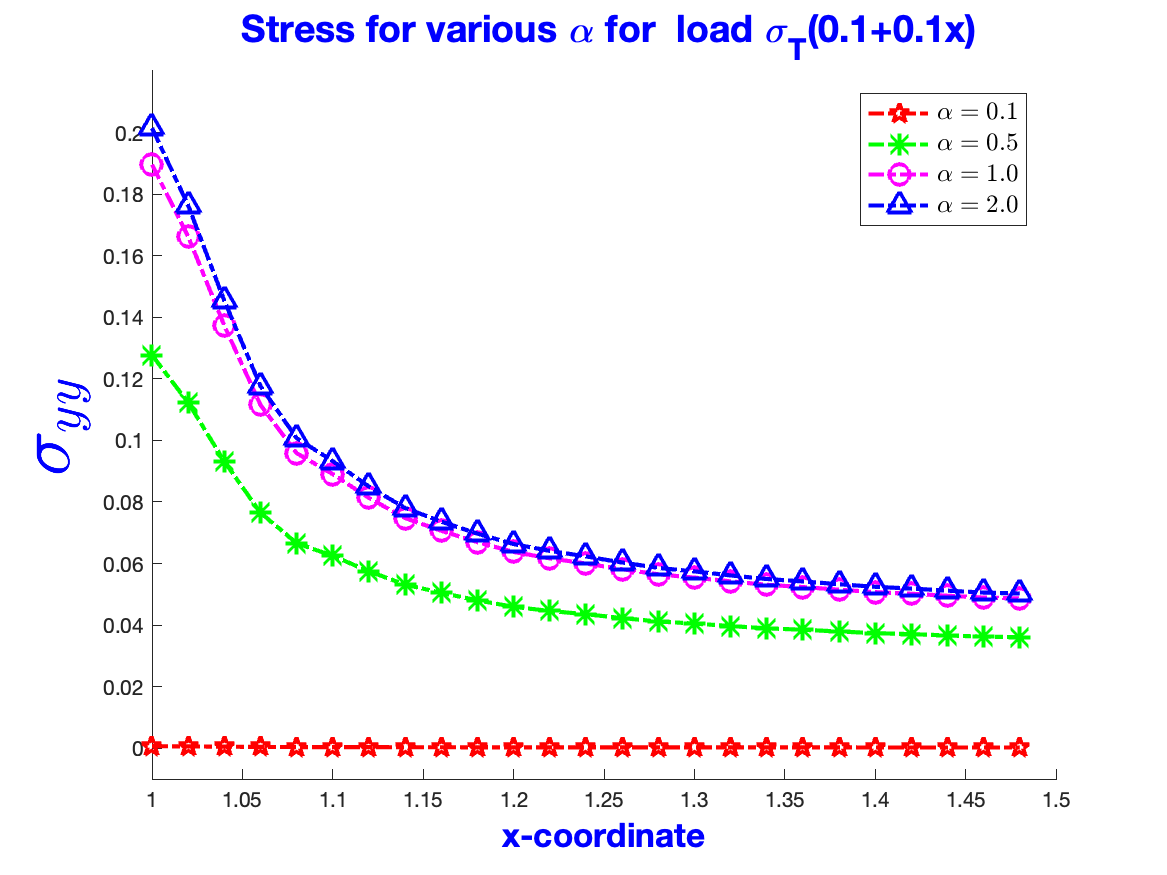}
        \caption{Stress for various $\alpha$ for $\beta = 1.0$ and $\sigma_{T} = 0.1$}
        \label{fig:stress_alpha}
    \end{subfigure}
    \hfill
    \begin{subfigure}{0.3\linewidth}
        \centering
        \includegraphics[width=\linewidth]{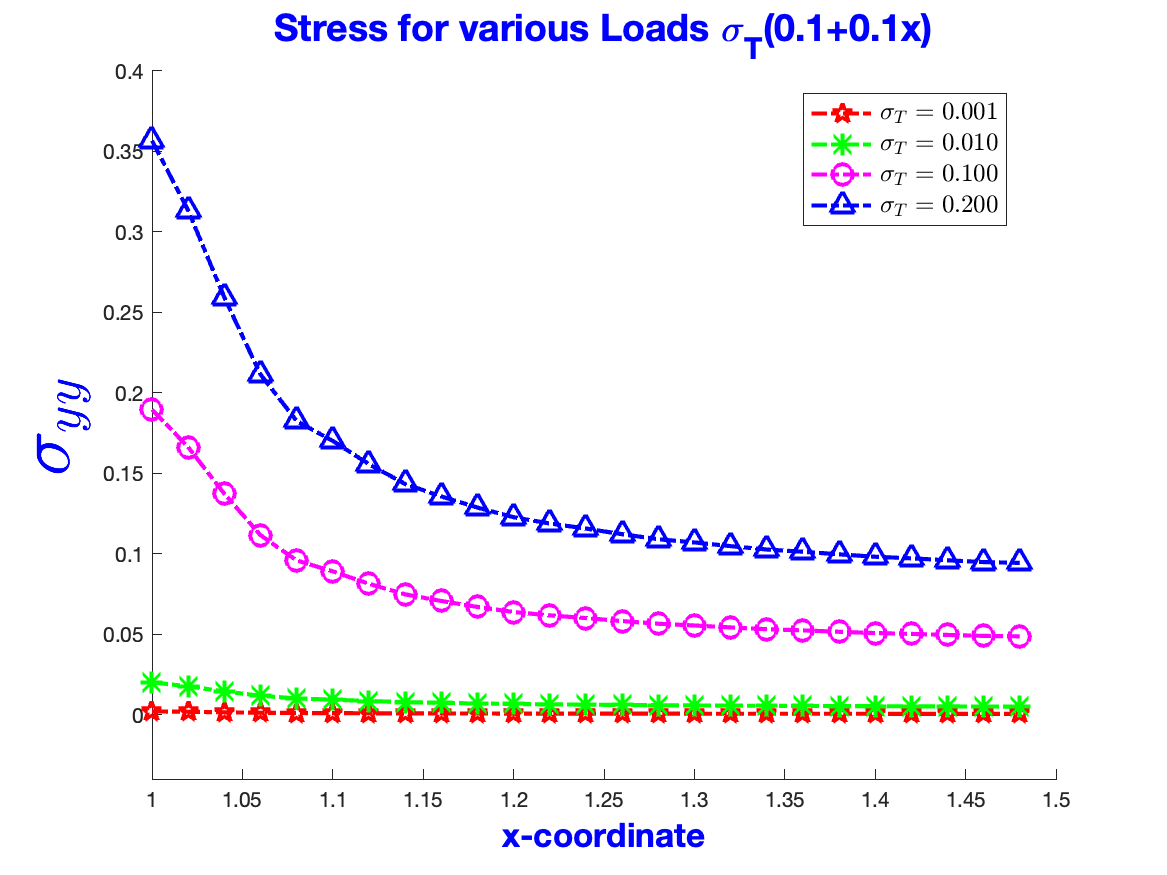}
        \caption{Stress for various $\sigma_{T}$ for $\beta = 1.0$ and $\alpha = 1.0$}
        \label{fig:stress_sigma}
    \end{subfigure}
    \caption{Stress plots for different parameter variations for slope loads $x$-direction.}
    \label{fig:stress_case1a}
\end{figure}

\begin{figure}[H]
    \centering
    \begin{subfigure}{0.3\linewidth}
        \centering
        \includegraphics[width=\linewidth]{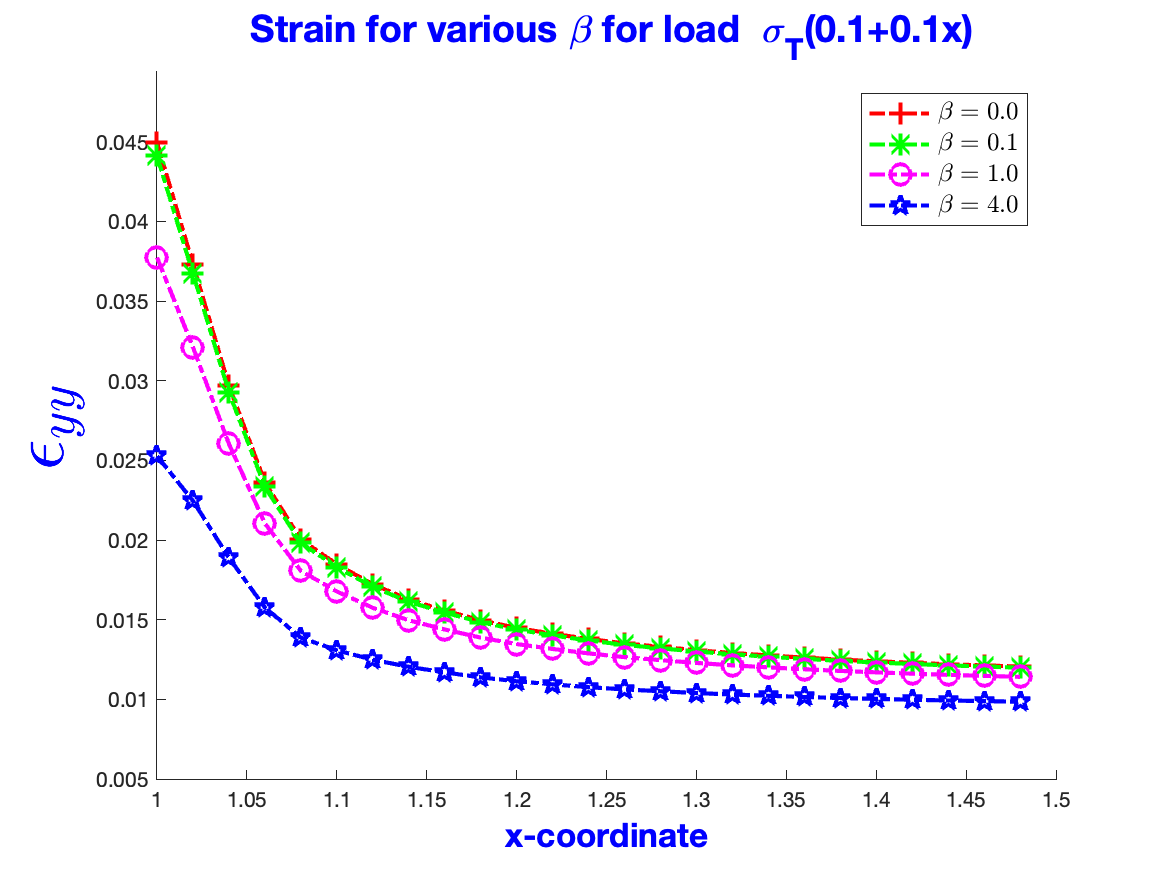}
        \caption{Strain for various $\beta$ for $\alpha = 1.0$ and $\sigma_{T} = 0.1$}
        \label{fig:strain_beta}
    \end{subfigure}
    \hfill
    \begin{subfigure}{0.3\linewidth}
        \centering
        \includegraphics[width=\linewidth]{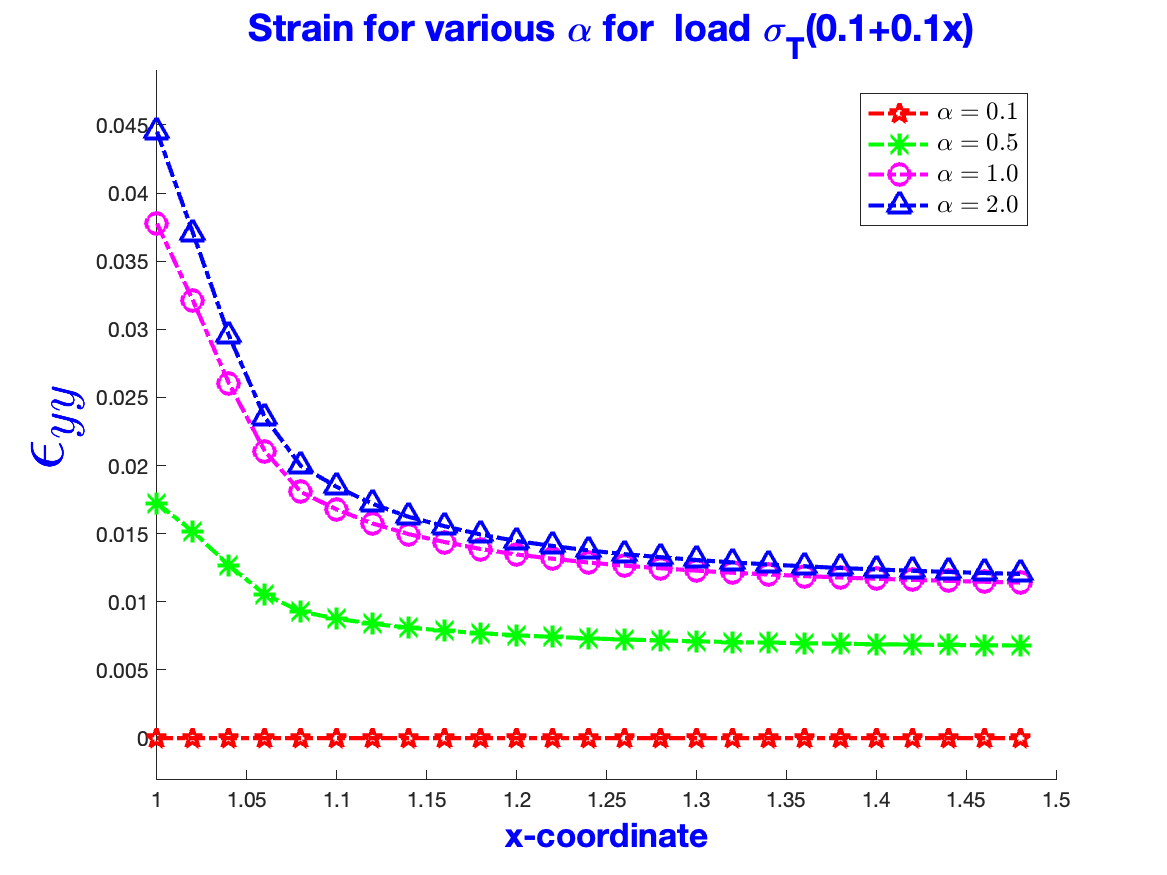}
        \caption{Strain for various $\alpha$ for $\beta = 1.0$ and $\sigma_{T} = 0.1$}
        \label{fig:strain_alpha}
    \end{subfigure}
    \hfill
    \begin{subfigure}{0.3\linewidth}
        \centering
        \includegraphics[width=\linewidth]{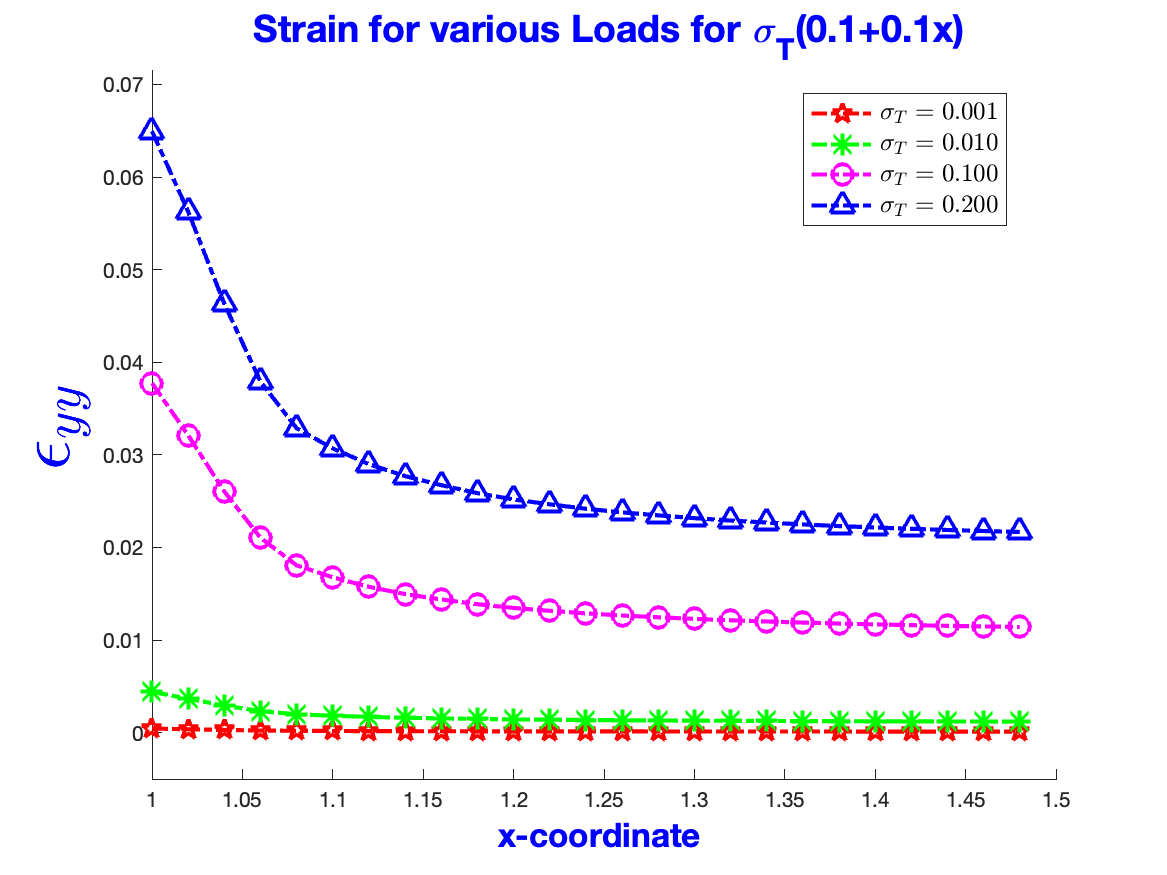}
        \caption{Strain for various $\sigma_{T}$ for $\beta = 1.0$ and $\alpha = 1.0$}
        \label{fig:strain_sigma}
    \end{subfigure}
    \caption{Strain plots for different parameter variations for slope loads $x$-direction.}
    \label{fig:strain_case1a}
\end{figure}

\begin{figure}[H]
    \centering
    \begin{subfigure}{0.3\linewidth}
        \centering
        \includegraphics[width=\linewidth]{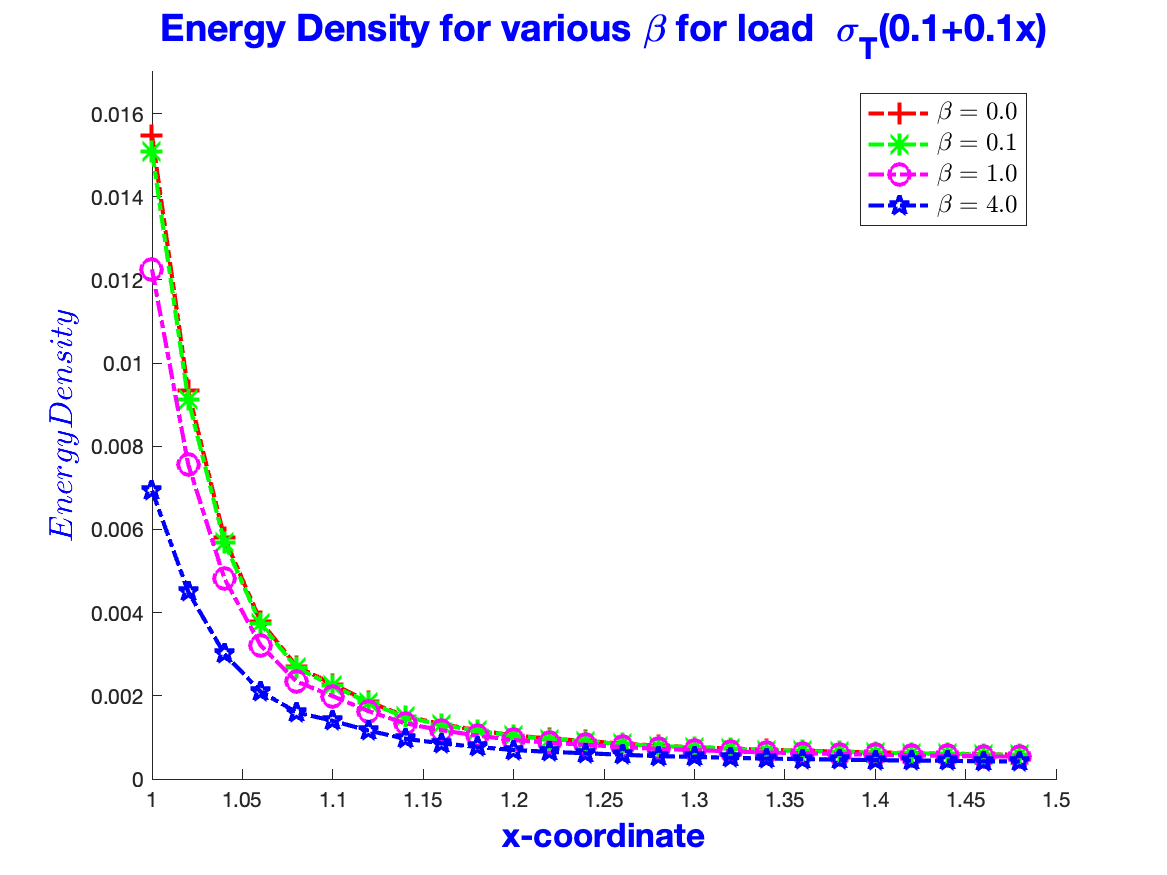}
        \caption{Energy density for various $\beta$ for $\alpha = 1.0$ and $\sigma_{T} = 0.1$}
        \label{fig:energy_density_beta}
    \end{subfigure}
    \hfill
    \begin{subfigure}{0.3\linewidth}
        \centering
        \includegraphics[width=\linewidth]{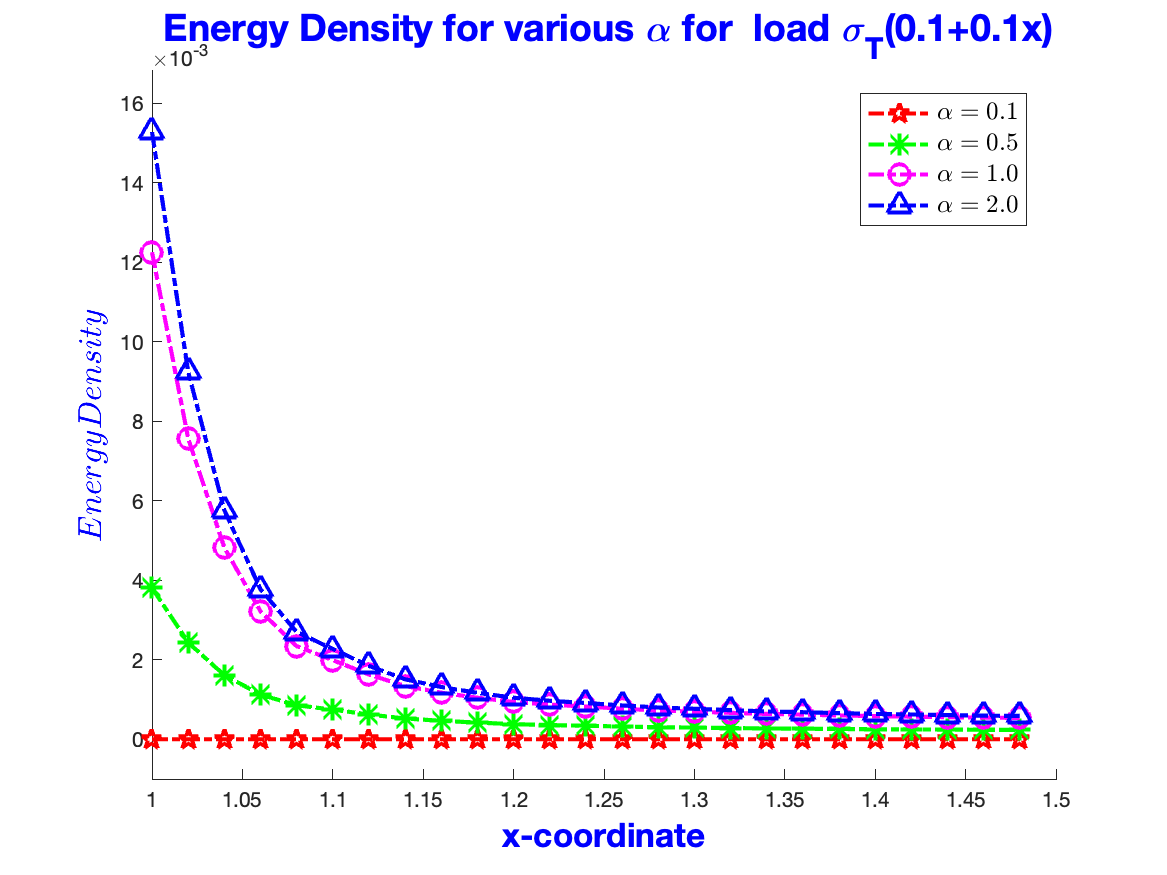}
        \caption{Energy density for various $\alpha$ for $\beta = 1.0$ and $\sigma_{T} = 0.1$}
        \label{fig:energy_density_alpha}
    \end{subfigure}
    \hfill
    \begin{subfigure}{0.3\linewidth}
        \centering
        \includegraphics[width=\linewidth]{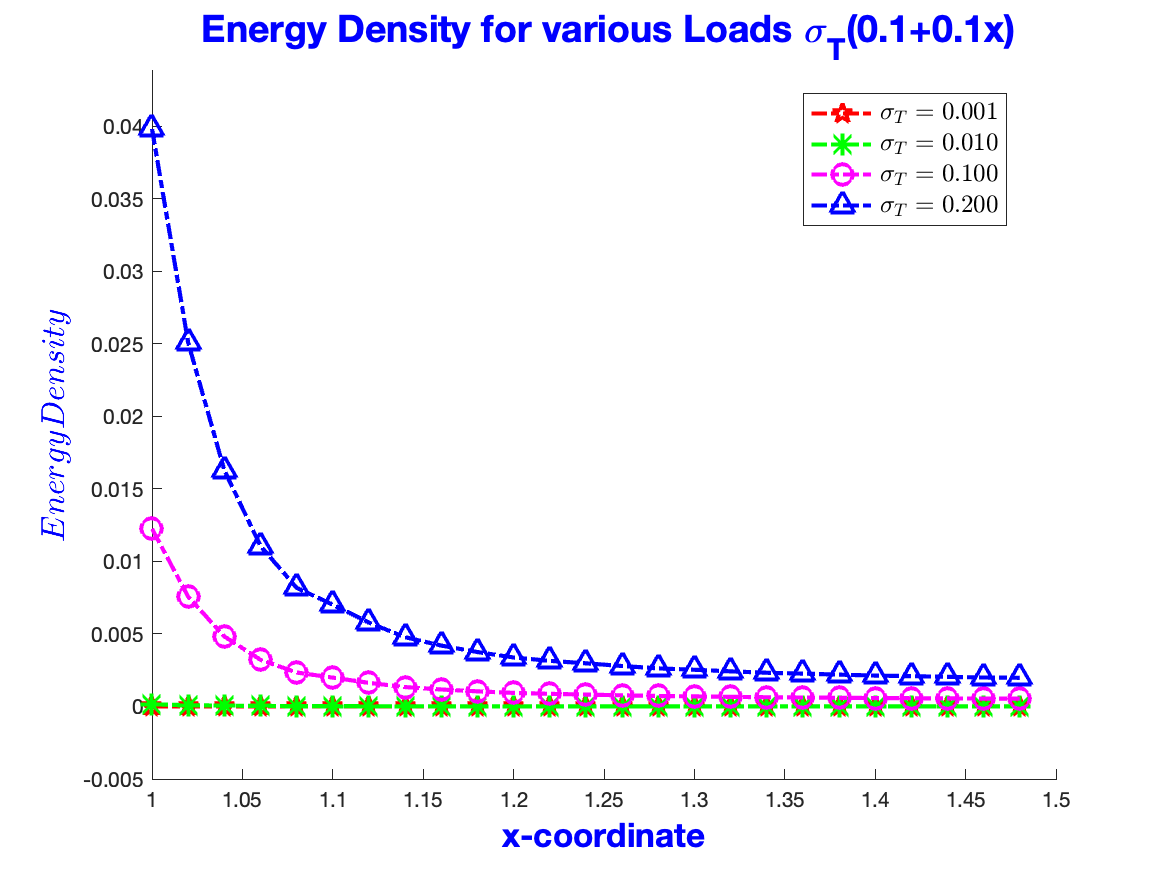}
        \caption{Energy density for various $\sigma_{T}$ for $\beta = 1.0$ and $\alpha = 1.0$}
        \label{fig:energy_density_sigma}
    \end{subfigure}
    \caption{Strain energy density plots for different parameter variations for the slope load in $x$-direction.}
    \label{fig:energy_density_1a}
\end{figure}

Figures~\ref{fig:stress_case1a}, \ref{fig:strain_case1a}, and \ref{fig:energy_density_1a} indicates that the stress, strain distributions and strain energy density ($\bfa{T} \colon \bfa{\epsilon}$) along a line extending to the crack tip for various values of the parameters~$\beta$, $\alpha$, and the top-load parameter $\sigma_{T}$.  A clear trend is observed in these plots: as the value of $\beta$ increases, the corresponding peak values of stress, strain, and strain energy density at the crack tip decrease moderately. This inverse relationship demonstrates that $\beta$ acts as a toughening or crack-mitigating parameter. Physically, this suggests that the mechanism governed by $\beta$—such as fiber bridging or localized plasticity—is effectively shielding the crack tip from the applied load. Consequently, increasing $\beta$ enhances the material's resilience to fracture by reducing the severity of the stress concentration at this critical point.  However, contrary results can be seen for the increasing values of $\alpha$ and $\sigma_T$, this result means the material's ability to resist crack propagation is decreasing as $\alpha$ increases. The higher stress and strain values indicate that a lower external force is required to reach the critical threshold for catastrophic failure. From a design and safety perspective, high values of $\alpha$ are undesirable and potentially dangerous.

\begin{figure}[H]
    \centering
    \begin{subfigure}{0.45\linewidth}
        \centering
        \includegraphics[width=\linewidth]{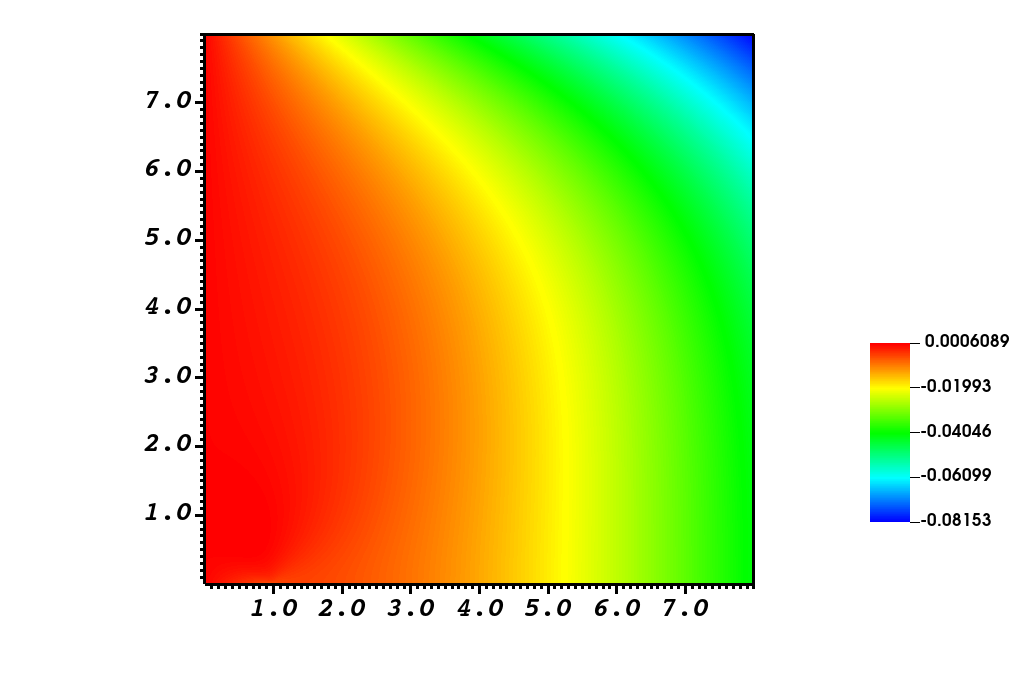}
        \caption{$x$-displacement for $\alpha = 1.0$, $\sigma_{T} = 0.1$, and $\beta = 1.0$}
        \label{fig:x_displacement}
    \end{subfigure}
    \hfill
    \begin{subfigure}{0.45\linewidth}
        \centering
        \includegraphics[width=\linewidth]{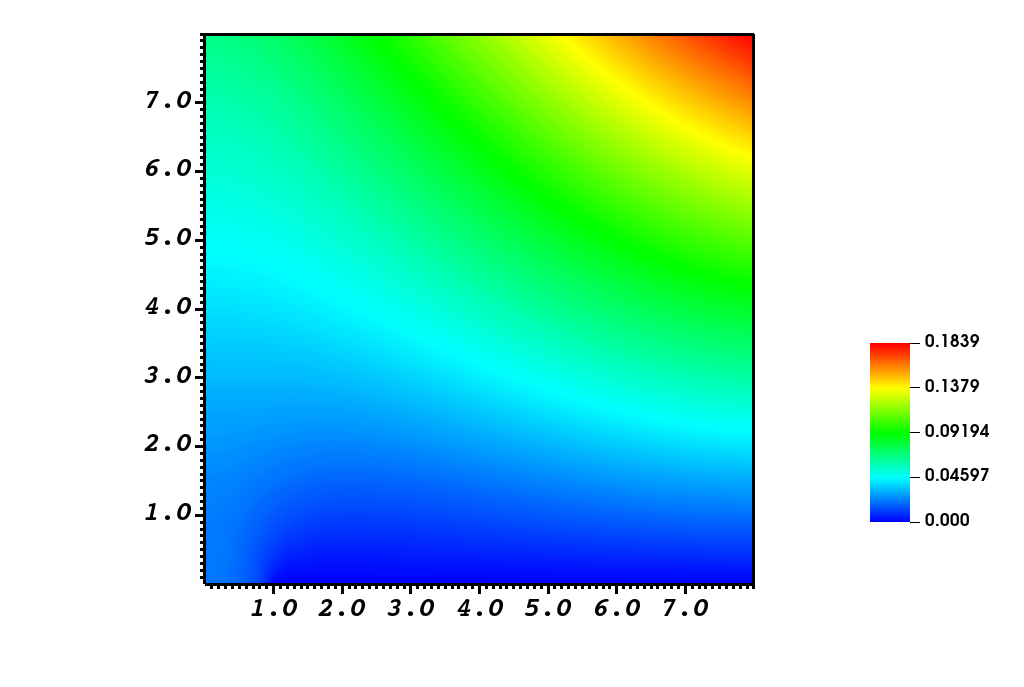}
        \caption{$y$-displacement for $\alpha = 1.0$, $\sigma_{T} = 0.1$, and $\beta = 1.0$}
        \label{fig:y_displacement}
    \end{subfigure}
    \vspace{1em}
    \begin{subfigure}{0.45\linewidth}
        \centering
        \includegraphics[width=\linewidth]{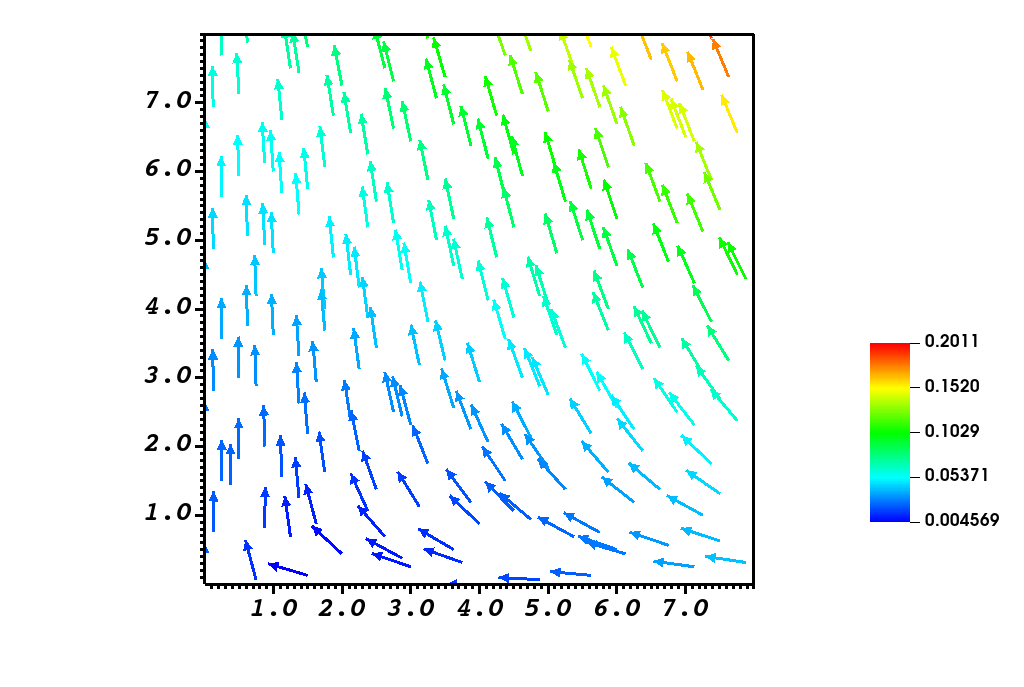}
        \caption{Vector-displacement for $\alpha = 1.0$, $\sigma_{T} = 0.1$, and $\beta = 1.0$}
        \label{fig:vector_displacement}
    \end{subfigure}
    \hfill
    \begin{subfigure}{0.45\linewidth}
        \centering
        \includegraphics[width=\linewidth]{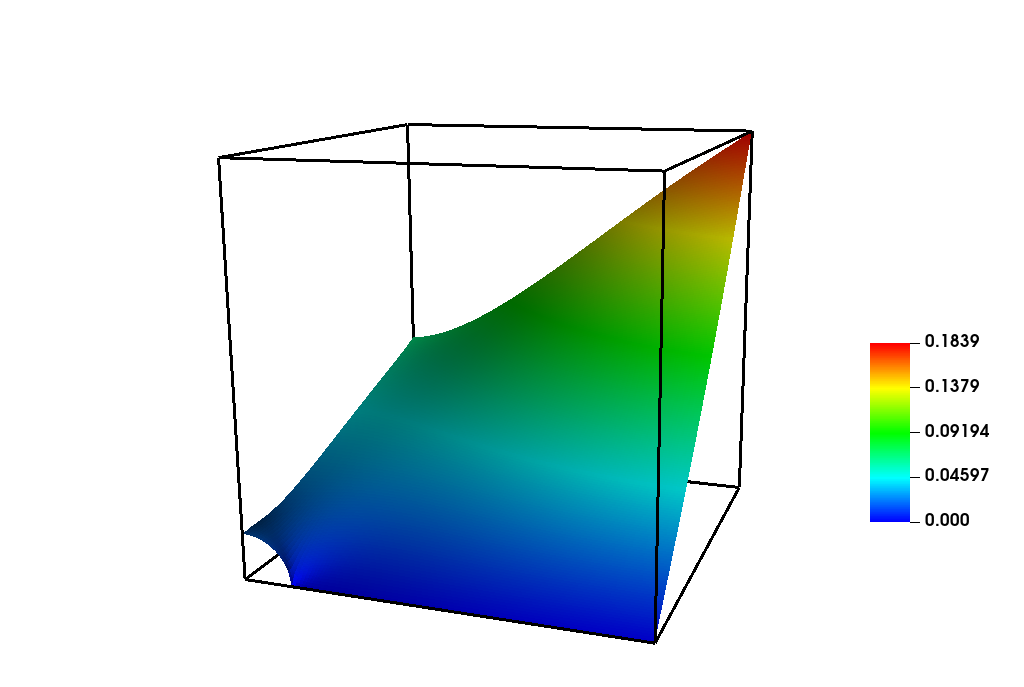}
        \caption{$y$-displacement at $\alpha = 1.0$, $\sigma_{T} = 0.1$, and $\beta = 1.0$}
        \label{fig:y_3Dcase1a}
    \end{subfigure}
    \caption{Displacement plots for $\alpha = 1.0$, $\sigma_{T} = 0.1$, and $\beta = 1.0$ for loads $\sigma_{T}(0.1+0.1x).$ in x-fiber direction.}
    \label{fig:displacement_combined_slope_M1}
\end{figure}

Figure~\ref{fig:displacement_combined_slope_M1} shows the displacements (both $x$ and $y$), both individually and as a vector. Figure~\ref{fig:y_3Dcase1a} illustrates the resulting crack opening profile. The shape is globally elliptical, consistent with classical fracture mechanics, but also exhibits significant blunting directly at the crack tip, a characteristic feature of plastic deformation.

\subsection{Case 1(b)-Slope load: Fibers aligned orthogonal to the crack}

In this part, we present the simulation results for an orthotropic material with an edge crack under mode-I type loading with the fibers aligned orthogonal to the crack plane. The numerical simulation of a Mode I crack problem where reinforcing fibers are aligned orthogonal to the crack plane is of significant scientific interest. This orientation represents the optimal configuration for maximizing fracture toughness, as the fibers act as the principal load-bearing constituents spanning the crack faces. This phenomenon, known as fiber bridging, is a primary toughening mechanism that shields the crack tip from the entire applied stress. Computational modeling is therefore essential for quantitatively evaluating the extent of this crack-shielding effect. Such simulations allow for the accurate prediction of the material's enhanced fracture resistance and fatigue life under tensile loading. Consequently, a thorough understanding of this configuration is crucial for the design and analysis of damage-tolerant composite structures in high-performance engineering applications.

In this numerical investigation, we analyze the critical case where the reinforcing fibers are oriented orthogonal to the crack plane. This configuration, representing the direction of maximum fracture toughness, is modeled by aligning the principal axis of material stiffness along the $y$-axis, while the crack resides along the $x$-axis. This physical setup is implemented computationally by defining the structural tensor as $\mathbf{M} = \mathbf{e}_2 \otimes \mathbf{e}_2$, where $\mathbf{e}_2$ is the unit vector in the $y$-direction. Such an alignment is of paramount practical importance as it directly resists Mode I crack opening through fiber bridging, fundamentally altering the stress distribution at the crack tip. For the numerical solution, the material's Poisson's ratio was held constant, and the resulting nonlinear system of equations was solved using the Picard iterative scheme. The solver demonstrated robust and efficient performance for this orientation, achieving rapid monotonic convergence as documented in Table~\ref{table2}. This result confirms the stability and effectiveness of our iterative method for this specific fiber-crack configuration. The results shown in the table are for $\alpha=1.0$, $\beta=1.0$, and $\sigma_T=0.1$. 

\begin{table}[H]
\centering
\begin{tabular}{|c|c|}
\hline 
Iteration No & Residual\tabularnewline
\hline 
\hline 
1 & 1.9708e-05\tabularnewline
\hline 
2 & 1.2046e-06\tabularnewline
\hline 
3 & 2.14193e-07\tabularnewline
\hline 
4 & 1.91804e-07\tabularnewline
\hline 
5 & 1.92496e-07\tabularnewline
\hline 
6 & 1.92453e-07\tabularnewline
\hline 
7 & 1.92456e-07\tabularnewline
\hline 
\end{tabular}
\caption{Value of the residual computed for each iteration for case-1(b). }
   \label{table2}
\end{table}

The convergence history reveals that the residual norm stagnated after the 7th iteration, remaining constant for all subsequent cycles up to the 100th iteration. This behavior indicates that while the current Picard iterative scheme is stable, it may not be sufficient to achieve a more stringent residual tolerance for the problem at hand, likely due to substantial material nonlinearities. To overcome this limitation and achieve a higher degree of accuracy, future work could explore more advanced numerical strategies.

\begin{figure}[H]
    \centering
    \begin{subfigure}{0.3\linewidth}
        \centering
        \includegraphics[width=\linewidth]{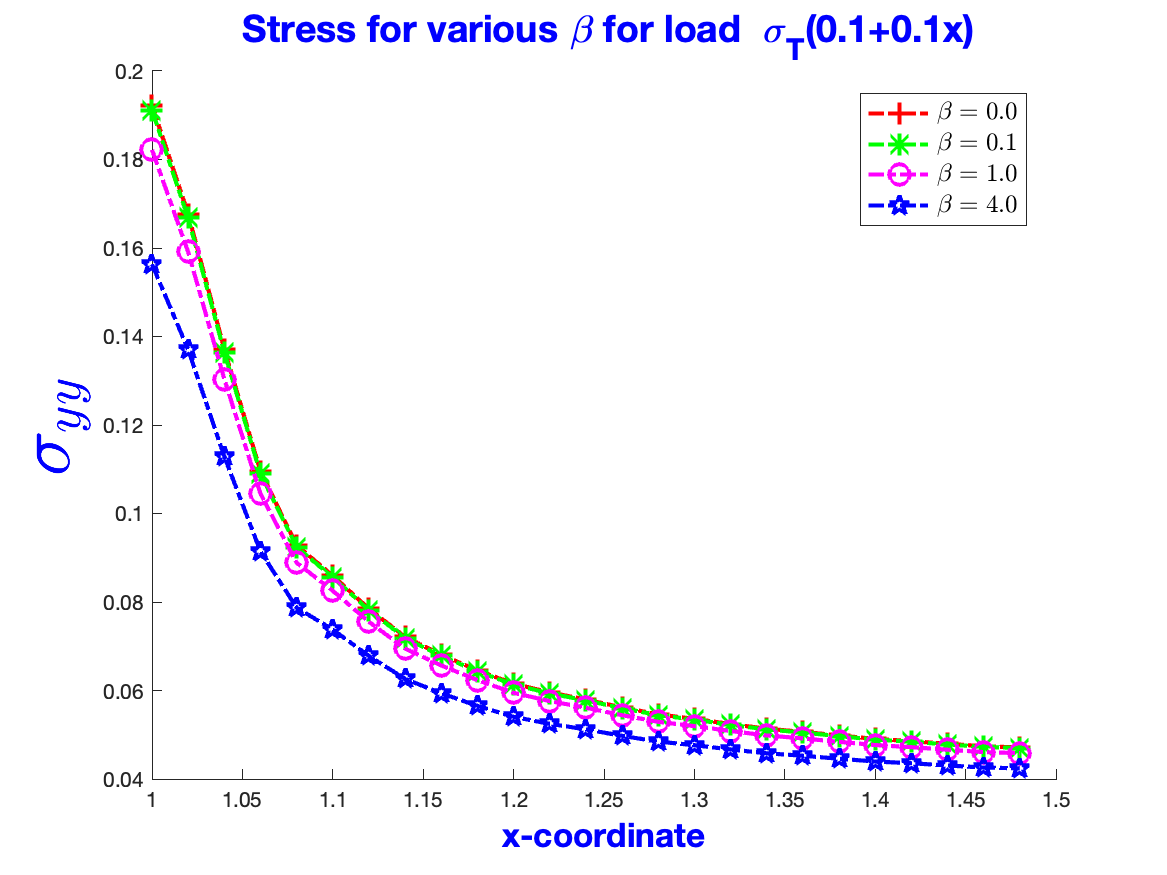}
        \caption{Stress for various $\beta$ for $\alpha = 1.0$ and $\sigma_{T} = 0.1$}
        \label{fig:stress_beta}
    \end{subfigure}
    \hfill
    \begin{subfigure}{0.3\linewidth}
        \centering
        \includegraphics[width=\linewidth]{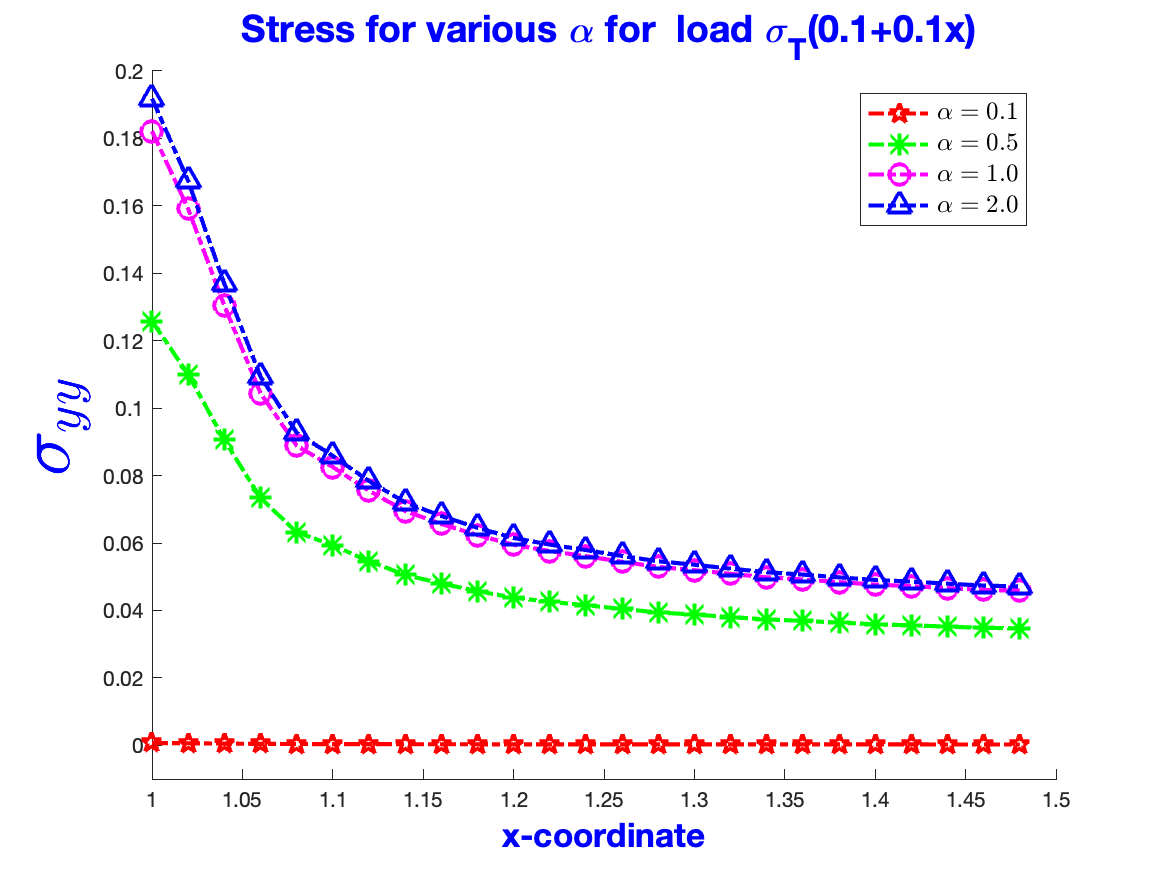}
        \caption{Stress for various $\alpha$ for $\beta = 1.0$ and $\sigma_{T} = 0.1$}
        \label{fig:stress_alpha}
    \end{subfigure}
    \hfill
    \begin{subfigure}{0.3\linewidth}
        \centering
        \includegraphics[width=\linewidth]{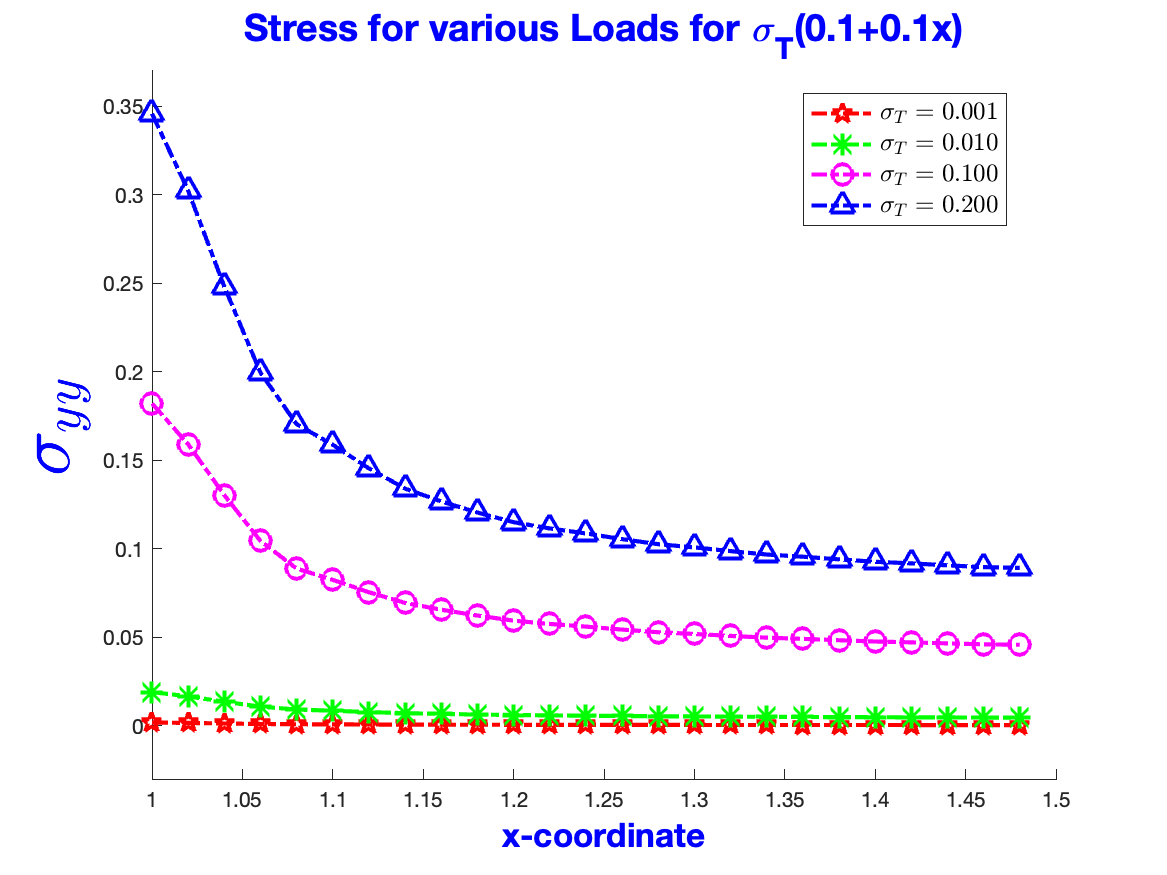}
        \caption{Stress for various $\sigma_{T}$ for $\beta = 1.0$ and $\alpha = 1.0$}
        \label{fig:stress_sigma}
    \end{subfigure}
    \caption{Stress distributions for various values of the parameters under a slope load applied in the $y$-direction.}
    \label{fig:stress_case1b}
\end{figure}

\begin{figure}[H]
    \centering
    \begin{subfigure}{0.3\linewidth}
        \centering
        \includegraphics[width=\linewidth]{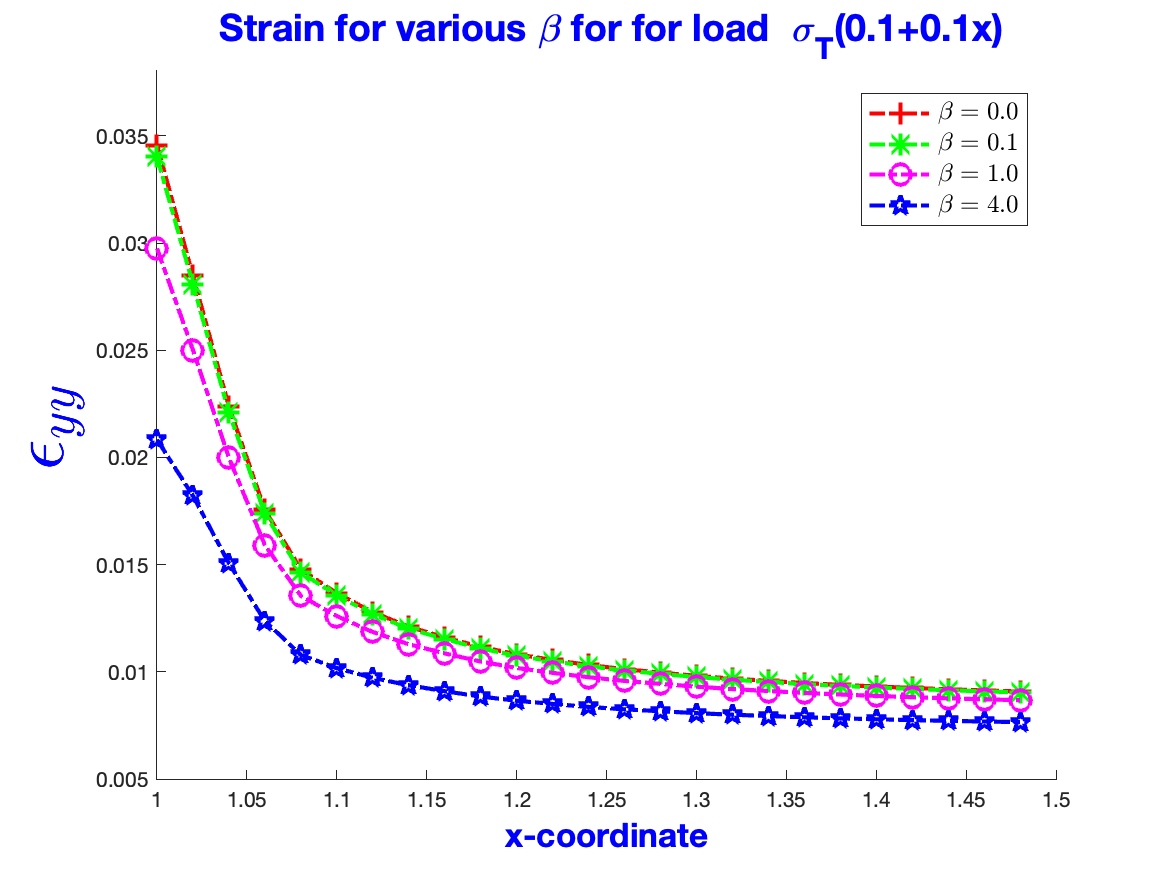}
        \caption{Strain for various $\beta$ for $\alpha = 1.0$ and $\sigma_{T} = 0.1$}
        \label{fig:strain_beta}
    \end{subfigure}
    \hfill
    \begin{subfigure}{0.3\linewidth}
        \centering
        \includegraphics[width=\linewidth]{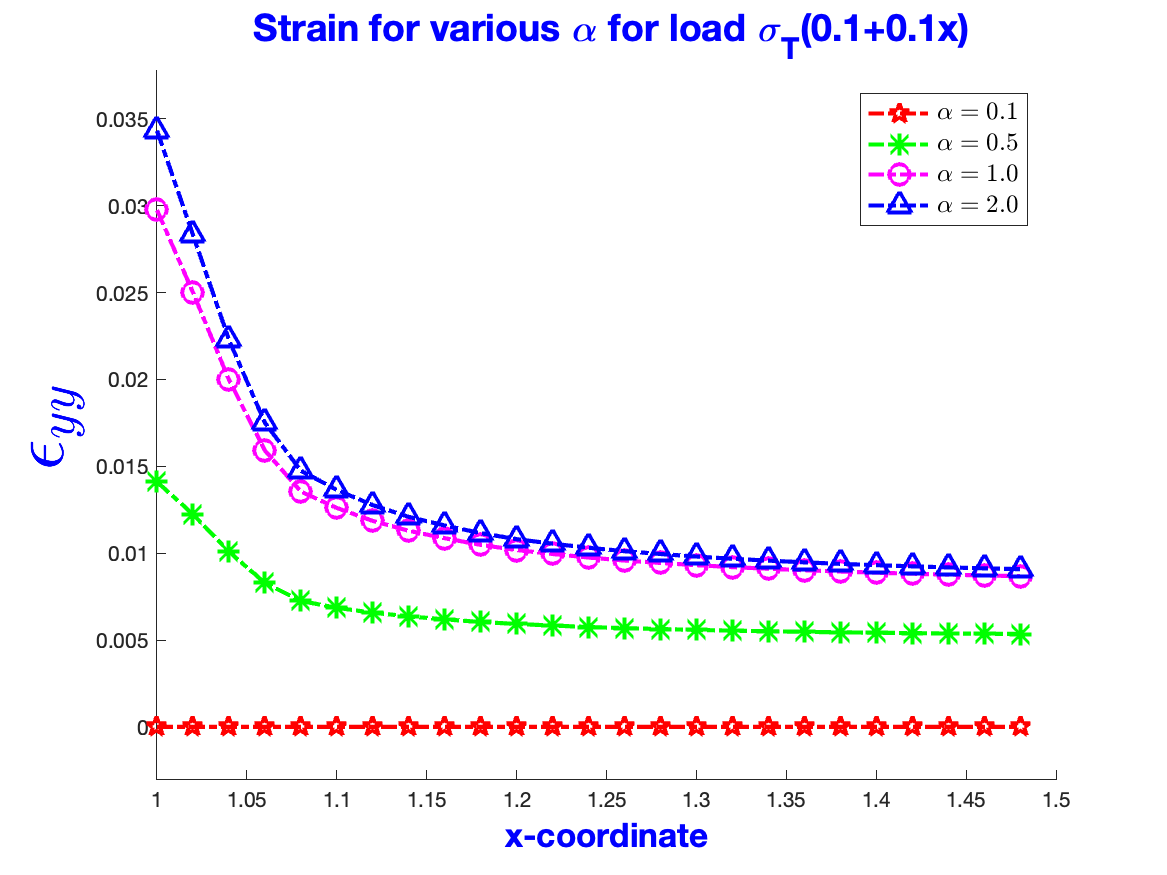}
        \caption{Strain for various $\alpha$ for $\beta = 1.0$ and $\sigma_{T} = 0.1$}
        \label{fig:strain_alpha}
    \end{subfigure}
    \hfill
    \begin{subfigure}{0.3\linewidth}
        \centering
        \includegraphics[width=\linewidth]{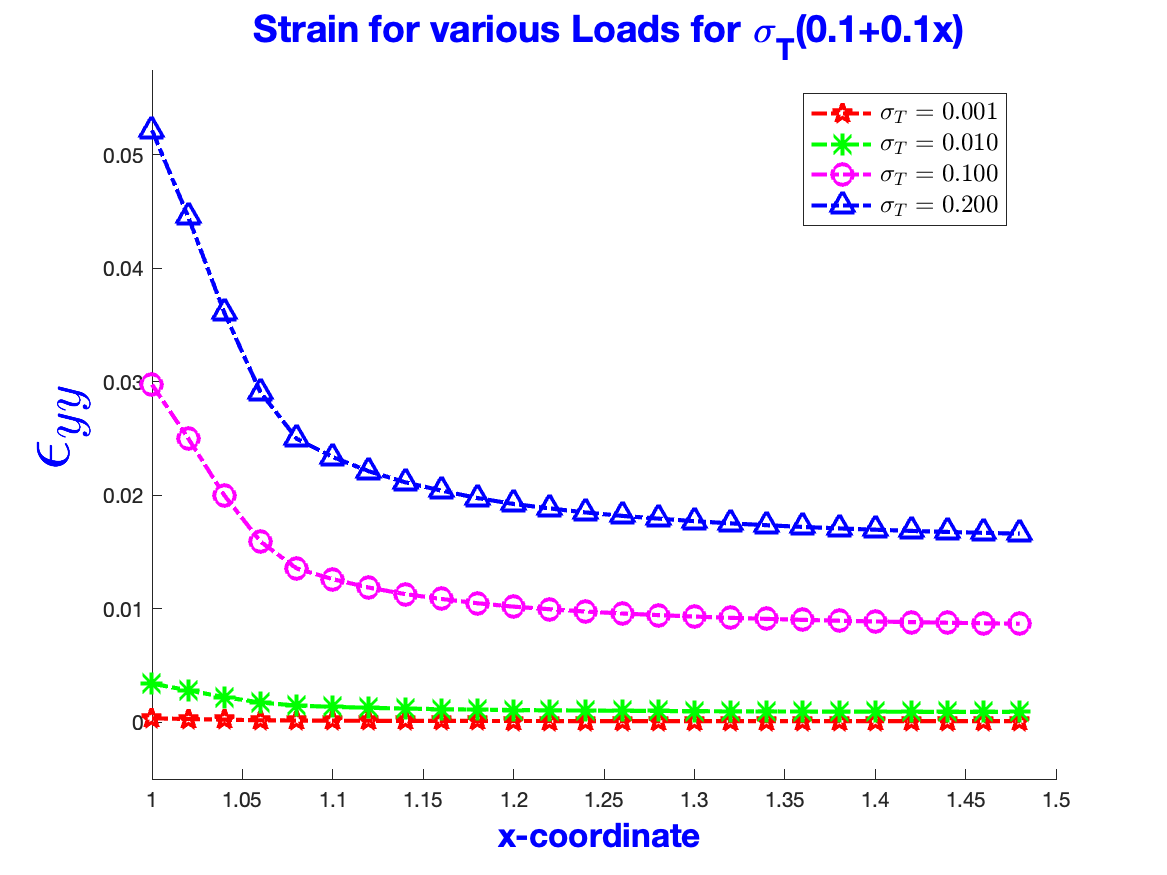}
        \caption{Strain for various $\sigma_{T}$ for $\beta = 1.0$ and $\alpha = 1.0$}
        \label{fig:strain_sigma}
    \end{subfigure}
    \caption{Strain plots for different parameter variations for slope loads y direction.}
    \label{fig:strain_case1b}
\end{figure}

\begin{figure}[H]
    \centering
    \begin{subfigure}{0.3\linewidth}
        \centering
        \includegraphics[width=\linewidth]{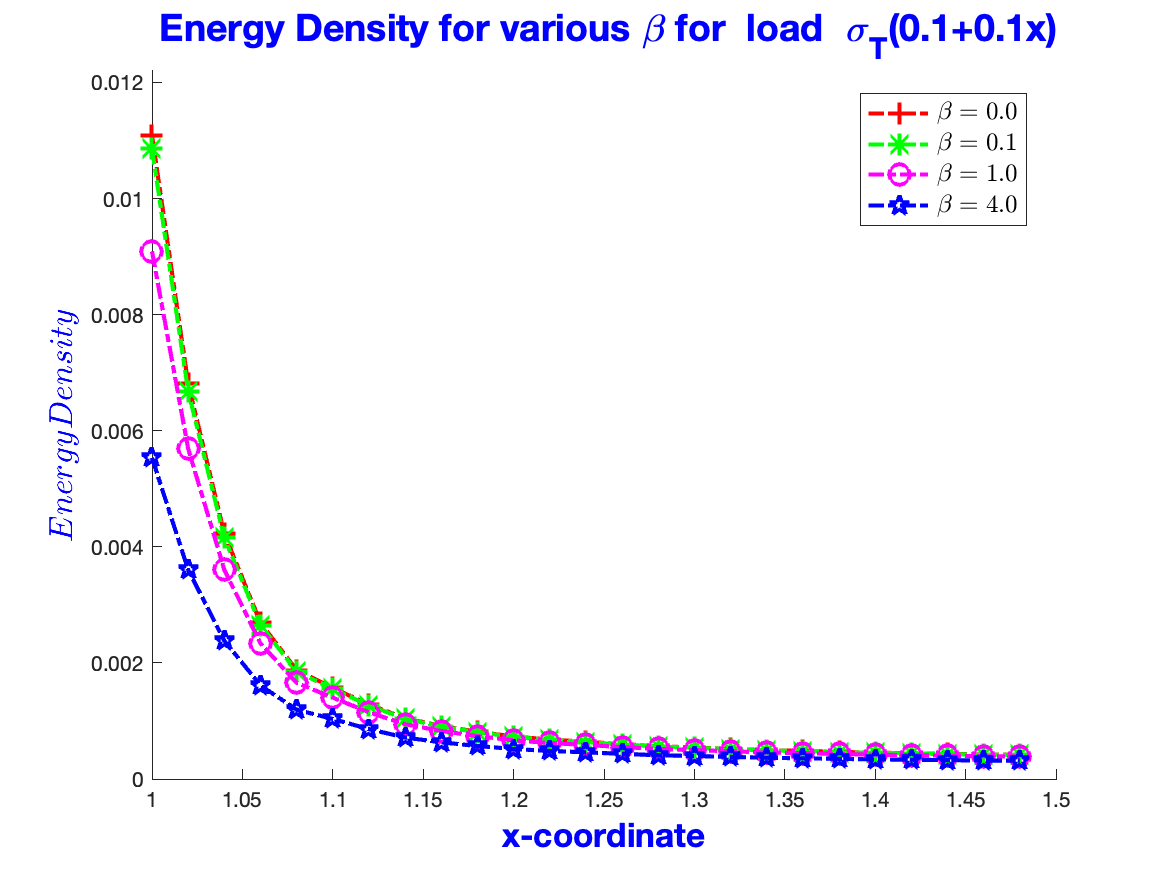}
        \caption{Energy density for various $\beta$ for $\alpha = 1.0$ and $\sigma_{T} = 0.1$}
        \label{fig:energy_density_beta}
    \end{subfigure}
    \hfill
    \begin{subfigure}{0.3\linewidth}
        \centering
        \includegraphics[width=\linewidth]{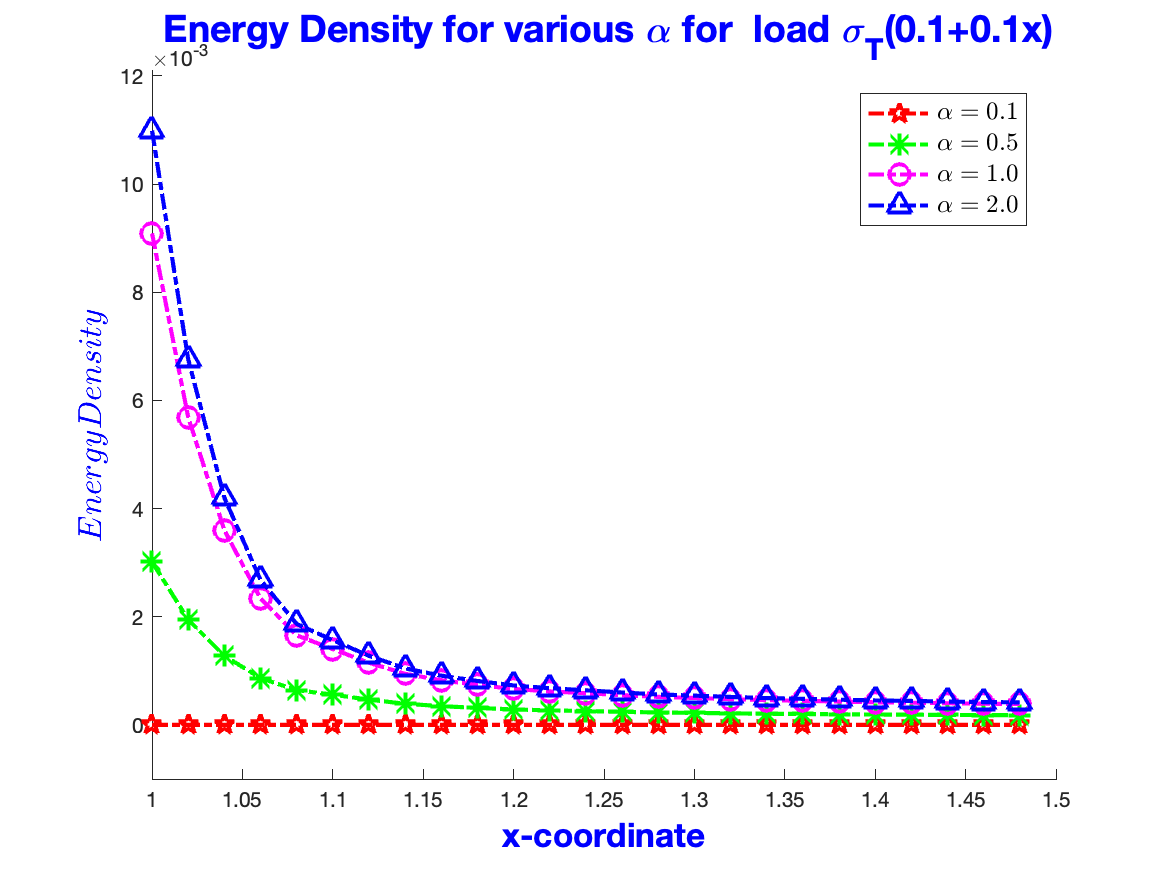}
        \caption{Energy density for various $\alpha$ for $\beta = 1.0$ and $\sigma_{T} = 0.1$}
        \label{fig:energy_density_alpha}
    \end{subfigure}
    \hfill
    \begin{subfigure}{0.3\linewidth}
        \centering
        \includegraphics[width=\linewidth]{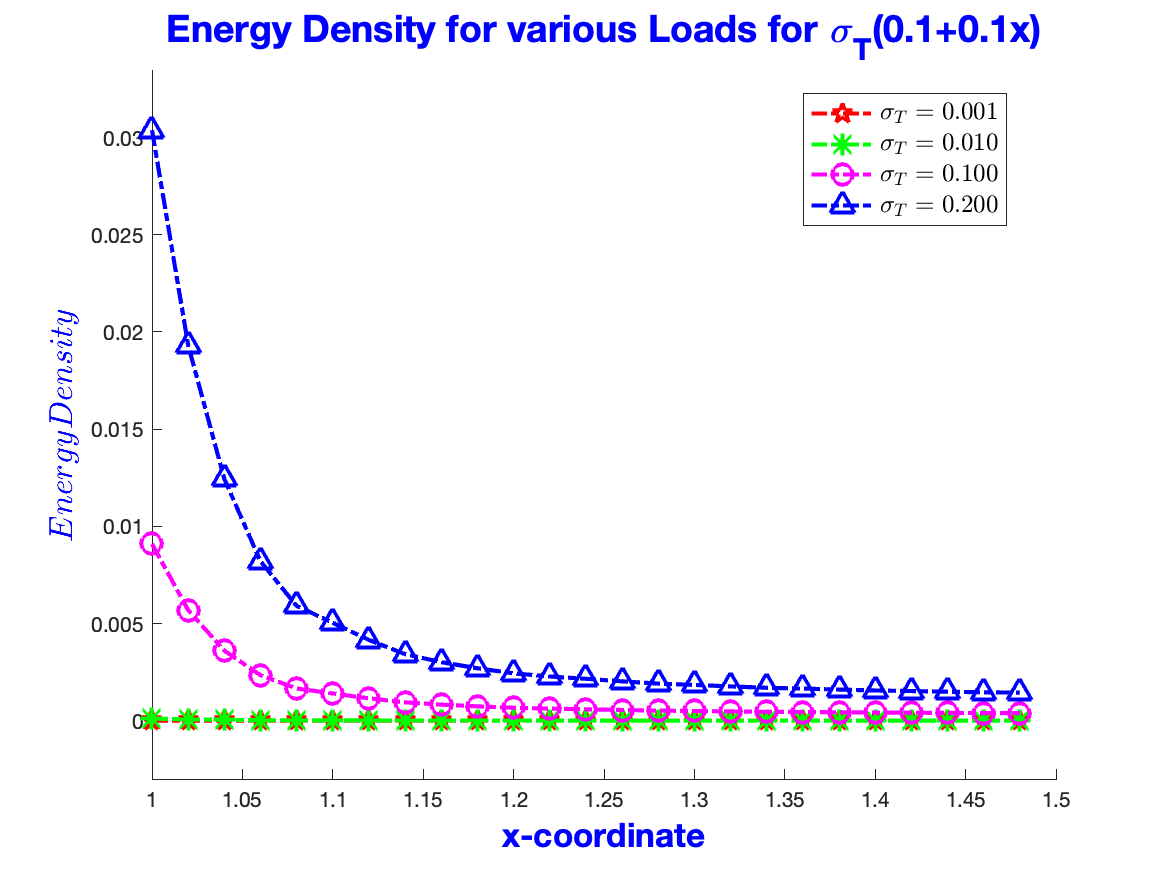}
        \caption{Energy density for various $\sigma_{T}$ for $\beta = 1.0$ and $\alpha = 1.0$ }
        \label{fig:energy_density_sigma}
    \end{subfigure}
    \caption{Energy density plots for different parameter variations for slope loads y direction.}
    \label{fig:energy_density_case1b}
\end{figure}

Figures~\ref{fig:stress_case1b}, \ref{fig:strain_case1b}, and \ref{fig:energy_density_case1b} depict the influence of the key parameters on the stress, strain, and strain energy density along a line directly ahead of the crack-tip.  The subsequent analysis examines the impact of key model parameters on the mechanical state at the crack tip. Specifically, we investigate the distributions of the stress tensor ($\bfa{T}$), strain tensor ($\bfa{\epsilon}$), and the strain energy density ($\bfa{T} \colon \bfa{\epsilon}$) along a line approaching the crack front. The study focuses on the effects of a toughening parameter ($\beta$), a potential degradation parameter ($\alpha$), and the magnitude of the applied top load ($\sigma_T$).  The results demonstrate a clear inverse relationship between the parameter $\beta$ and the severity of the crack-tip fields. As $\beta$ is increased, the peak values of stress, strain, and strain energy density exhibit a moderate decrease. This confirms that $\beta$ functions as a crack-mitigating parameter, wherein its underlying physical mechanism—such as fiber bridging or localized plasticity—acts to shield the crack tip. This shielding effectively dissipates energy and reduces the local stress concentration, thereby enhancing the material's overall fracture resilience. In stark contrast, increasing the values of the parameter $\alpha$ and the applied load $\sigma_T$ leads to an amplification of the crack-tip fields. The intensification due to $\sigma_T$ is an expected outcome of increasing the external load. However, the similar effect of $\alpha$ is more revealing: it signifies that this parameter governs a mechanism that is detrimental to the material's integrity. An increasing $\alpha$ heightens the stress concentration, indicating a reduction in the material's ability to resist crack propagation. From a design standpoint, this identifies $\alpha$ as a critical parameter whose value must be carefully controlled to ensure structural safety and avoid premature failure.

\begin{figure}[H]
    \centering
    \begin{subfigure}{0.45\linewidth}
        \centering
        \includegraphics[width=\linewidth]{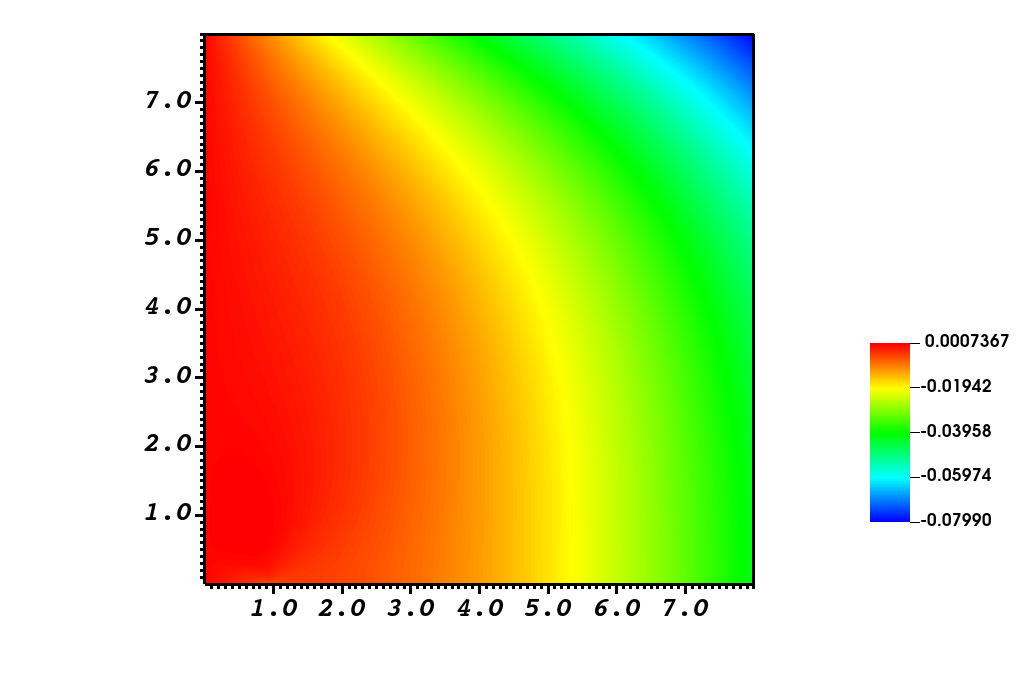}
        \caption{X-displacement at $\alpha = 1.0$, $\sigma_{T} = 0.1$, and $\beta = 1.0$}
        \label{fig:x_displacement}
    \end{subfigure}
    \hfill
    \begin{subfigure}{0.45\linewidth}
        \centering
        \includegraphics[width=\linewidth]{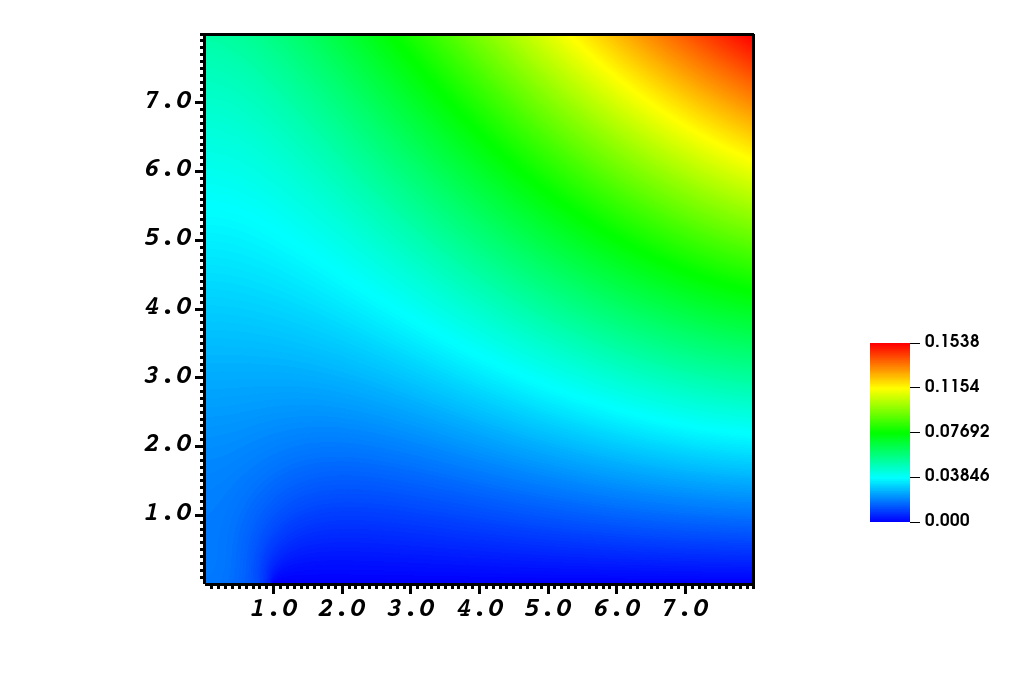}
        \caption{Y-displacement at $\alpha = 1.0$, $\sigma_{T} = 0.1$, and $\beta = 1.0$}
        \label{fig:y_displacement}
    \end{subfigure}
    
    \vspace{1em}
    
    \begin{subfigure}{0.45\linewidth}
        \centering
        \includegraphics[width=\linewidth]{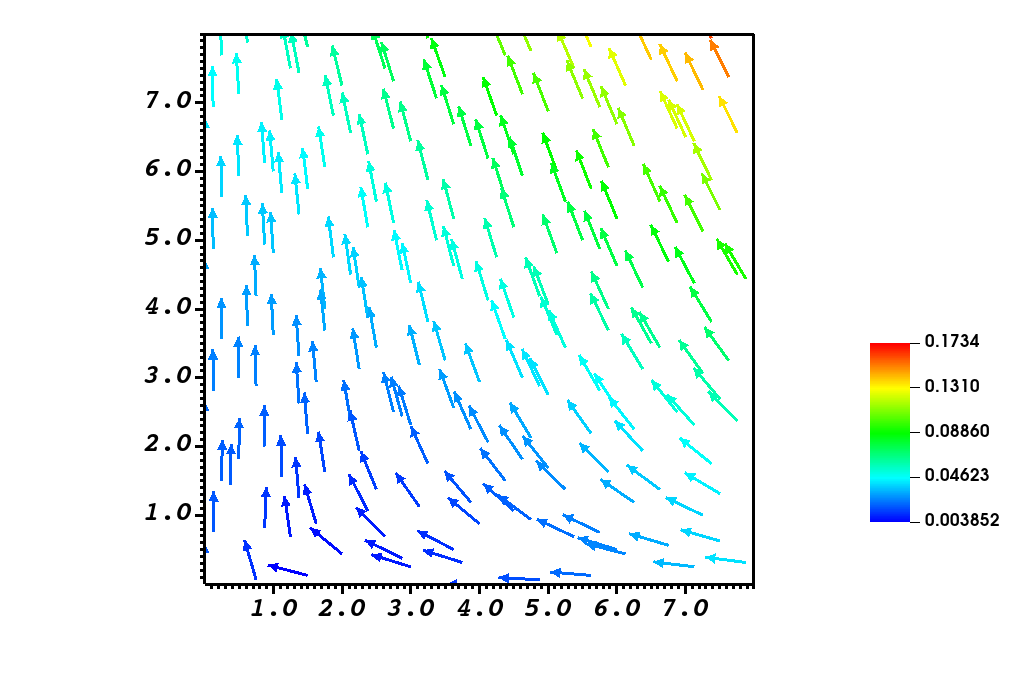}
        \caption{Vector-displacement at $\alpha = 1.0$, $\sigma_{T} = 0.1$, and $\beta = 1.0$}
        \label{fig:vector_displacement}
    \end{subfigure}
    \hfill
    \begin{subfigure}{0.45\linewidth}
        \centering
        \includegraphics[width=\linewidth]{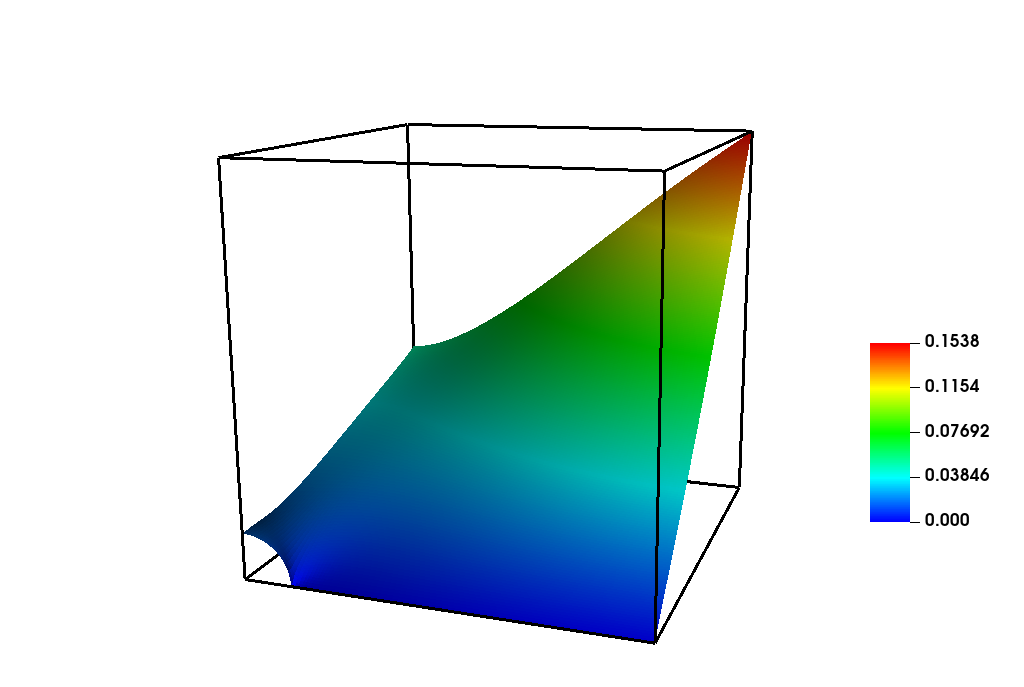}
        \caption{Y-displacement at $\alpha = 1.0$, $\sigma_{T} = 0.1$, and $\beta = 1.0$ (3D view)}
        \label{fig:y_displacement_3D_case1b}
    \end{subfigure}

\caption{Displacement plot for the y-direction fiber orientation under a linear slope load, using the parameter set: $\alpha = 1.0$, $\beta = 1.0$, and $\sigma_{T} = 0.1$.}
    \label{fig:displacement_3D_case1b}
\end{figure}

Figure~\ref{fig:displacement_3D_case1b} shows the displacements (both $x$ and $y$), both individually and as a vector. Figure~\ref{fig:y_displacement_3D_case1b} illustrates the resulting crack opening profile. The crack opening exhibits two distinct features: a globally elliptical profile consistent with linear elastic theory, and a pronounced blunting at the crack tip characteristic of plastic yielding.

\subsection{Case 2(a)-Non-uniform load: Fibers aligned with the crack  plane}
This numerical investigation examines a cracked material reinforced with fibers aligned with the $x$-axis (parallel to the crack) that is subjected to a non-uniform load on its top surface.
Non-uniform loads, such as a sinusoidal distribution represented by $\sigma_T \left( \dfrac{\sin (\pi \, x)}{8}  \right)$, are crucial for studying crack tip fields because they more accurately model realistic and complex loading scenarios compared to simple uniform tension. In practice, loads on structures are rarely uniform due to factors such as thermal gradients, residual stresses from processes like welding, or complex pressure distributions on components, such as aerospace fuselages or pressure vessels. Analyzing a crack's response to such a distribution provides a more precise understanding of stress intensity factors and the conditions for crack propagation. This enables engineers to more accurately predict the fatigue life and fracture toughness of materials in real-world applications, resulting in safer and more reliable designs.

The nonlinear system of equations was solved numerically using the combination of Picard's iterative scheme and finite element method, assuming a constant Poisson's ratio for the material. Our solver proved to be highly stable and efficient for this fiber-crack configuration, achieving rapid monotonic convergence as detailed in Table~\ref{table3}. The results presented in the table, which confirm the method's effectiveness, correspond to parameter values of $\alpha=1.0$, $\beta=1.0$, and $\sigma_T=0.1$.

\begin{table}[H]
\centering
\begin{tabular}{|c|c|}
\hline 
Iteration No & Residual\tabularnewline
\hline 
\hline 
1 & 0.000181143\tabularnewline
\hline 
2 & 2.33081e-05\tabularnewline
\hline 
3 & 2.94533e-06\tabularnewline
\hline 
4 & 6.15349e-07\tabularnewline
\hline 
5 & 6.04267e-07\tabularnewline
\hline 
6 & 5.94436e-07\tabularnewline
\hline 
7 & 5.95461e-07\tabularnewline
\hline 
8 & 5.95343e-07\tabularnewline
\hline 
9 & 5.95356e-07\tabularnewline
\hline 
10 & 5.95355e-07\tabularnewline
\hline 
\end{tabular}
\caption{Iterative solver performance for non-uniform load with fibers parallel to the crack.}
\label{table3}
\end{table}

\begin{figure}[H]
    \centering
    \begin{subfigure}{0.3\linewidth}
        \centering
        \includegraphics[width=\linewidth]{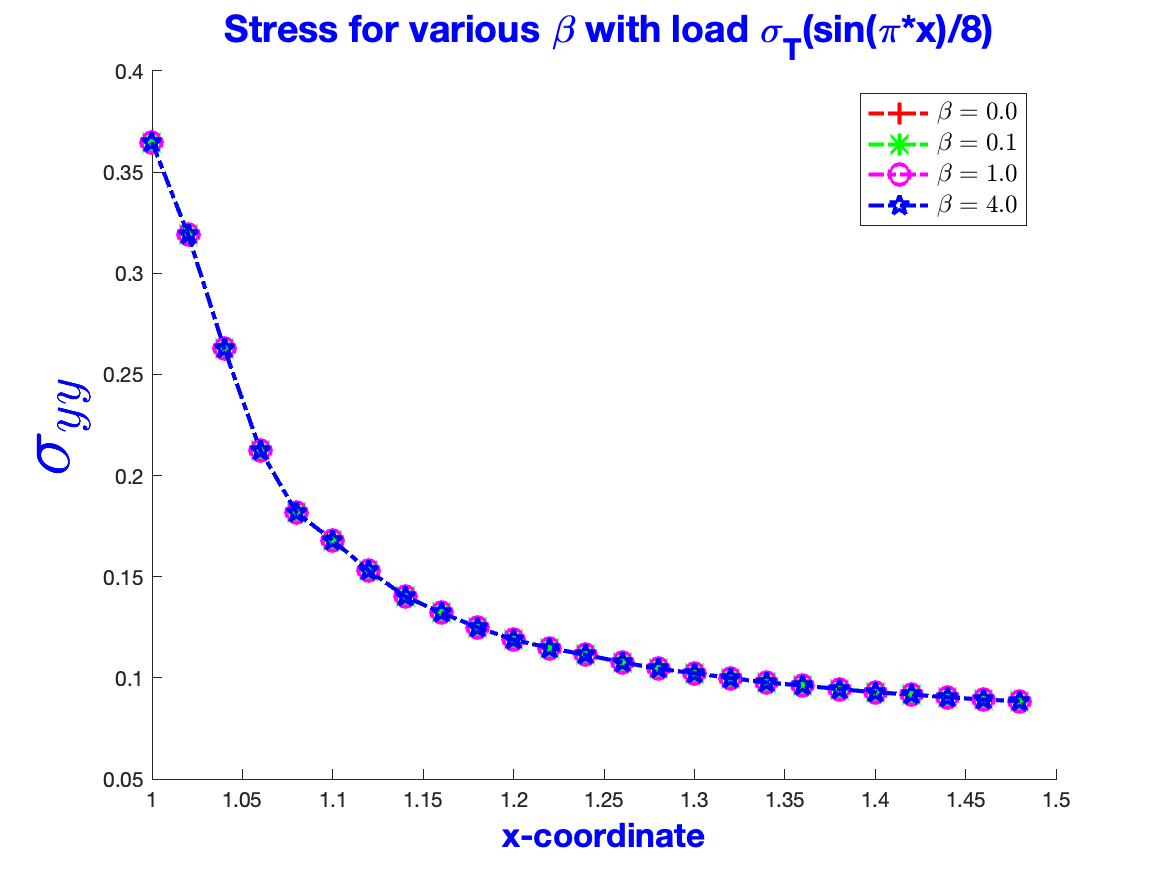}
        \caption{Stress for various $\beta$ for $\alpha = 1.0$ and $\sigma_{T} = 0.1$}
        \label{fig:strain_beta}
    \end{subfigure}
    \hfill
    \begin{subfigure}{0.3\linewidth}
        \centering
        \includegraphics[width=\linewidth]{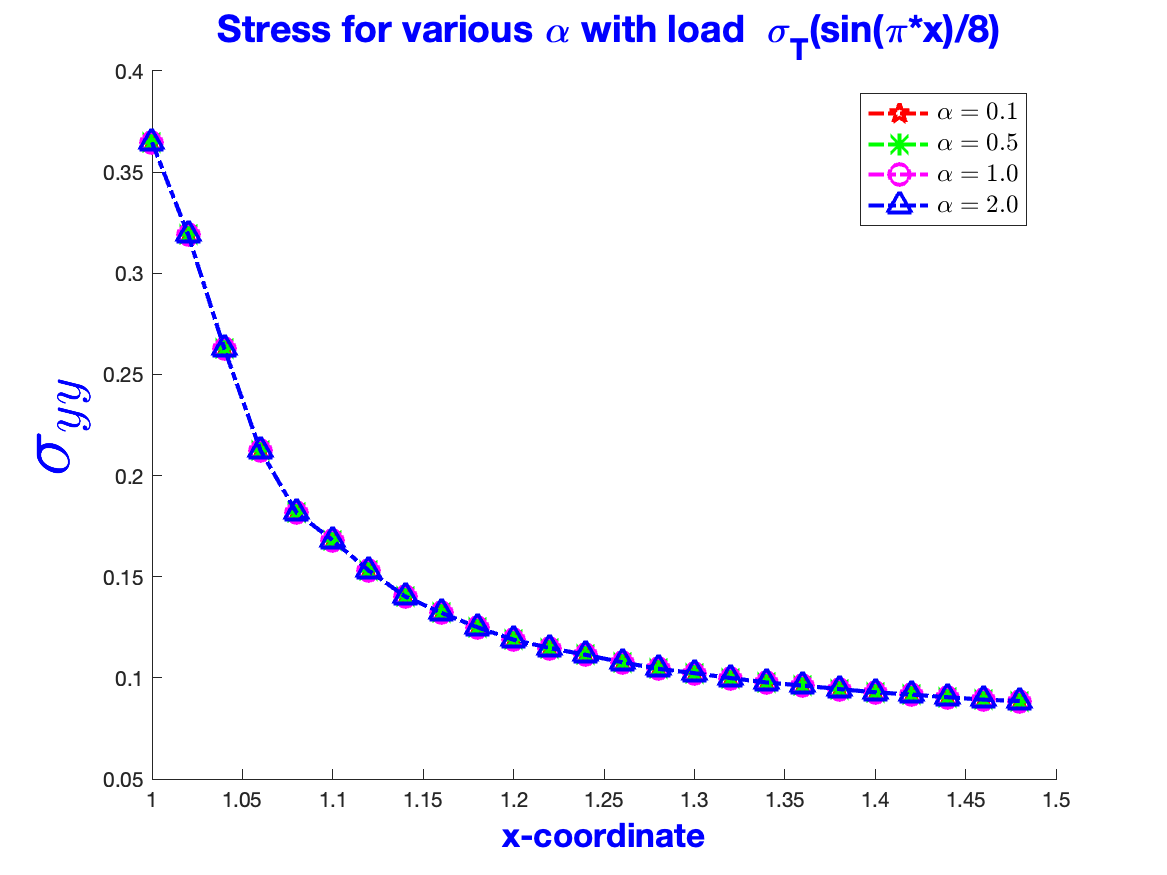}
        \caption{Stress for various $\alpha$ for $\beta = 1.0$ and $\sigma_{T} = 0.1$}
        \label{fig:strain_alpha}
    \end{subfigure}
    \hfill
    \begin{subfigure}{0.3\linewidth}
        \centering
        \includegraphics[width=\linewidth]{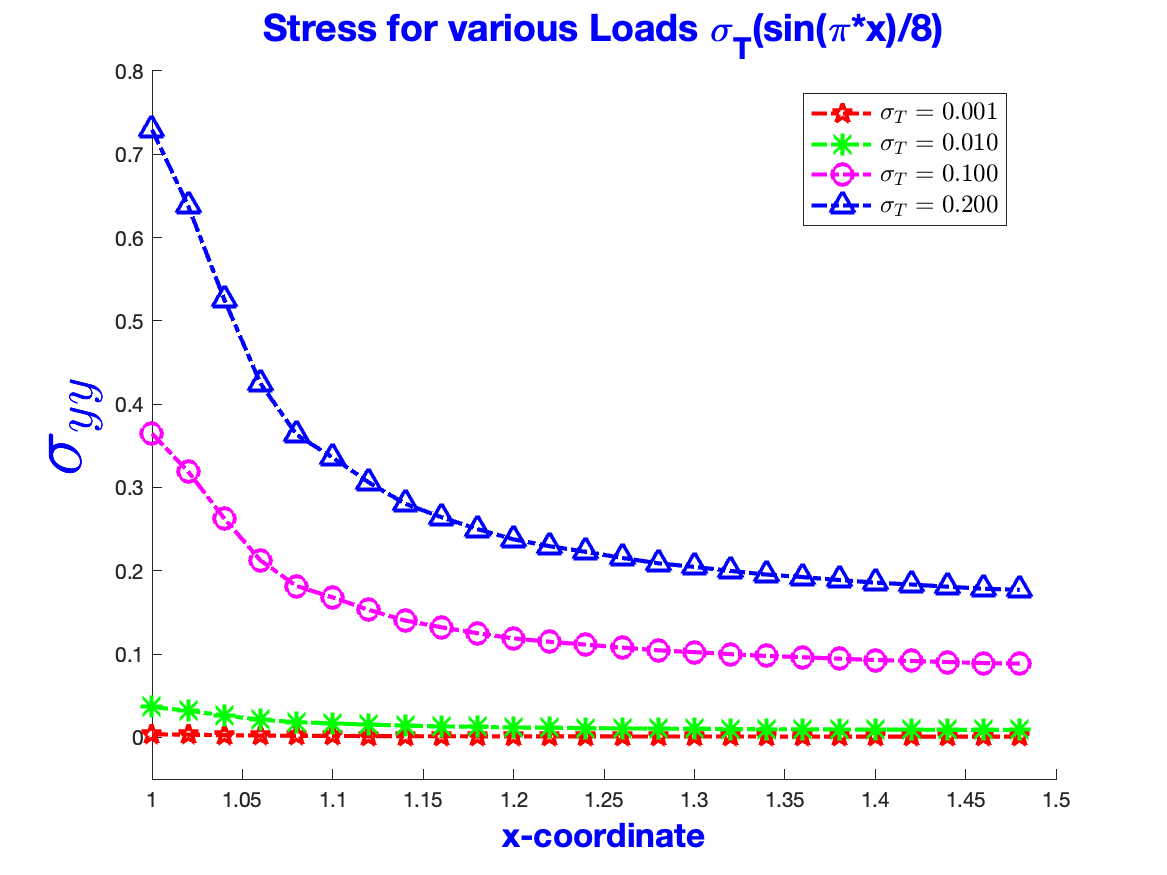}
        \caption{Stress for various $\sigma_{T}$ for $\beta = 1.0$ and $\alpha = 1.0$}
        \label{fig:strain_sigma}
    \end{subfigure}
    \caption{Stress plots for different parameter variations for sine loads.}
    \label{fig:stress_2a}
\end{figure}

\begin{figure}[H]
    \centering
    \begin{subfigure}{0.3\linewidth}
        \centering
        \includegraphics[width=\linewidth]{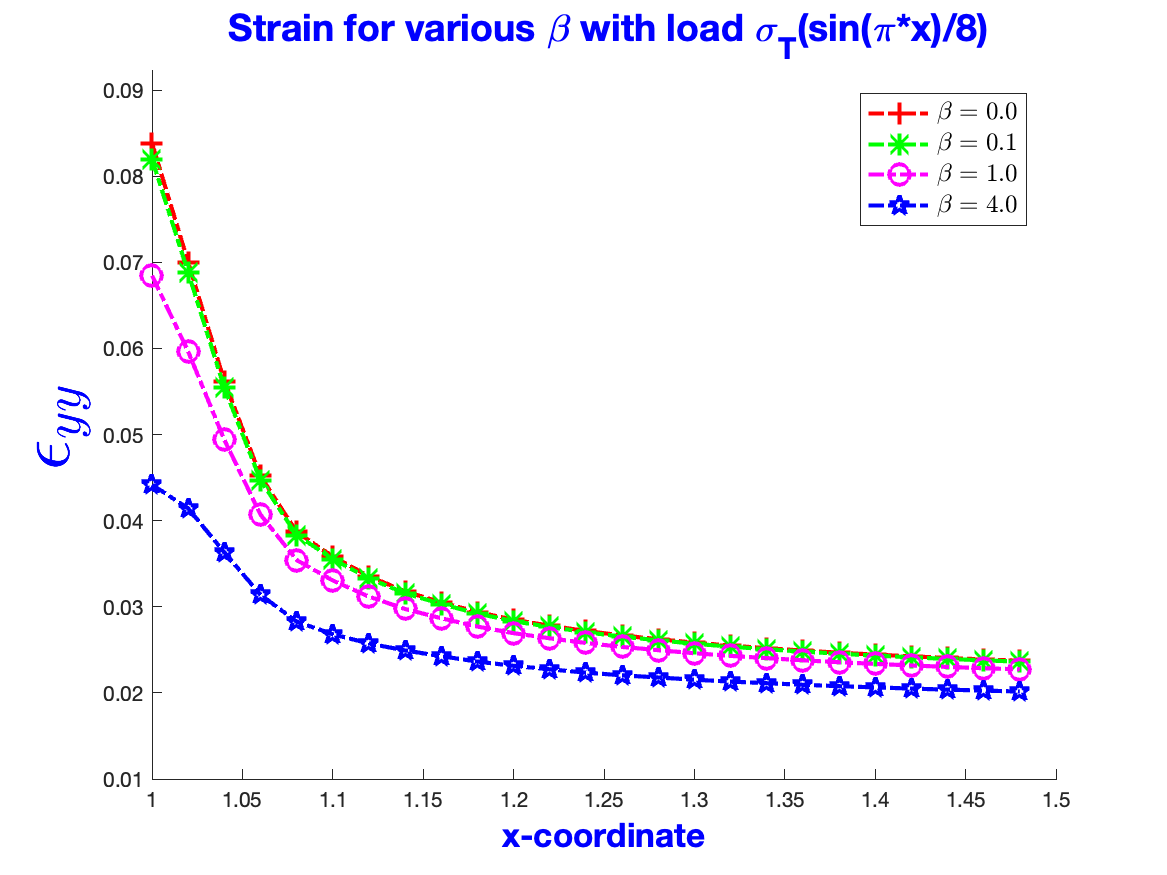}
        \caption{Strain for various $\beta$ for $\alpha = 1.0$ and $\sigma_{T} = 0.1$}
        \label{fig:strain_beta}
    \end{subfigure}
    \hfill
    \begin{subfigure}{0.3\linewidth}
        \centering
        \includegraphics[width=\linewidth]{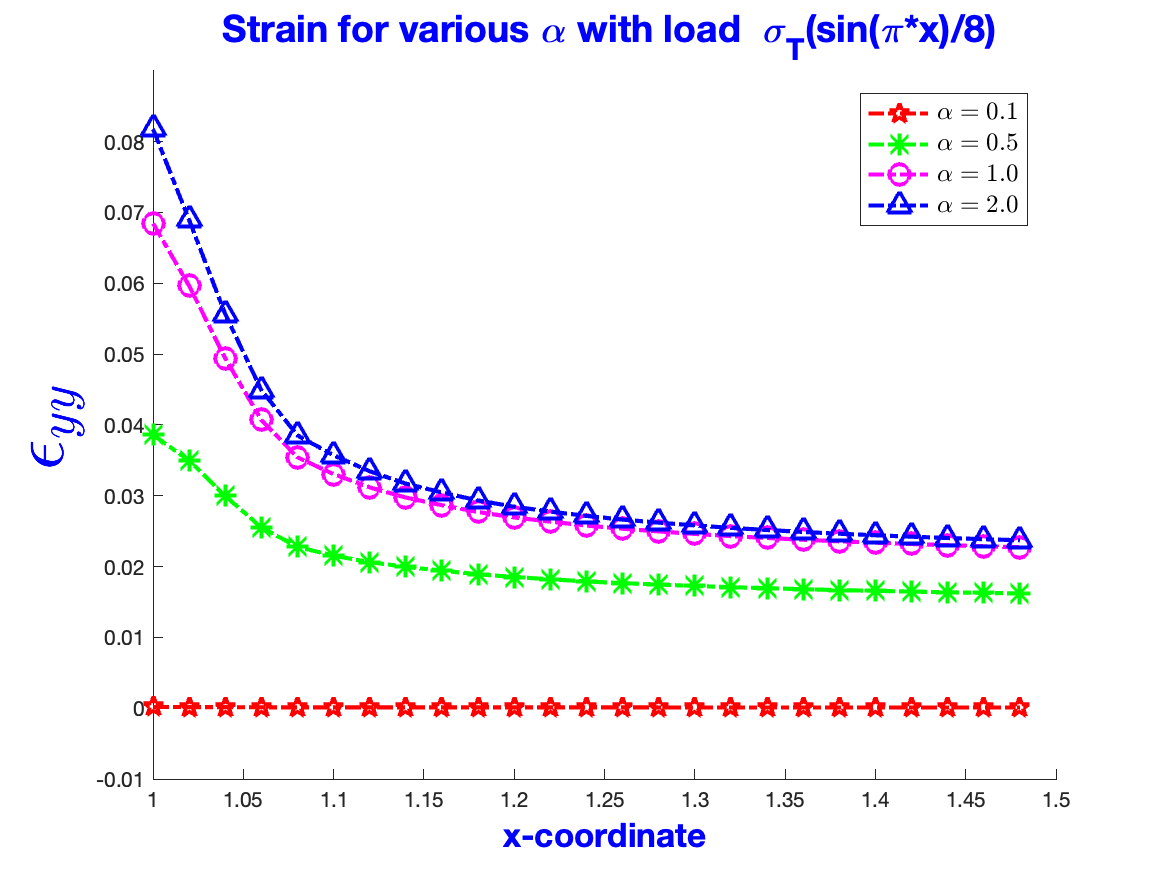}
        \caption{Strain for various $\alpha$ for $\beta = 1.0$ and $\sigma_{T} = 0.1$}
        \label{fig:strain_alpha}
    \end{subfigure}
    \hfill
    \begin{subfigure}{0.3\linewidth}
        \centering
        \includegraphics[width=\linewidth]{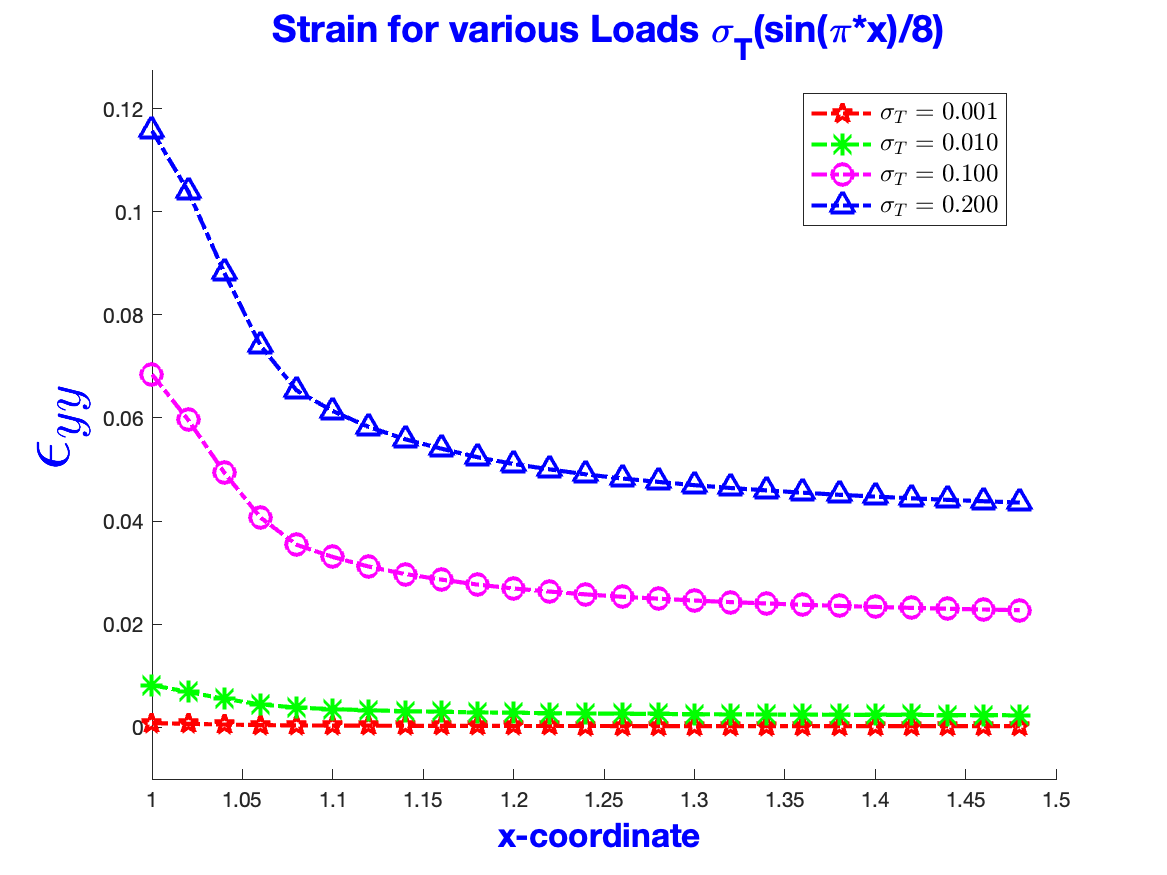}
        \caption{Strain for various $\sigma_{T}$ for $\beta = 1.0$ and $\alpha = 1.0$}
        \label{fig:strain_sigma}
    \end{subfigure}
    \caption{Strain plots for different parameter variations for sine load.}
    \label{fig:strain_2a}
\end{figure}

\begin{figure}[H]
    \centering
    \begin{subfigure}{0.3\linewidth}
        \centering
        \includegraphics[width=\linewidth]{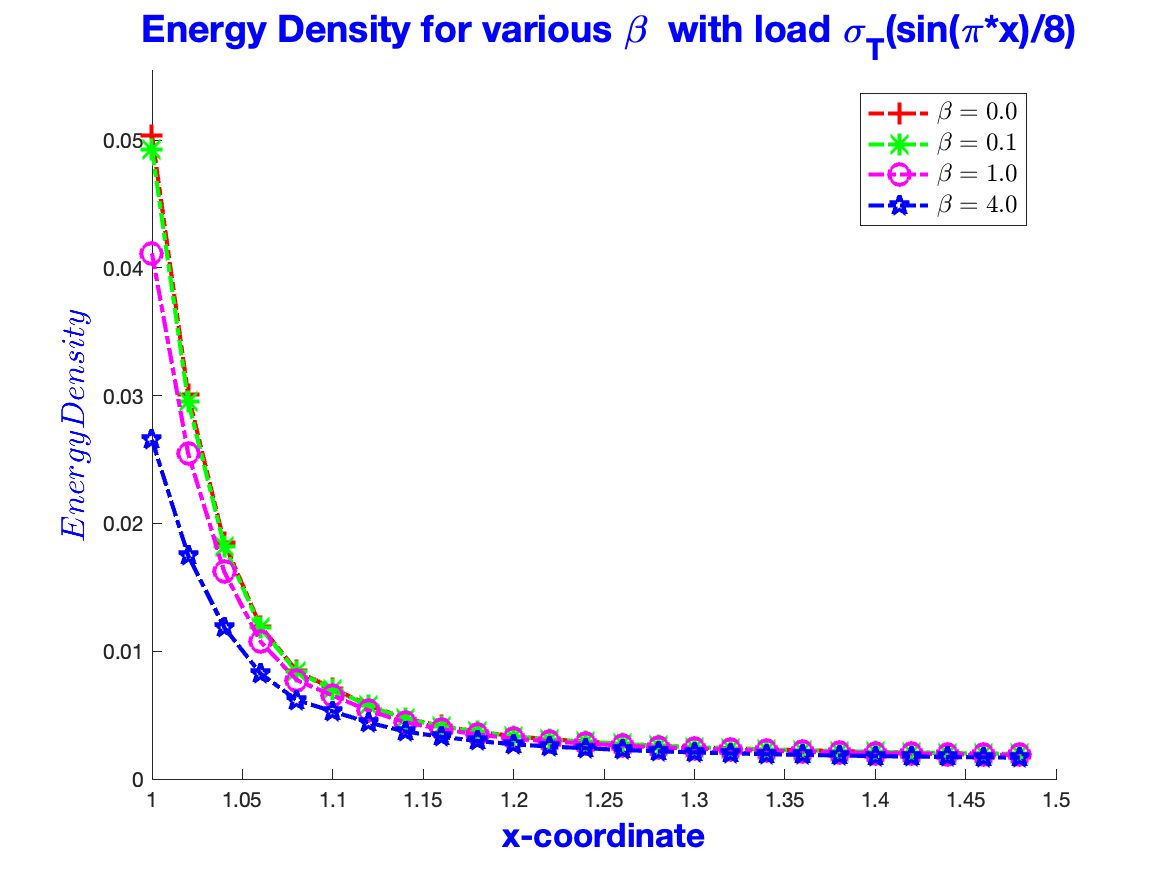}
        \caption{Energy density for various $\beta$ for $\alpha = 1.0$ and $\sigma_{T} = 0.1$}
        \label{fig:energy_density_beta}
    \end{subfigure}
    \hfill
    \begin{subfigure}{0.3\linewidth}
        \centering
        \includegraphics[width=\linewidth]{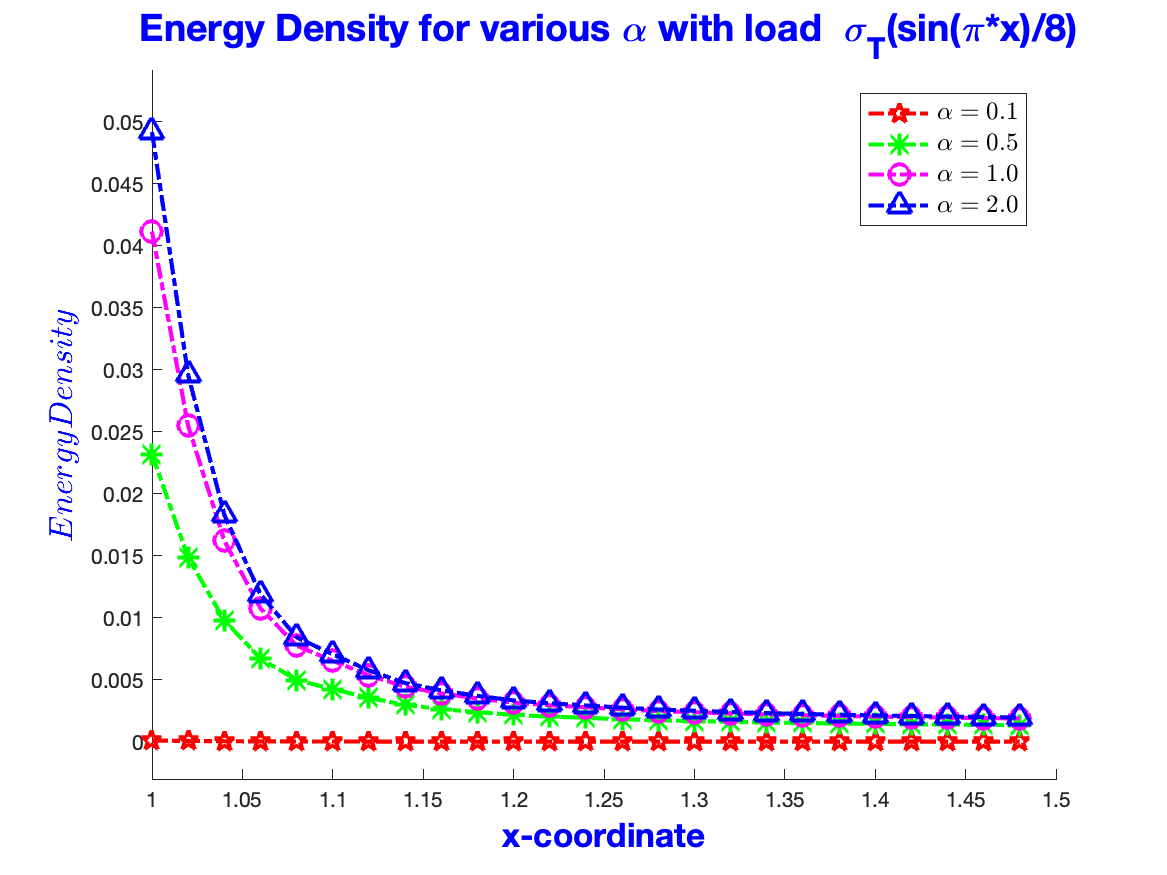}
        \caption{Energy density for various $\alpha$ for $\beta = 1.0$ and $\sigma_{T} = 0.1$}
        \label{fig:energy_density_alpha}
    \end{subfigure}
    \hfill
    \begin{subfigure}{0.3\linewidth}
        \centering
        \includegraphics[width=\linewidth]{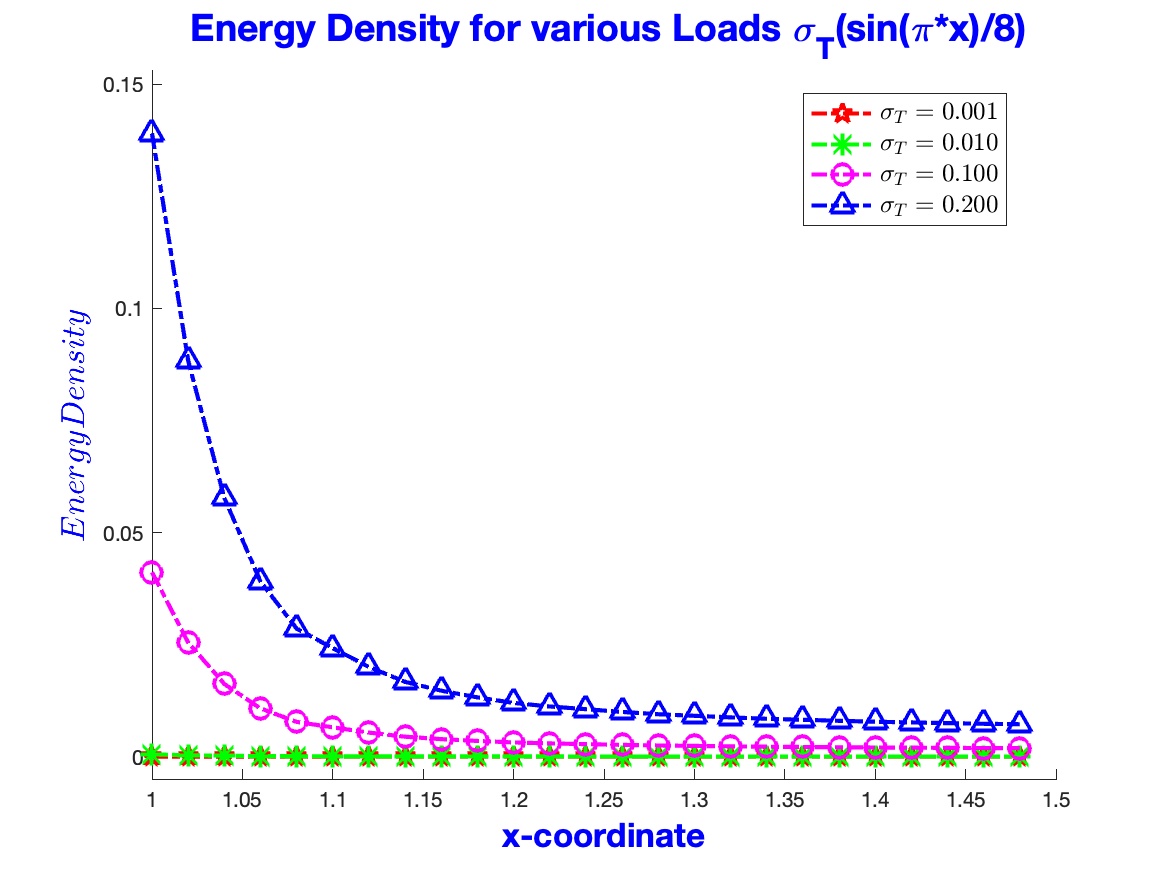}
        \caption{Energy density for various $\sigma_{T}$ for $\beta = 1.0$ and $\alpha = 1.0$}
        \label{fig:energy_density_sigma}
    \end{subfigure}
    \caption{Energy density plots for different parameter variations in case 2(a).}
    \label{fig:energy_2a}
\end{figure}

A parametric study was conducted to investigate the sensitivity of the local crack-tip response to variations in the model parameters $\beta$, $\alpha$, and $\sigma_T$. The resulting distributions for key mechanical fields, including stress, strain, and strain energy density, are illustrated in Figures~\ref{fig:stress_2a}--\ref{fig:energy_2a}. While these results generally demonstrate a clear dependence on the chosen parameters, a noteworthy exception is observed for the crack opening stress, $\sigma_{yy}$. The analysis reveals a remarkable finding: the distribution of this specific stress component ahead of the crack tip appears to be entirely independent of the values of $\beta$, $\alpha$, and $\sigma_T$ within the ranges investigated. Both crack-tip strains and strain energy density decrease with increasing values of $\beta$, and this result indicates that the parameter $\beta$ plays a crucial role in mitigating the severity of the local deformation at the crack tip. The reduction in strain energy density is particularly significant, as it implies that more energy is required from the external load to advance the crack, thus increasing the material's overall fracture toughness. Consequently, $\beta$ can be interpreted as a key toughening parameter in the model, where higher values correspond to enhanced resistance against crack propagation. This insight is valuable from a material design standpoint, as it identifies a direct pathway for improving the durability and failure resistance of the material. In contrast, an opposite trend is observed for the parameters $\alpha$ and $\sigma_T$. These parameters have an amplifying effect on the mechanical fields near the crack tip, exhibiting a direct correlation where an increase in their values results in higher magnitudes of local strain and strain energy density. This trend suggests these parameters contribute to a reduction in the material's overall fracture resistance.

\begin{figure}[H]
    \centering
    \begin{subfigure}{0.45\linewidth}
        \centering
        \includegraphics[width=\linewidth]{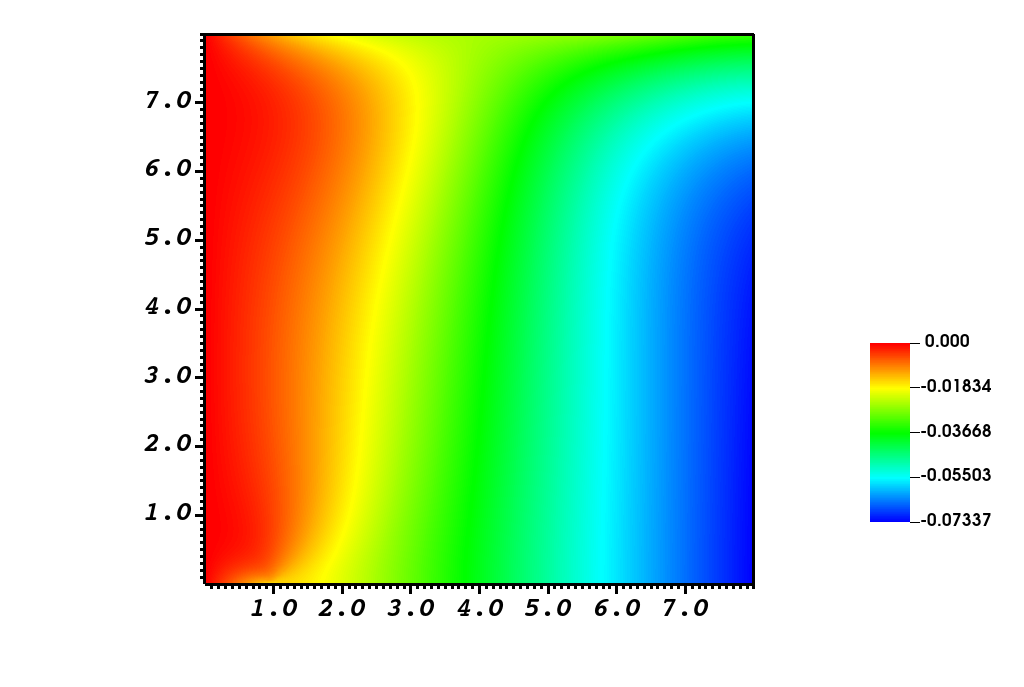}
        \caption{X-displacement at $\alpha = 1.0$, $\sigma_{T} = 0.1$, and $\beta = 1.0$}
        \label{fig:x_displacement}
    \end{subfigure}
    \hfill
    \begin{subfigure}{0.45\linewidth}
        \centering
        \includegraphics[width=\linewidth]{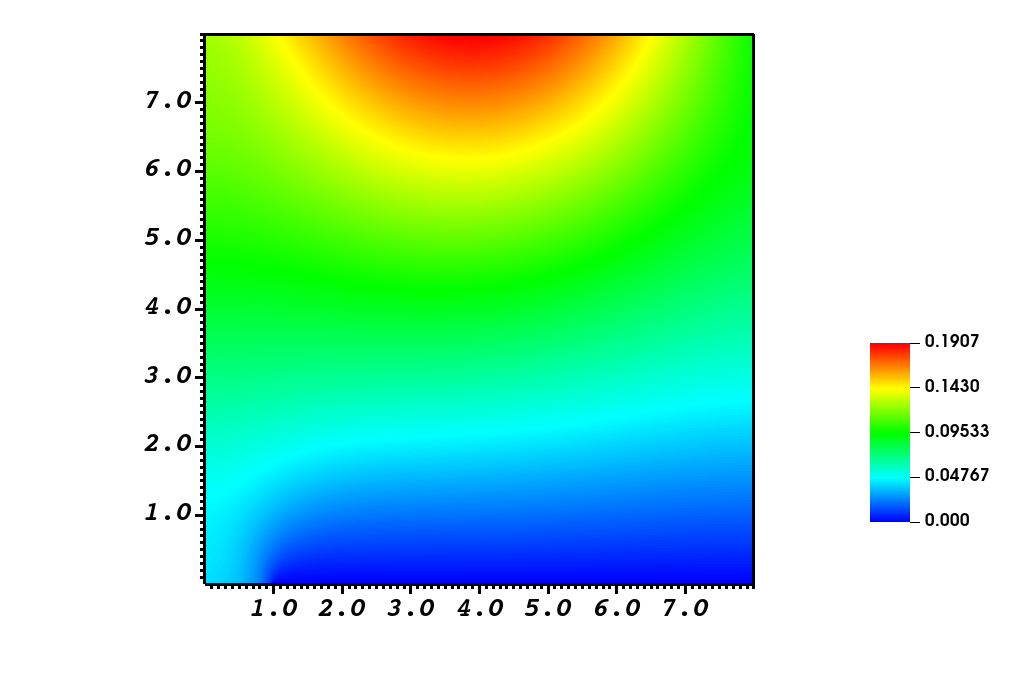}
        \caption{Y-displacement at $\alpha = 1.0$, $\sigma_{T} = 0.1$, and $\beta = 1.0$}
        \label{fig:y_displacement}
    \end{subfigure}
    \vspace{1em}
    \begin{subfigure}{0.45\linewidth}
        \centering
        \includegraphics[width=\linewidth]{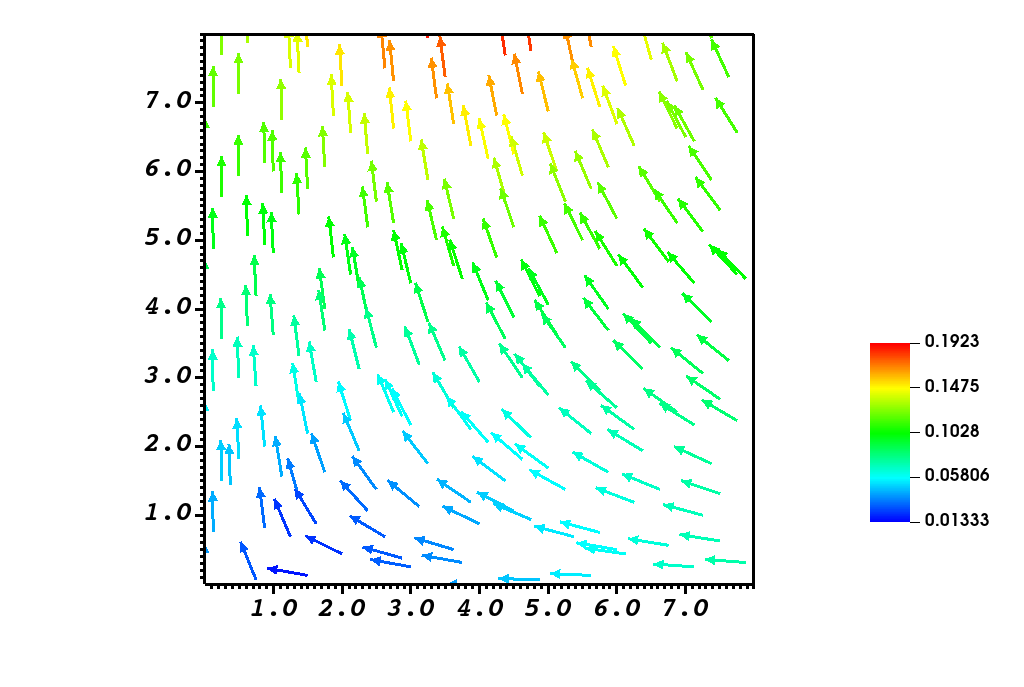}
        \caption{Vector-displacement at $\alpha = 1.0$, $\sigma_{T} = 0.1$, and $\beta = 1.0$}
        \label{fig:vector_displacement}
    \end{subfigure}
    \hfill
    \begin{subfigure}{0.45\linewidth}
        \centering
        \includegraphics[width=\linewidth]{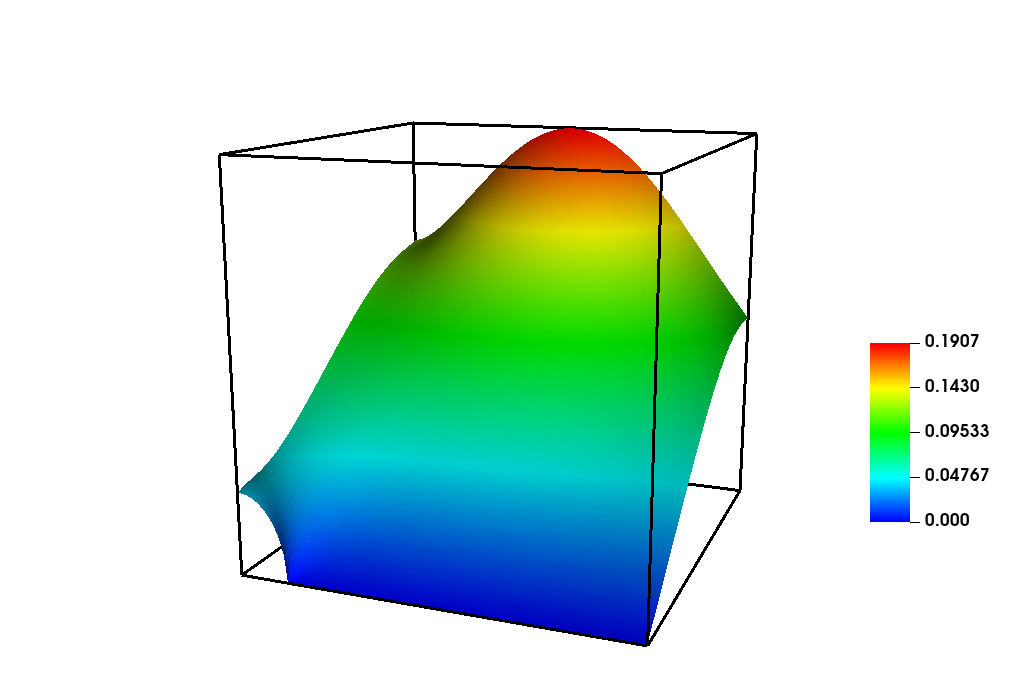}
        \caption{Y-displacement at $\alpha = 1.0$, $\sigma_{T} = 0.1$, and $\beta = 1.0$ (3D view)}
        \label{fig:y_displacement_3d}
    \end{subfigure}

    \caption{Displacement plots for $\alpha = 1.0$, $\sigma_{T} = 0.1$, and $\beta = 1.0$ for load $\sigma_{T}(\frac{sin(\pi x)}{8})$ in x-fiber direction}
    \label{fig:displacement_2a}
\end{figure}

The complete displacement field is presented in Figure~\ref{fig:displacement_2a}. This figure illustrates the individual displacement components in the $x$ and $y$ directions, alongside the overall displacement vector field. The vertical displacement along the crack faces is then extracted from this data to construct the crack opening profile, shown separately in Figure~\ref{fig:displacement_2a}. A detailed examination of this profile reveals two key characteristics that highlight the material's elastic-plastic response. On a global scale, the opening conforms to the classic {elliptical shape} predicted by Linear Elastic Fracture Mechanics (LEFM). 

\subsection{Case 2(b)-Non-uniform load: Fibers aligned orthogonal with the crack  plane}
This section presents the numerical simulation results for a benchmark case: an edge-cracked {orthotropic} material subjected to {Mode-I non-uniform loading}. A key aspect of this investigation is the fiber orientation, which is aligned perpendicular to the crack plane. This configuration is strategically chosen as it represents an optimal arrangement for maximizing {fracture toughness}, where the stiff fibers act as the primary load-bearing constituents, effectively bridging the crack faces and inhibiting their separation. For the numerical solution, the material's Poisson's ratio was held constant, and the resulting nonlinear system of equations was solved using a Picard iterative scheme. The solver demonstrated robust and efficient performance for this configuration, achieving rapid monotonic convergence. The convergence data presented in Table~\ref{table4}, corresponding to the aforementioned fiber-crack orientation, confirm the stability and effectiveness of our iterative method for this class of problems.  The computations done for the above table are for the parameters $\alpha=1.0$, $\beta=1.0$, and $\sigma_T=0.1$. 

\begin{table}[H]
\centering
\begin{tabular}{|c|c|}
\hline 
Iteration No & Residual\tabularnewline
\hline 
\hline 
1 & 0.000143355\tabularnewline
\hline 
2 & 1.68854e-05\tabularnewline
\hline 
3 & 2.07079e-06\tabularnewline
\hline 
4 & 5.87719e-07\tabularnewline
\hline 
5 & 6.10429e-07\tabularnewline
\hline 
6 & 6.04246e-07\tabularnewline
\hline 
7 & 6.04868e-07\tabularnewline
\hline 
8 & 6.04803e-07\tabularnewline
\hline 
9 & 6.0481e-07\tabularnewline
\hline 
10 & 6.04809e-07\tabularnewline
\hline 
\end{tabular}
\caption{Solver performance for non-uniform load with fibers orthogonal to the crack.}
\label{table4}
\end{table}

\begin{figure}[H]
    \centering
    \begin{subfigure}{0.3\linewidth}
        \centering
        \includegraphics[width=\linewidth]{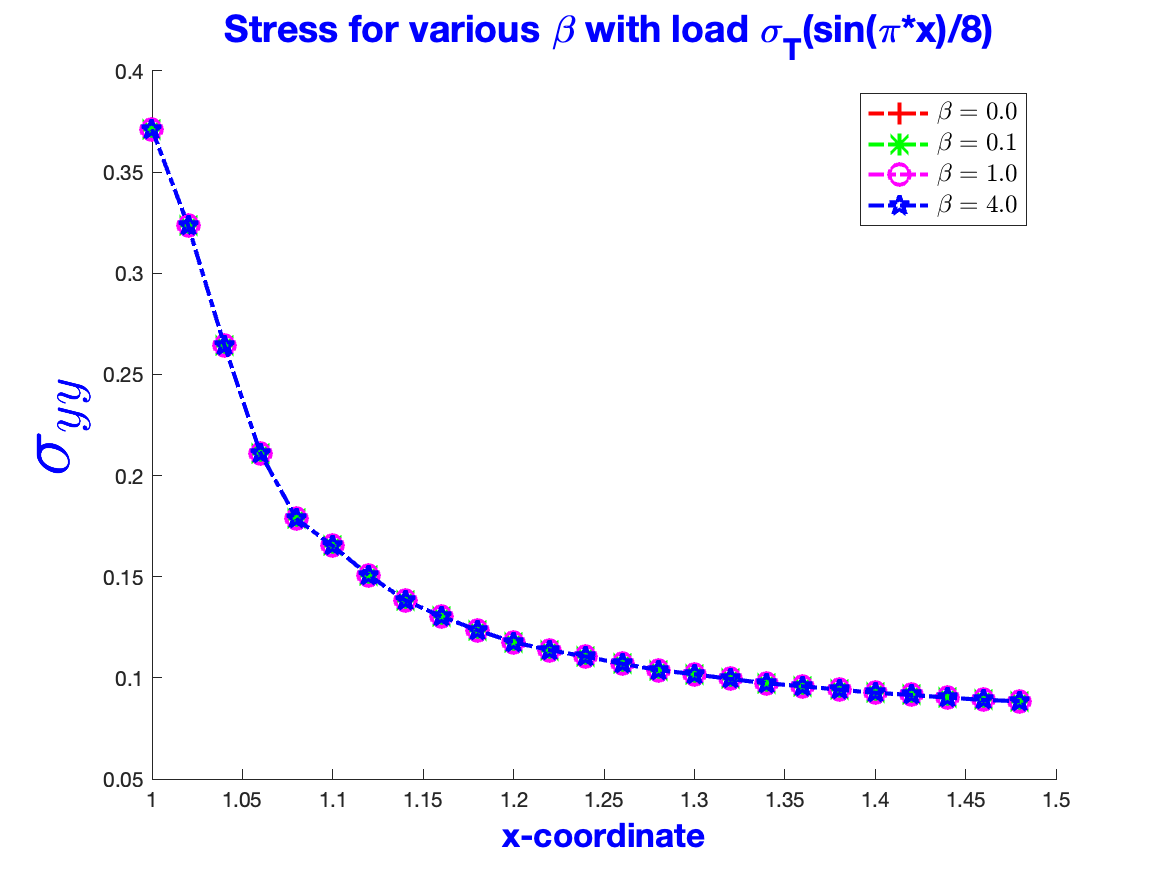}
        \caption{Stress for various $\beta$ for $\alpha = 1.0$ and $\sigma_{T} = 0.1$}
        \label{fig:strain_beta}
    \end{subfigure}
    \hfill
    \begin{subfigure}{0.3\linewidth}
        \centering
        \includegraphics[width=\linewidth]{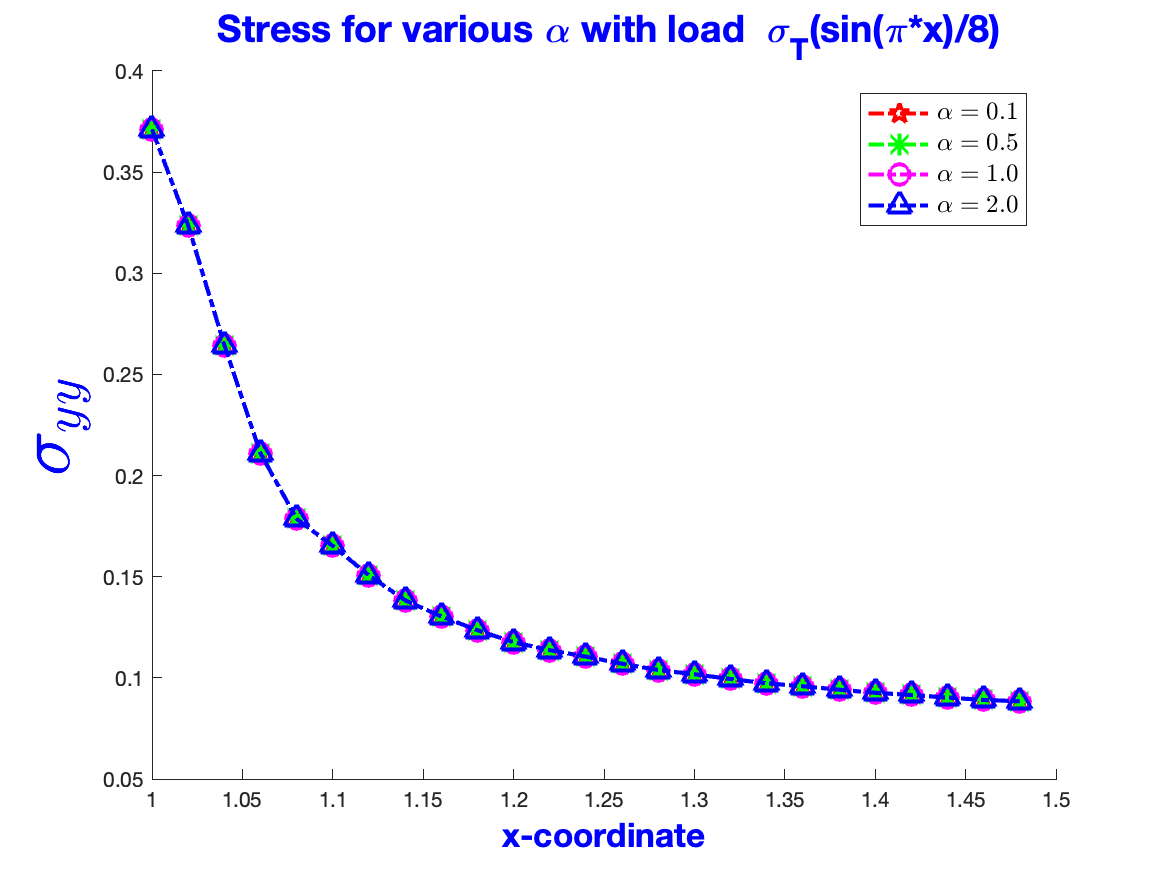}
        \caption{Stress for various $\alpha$ for $\beta = 1.0$ and $\sigma_{T} = 0.1$}
        \label{fig:strain_alpha}
    \end{subfigure}
    \hfill
    \begin{subfigure}{0.3\linewidth}
        \centering
        \includegraphics[width=\linewidth]{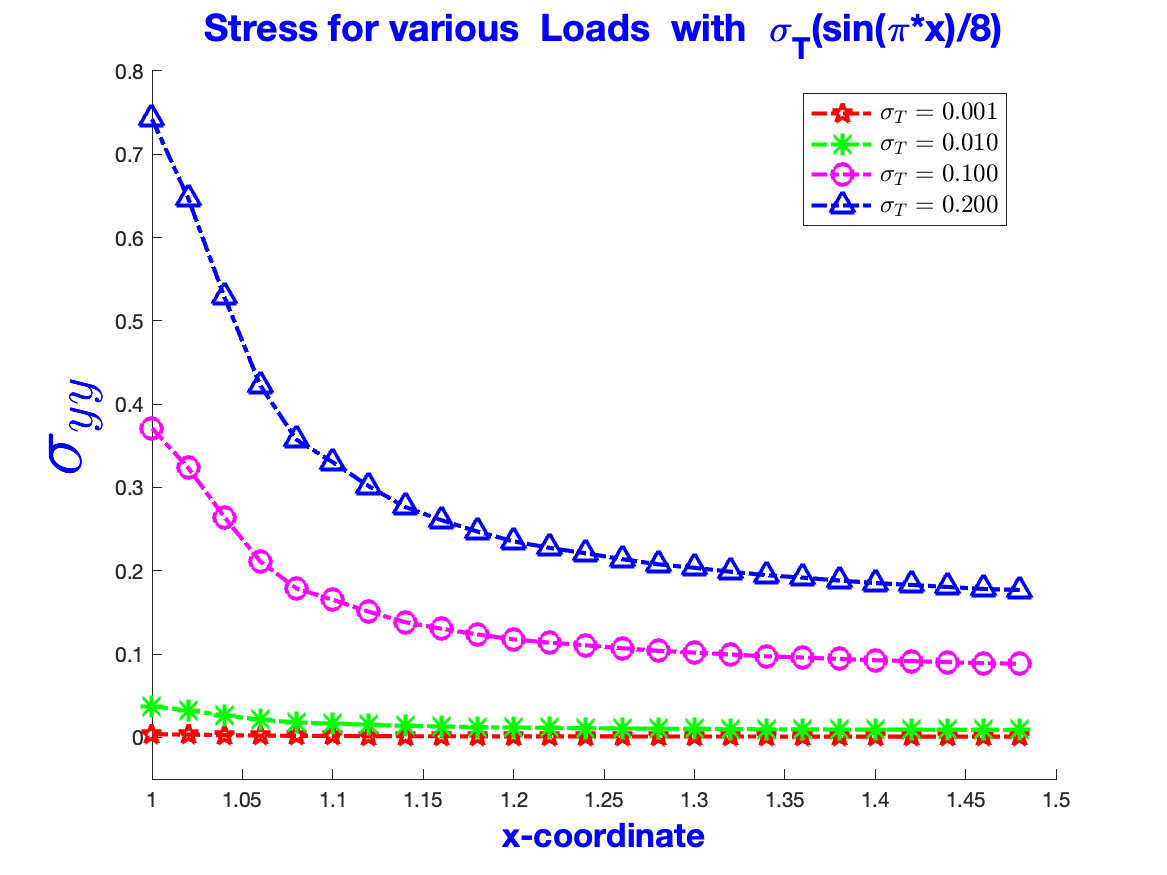}
        \caption{Stress for various $\sigma_{T}$ for $\beta = 1.0$ and $\alpha = 1.0$}
        \label{fig:strain_sigma}
    \end{subfigure}
    \caption{Stress plots for different parameter variations for non-uniform load with fibers orthogonal to the crack.}
    \label{fig:stress_2b}
\end{figure}
\begin{figure}[H]
    \centering
    \begin{subfigure}{0.3\linewidth}
        \centering
        \includegraphics[width=\linewidth]{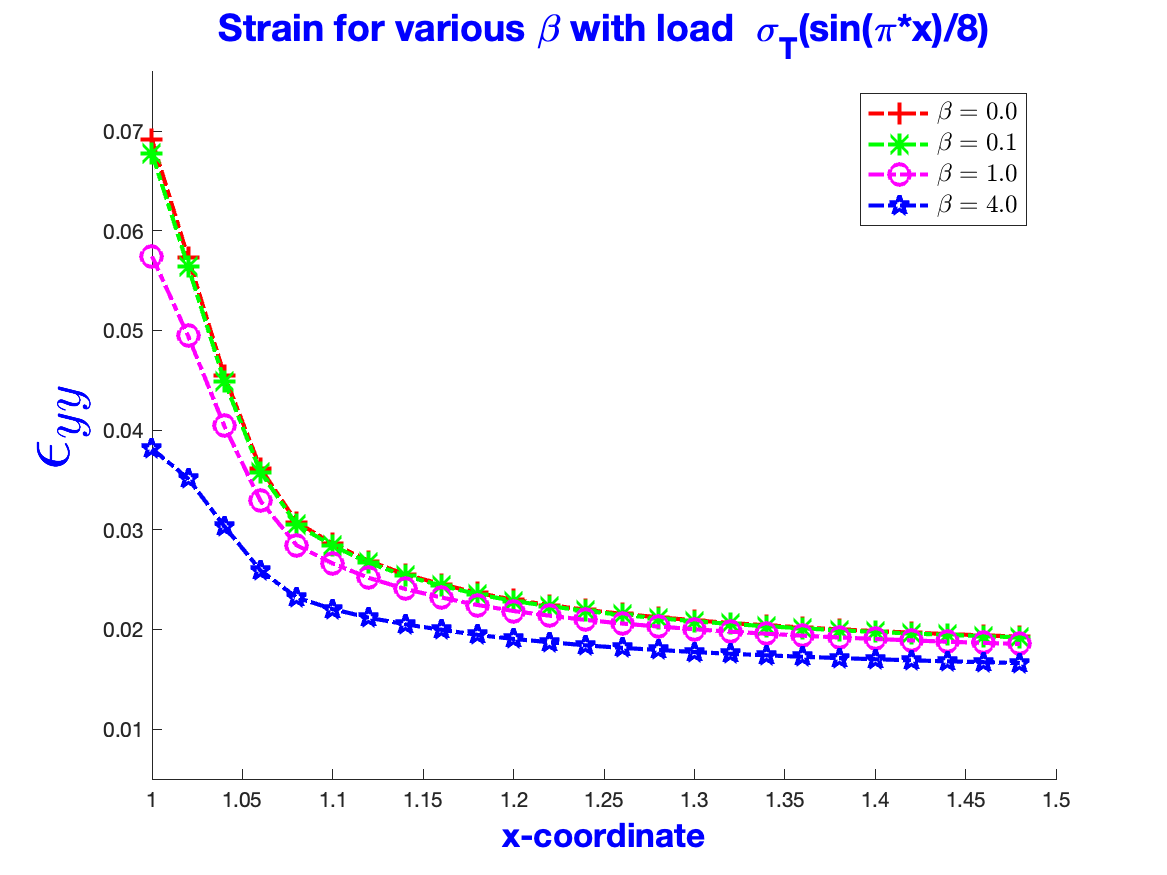}
        \caption{Strain for various $\beta$ for $\alpha = 1.0$ and $\sigma_{T} = 0.1$}
        \label{fig:strain_beta}
    \end{subfigure}
    \hfill
    \begin{subfigure}{0.3\linewidth}
        \centering
        \includegraphics[width=\linewidth]{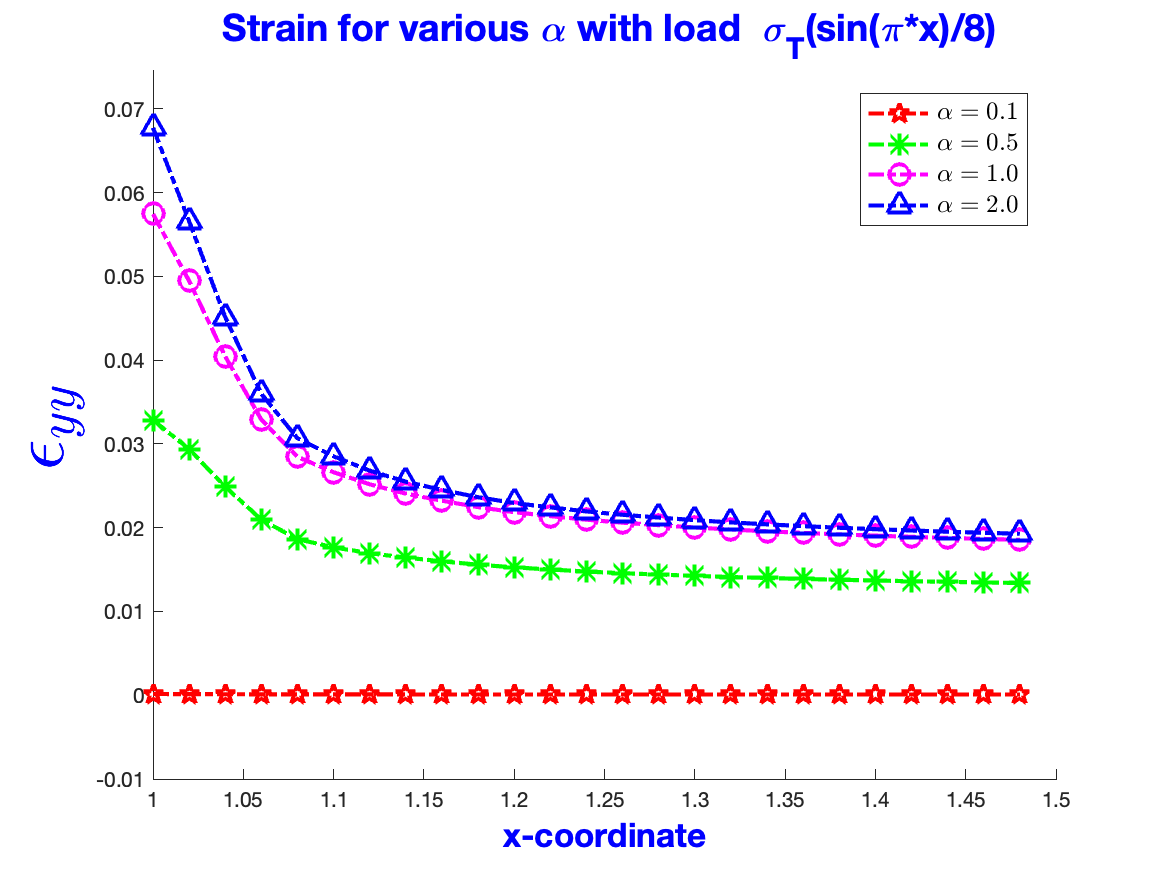}
        \caption{Strain for various $\alpha$ for $\beta = 1.0$ and $\sigma_{T} = 0.1$}
        \label{fig:strain_alpha}
    \end{subfigure}
    \hfill
    \begin{subfigure}{0.3\linewidth}
        \centering
        \includegraphics[width=\linewidth]{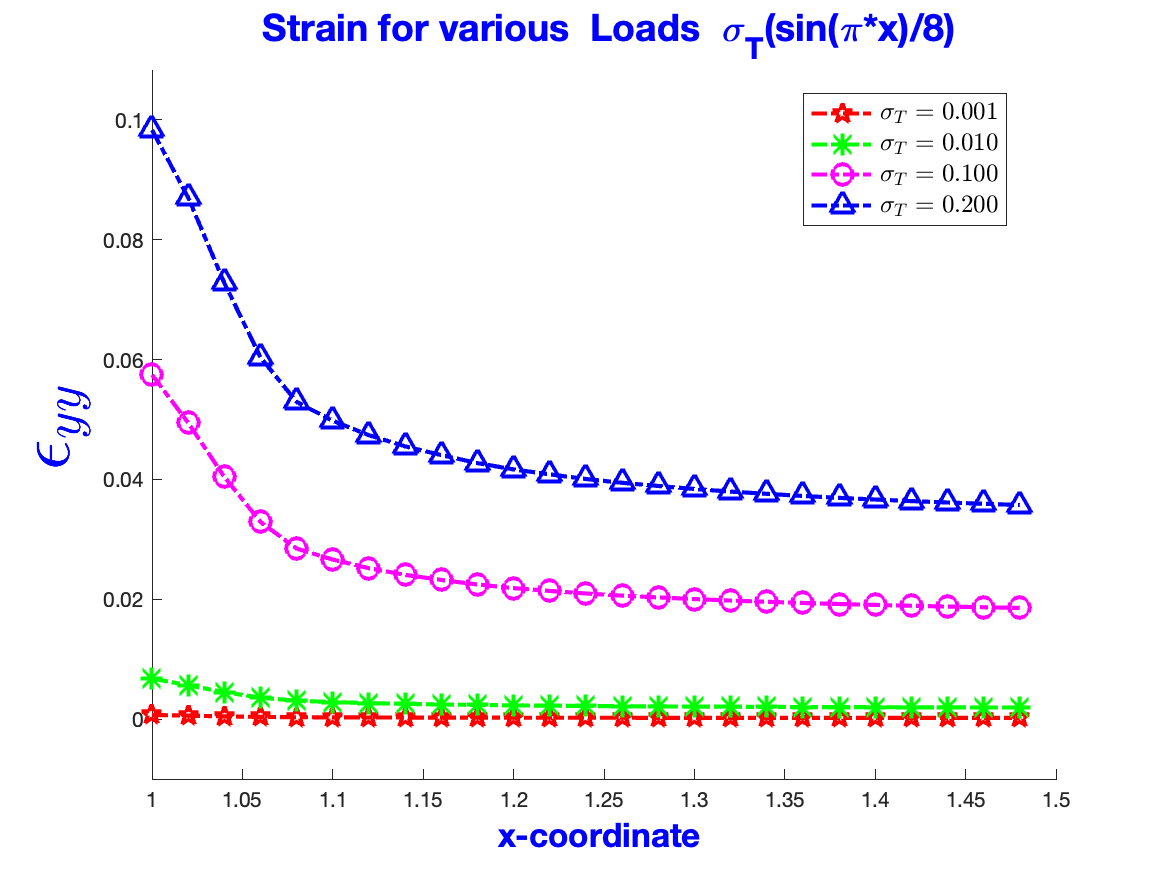}
        \caption{Strain for various $\sigma_{T}$ for $\beta = 1.0$ and $\alpha = 1.0$}
        \label{fig:strain_sigma}
    \end{subfigure}
    \caption{Strain plots for different parameter variations for non-uniform load with fibers orthogonal to the crack.}
    \label{fig:strain_2b}
\end{figure}

\begin{figure}[H]
    \centering
    \begin{subfigure}{0.3\linewidth}
        \centering
        \includegraphics[width=\linewidth]{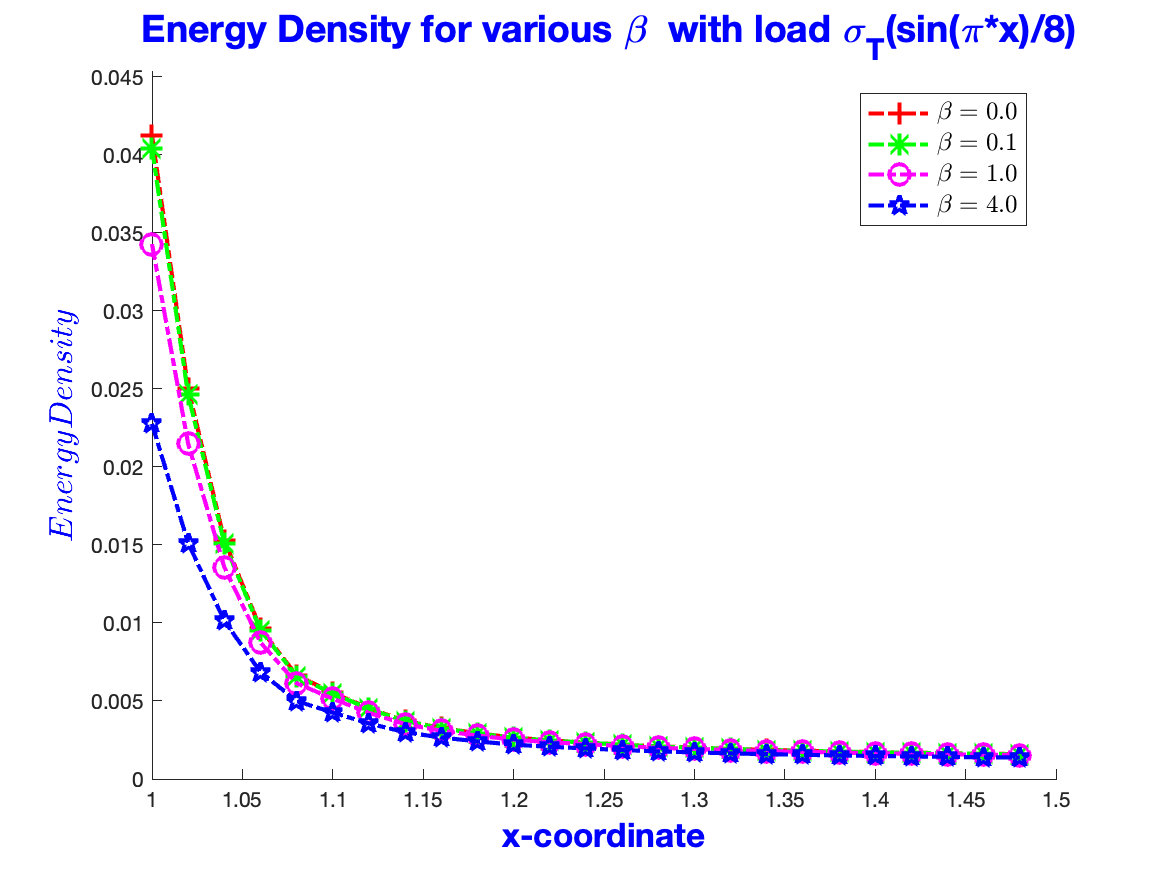}
        \caption{Energy density for various $\beta$ for $\alpha = 1.0$ and $\sigma_{T} = 0.1$}
        \label{fig:energy_density_beta}
    \end{subfigure}
    \hfill
    \begin{subfigure}{0.3\linewidth}
        \centering
        \includegraphics[width=\linewidth]{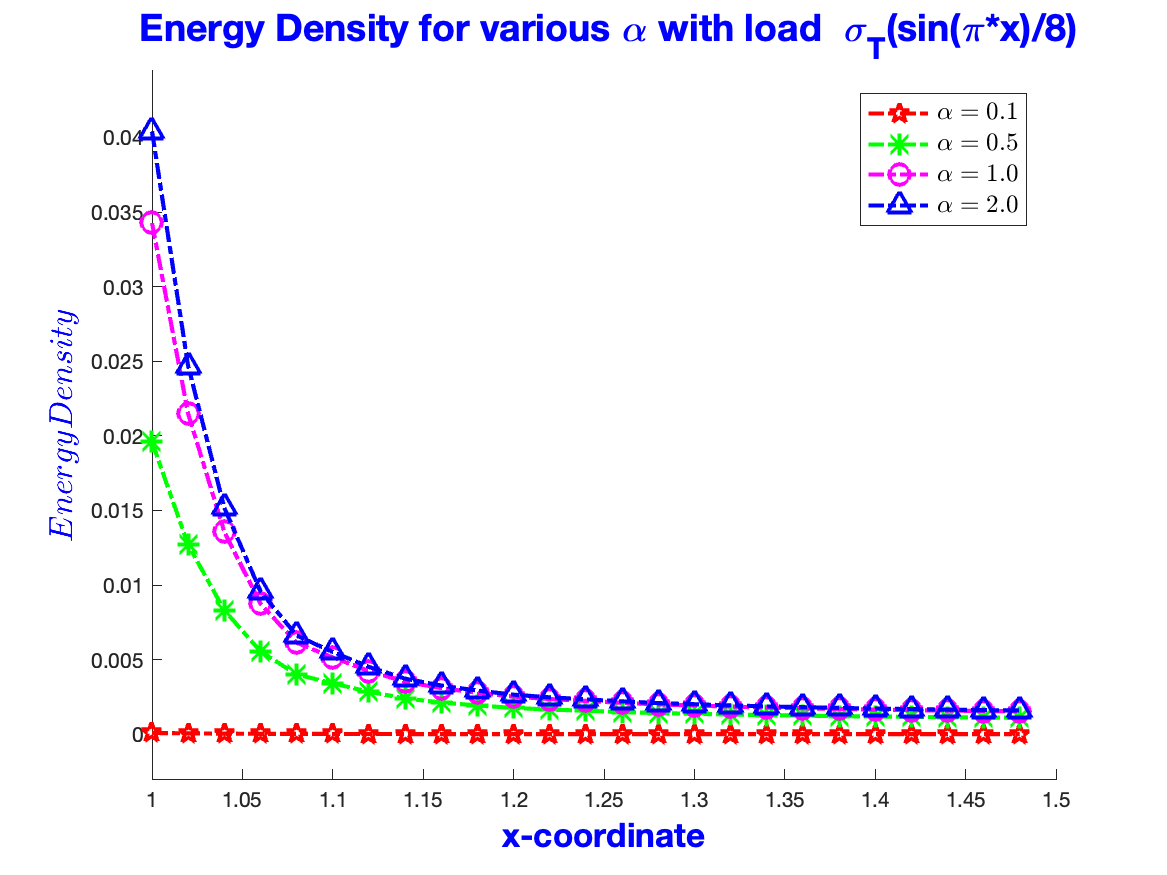}
        \caption{Energy density for various $\alpha$ for $\beta = 1.0$ and $\sigma_{T} = 0.1$}
        \label{fig:energy_density_alpha}
    \end{subfigure}
    \hfill
    \begin{subfigure}{0.3\linewidth}
        \centering
        \includegraphics[width=\linewidth]{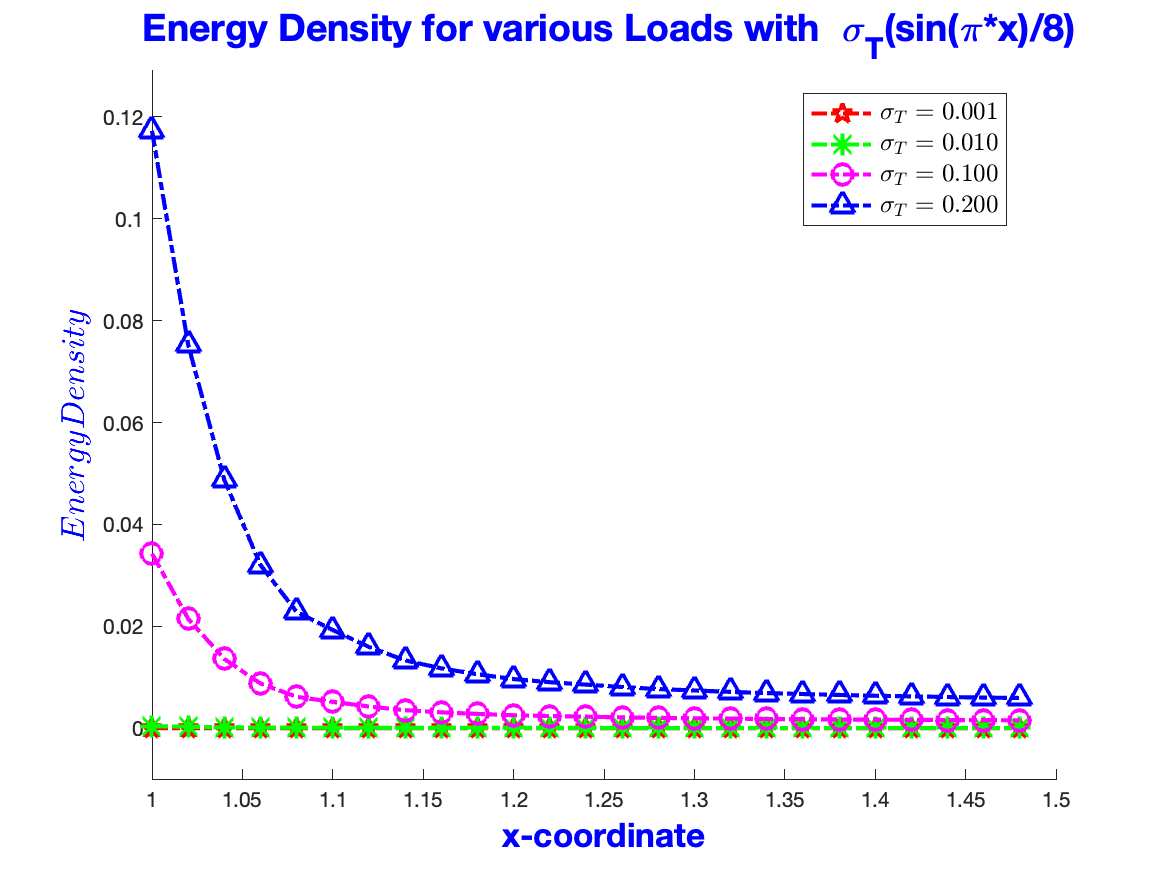}
        \caption{Energy density for various $\sigma_{T}$ for $\beta = 1.0$ and $\alpha = 1.0$}
        \label{fig:energy_density_sigma}
    \end{subfigure}
    \caption{Energy density plots for different parameter variations for non-uniform load with fibers orthogonal to the crack.}
    \label{fig:energy_density_2b}
\end{figure}

Figures~\ref{fig:stress_2b}, \ref{fig:strain_2b}, and \ref{fig:energy_density_2b} illustrate the distinct influence of the model parameters $\beta$, $\alpha$, and $\sigma_T$ on the mechanical fields at the crack tip. The parameter $\beta$ plays a crucial role in mitigating the severity of local deformation. While its impact on stress is minor, increasing the value of $\beta$ causes a significant decrease in both crack-tip strains and strain energy density. The reduction in strain energy density is particularly important, as it implies that more energy is required from the external load to advance the crack. This directly increases the material’s overall fracture toughness. Consequently, $\beta$ can be interpreted as a key toughening parameter in the model, where higher values correspond to enhanced resistance against crack propagation. In contrast, an opposite trend is observed for the parameters $\alpha$ and $\sigma_T$. These parameters have an amplifying effect on the mechanical fields near the crack tip. They exhibit a direct correlation where an increase in their values results in higher magnitudes of local strain and strain energy density. This trend suggests that $\alpha$ and $\sigma_T$ contribute to a reduction in the material’s overall fracture resistance, as they intensify the local conditions that drive crack growth.

\begin{figure}[H]
    \centering
    \begin{subfigure}{0.45\linewidth}
        \centering
        \includegraphics[width=\linewidth]{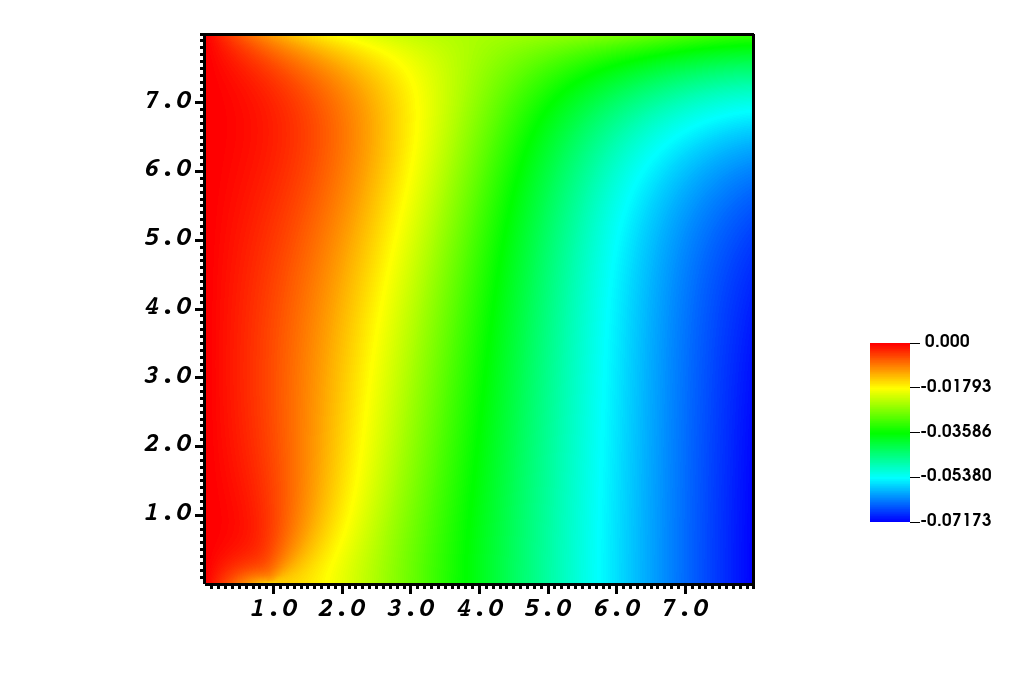}
        \caption{$x$-displacement for $\alpha = 1.0$, $\sigma_{T} = 0.1$, and $\beta = 1.0$}
        \label{fig:x_displacement}
    \end{subfigure}
    \hfill
    \begin{subfigure}{0.45\linewidth}
        \centering
        \includegraphics[width=\linewidth]{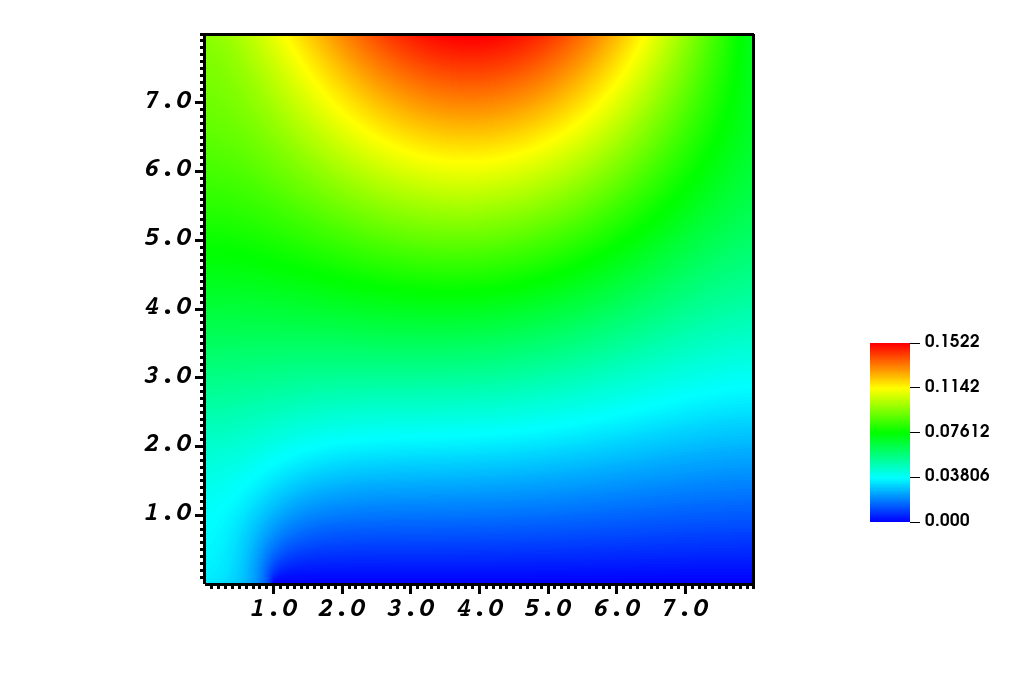}
        \caption{$y$-displacement for $\alpha = 1.0$, $\sigma_{T} = 0.1$, and $\beta = 1.0$}
        \label{fig:y_displacement}
    \end{subfigure}
    
    \vspace{1em}
    
    \begin{subfigure}{0.45\linewidth}
        \centering
        \includegraphics[width=\linewidth]{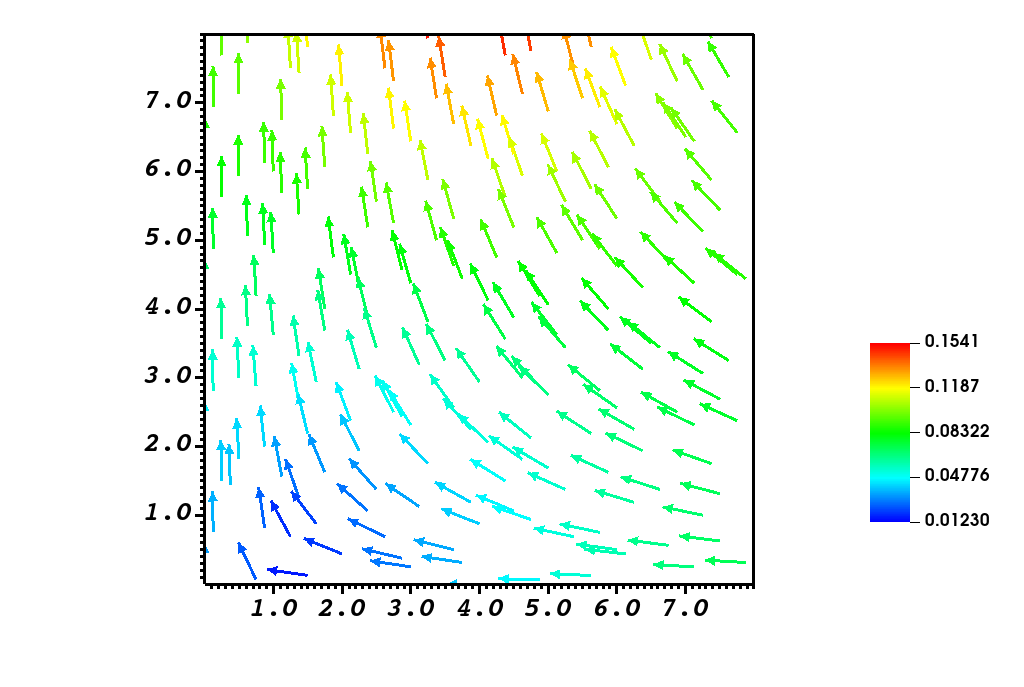}
        \caption{Vector-displacement for $\alpha = 1.0$, $\sigma_{T} = 0.1$, and $\beta = 1.0$}
        \label{fig:vector_displacement}
    \end{subfigure}
    \hfill
    \begin{subfigure}{0.45\linewidth}
        \centering
        \includegraphics[width=\linewidth]{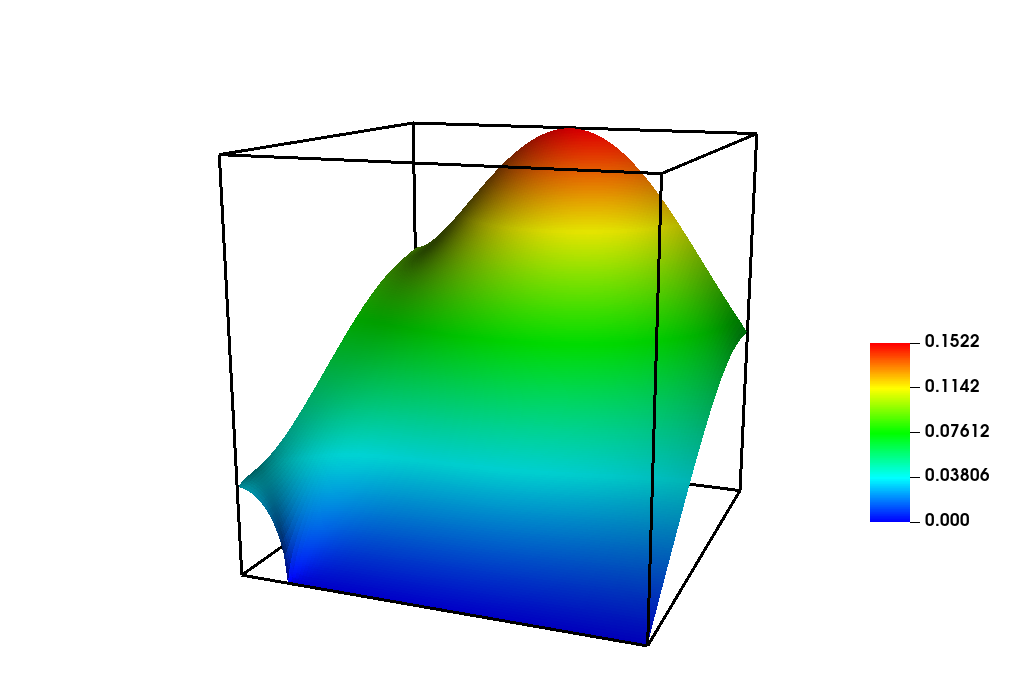}
        \caption{$y$-displacement for $\alpha = 1.0$, $\sigma_{T} = 0.1$, and $\beta = 1.0$ (3D view)}
        \label{fig:y_displacement_3d}
    \end{subfigure}

    \caption{Displacement plots for $\alpha = 1.0$, $\sigma_{T} = 0.1$, and $\beta = 1.0$ for load $\sigma_{T}(\frac{sin(\pi x)}{8})$ for material fibers assumed to in  $y$-direction. }
    \label{fig:displacement_2b}
\end{figure}

Figure~\ref{fig:displacement_2a} presents the comprehensive displacement field obtained from the simulation. This visualization provides a detailed breakdown of the individual displacement components along the horizontal ($u_x$) and vertical ($u_y$) axes, as well as the overall displacement vector field. From this full-field data, the vertical displacements ($u_y$) along the $y$-axis are extracted to construct the crack opening profile. This profile, also depicted in Figure~\ref{fig:displacement_2a}, is crucial for understanding the material's response. A detailed examination of its geometry reveals key characteristics of the underlying elastic-plastic behavior. On a global scale, the crack opening profile closely matches the classic {elliptical shape} predicted by LEFM. This agreement indicates that, far from the high-stress region of the crack tip, the material's deformation is governed primarily by its elastic properties.

\subsection{Effect of loading configuration on crack-tip fields}

This section investigates how different loading configurations affect the mechanical response near a crack tip. We analyze and contrast the effects of a {uniform tensile load} with a {non-uniform tensile load} applied to the top edge of the cracked body. Understanding the differences between these loading scenarios is crucial for several reasons. In many real-world engineering applications, structures are subjected to complex, non-uniform loads rather than idealized, uniform ones. For example, bending, thermal gradients, or contact pressures can induce linearly varying stress states. This study is therefore essential for bridging the gap between theoretical models, which often assume uniform loading, and the practical performance of materials in service. By examining both cases, we can more accurately predict fracture initiation and assess the structural integrity of components under realistic operational conditions.

To provide a comprehensive comparison, this study examines the material's response under several distinct loading conditions. We analyze three specific tensile load profiles: a linearly varying (slope) load defined by $\sigma_T(0.1 + 0.1x)$, a uniform load, and a cyclic load given by $\sigma_T\left(\frac{\sin(\pi x)}{8}\right)$. This comparative analysis is critical because it simulates a range of realistic loading environments beyond simple uniform tension. By investigating both linear and periodic load distributions, we can better understand how stress gradients influence fracture mechanics. Each of these loading scenarios was simulated for two principal fiber orientations to create a thorough comparison of performance.

\begin{figure}[H]
    \centering
    \begin{subfigure}[b]{0.32\linewidth}
        \centering
        \includegraphics[width=\linewidth]{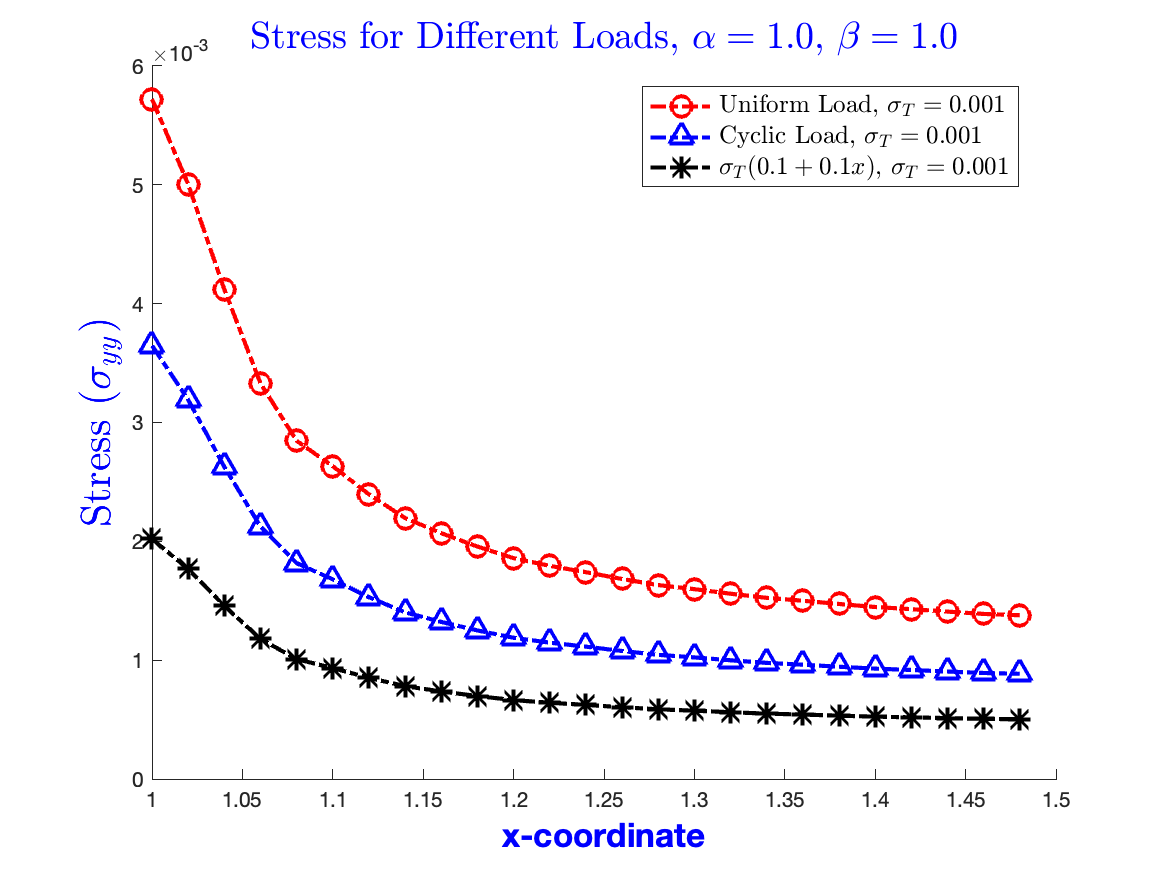}
        \caption{Load = 0.001}
        \label{fig:stress_lp001}
    \end{subfigure}
    \hfill
    \begin{subfigure}[b]{0.32\linewidth}
        \centering
        \includegraphics[width=\linewidth]{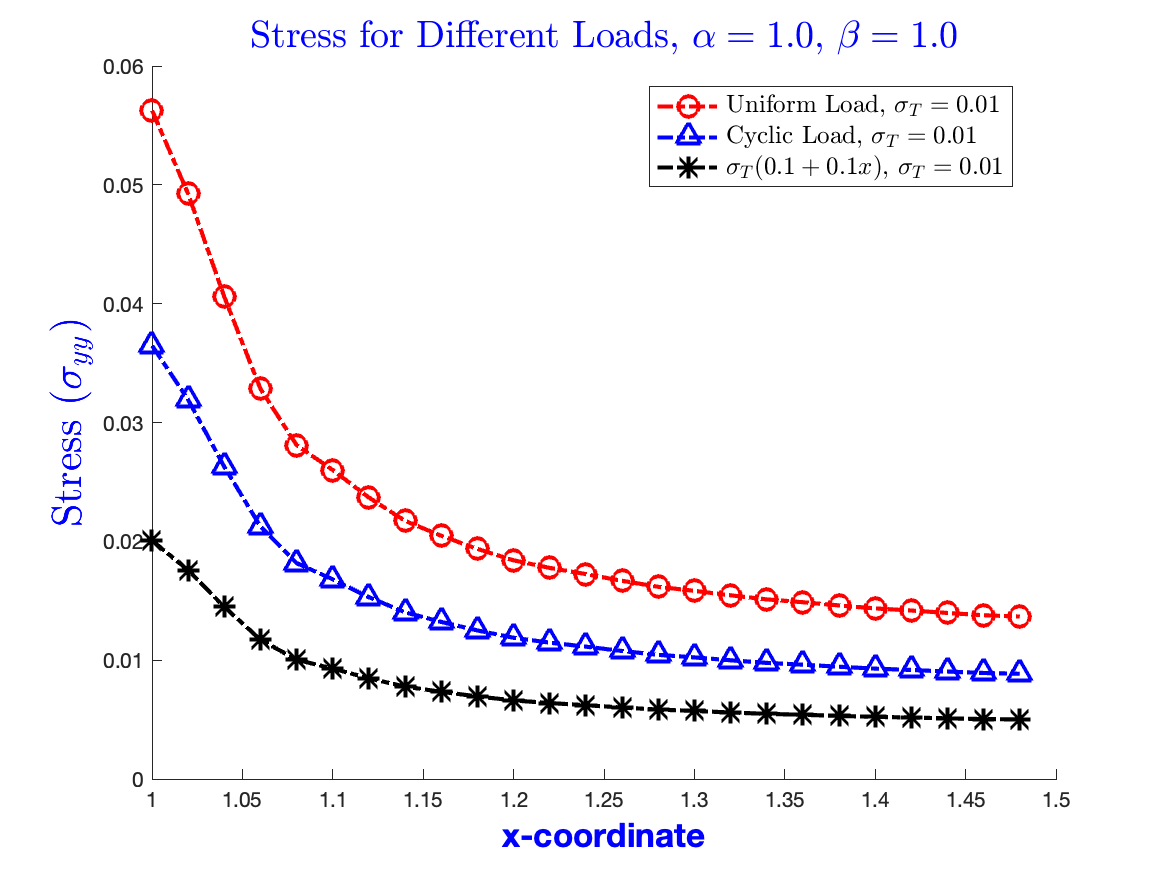}
        \caption{Load = 0.01}
        \label{fig:stress_lp01}
    \end{subfigure}
    \hfill
    \begin{subfigure}[b]{0.32\linewidth}
        \centering
        \includegraphics[width=\linewidth]{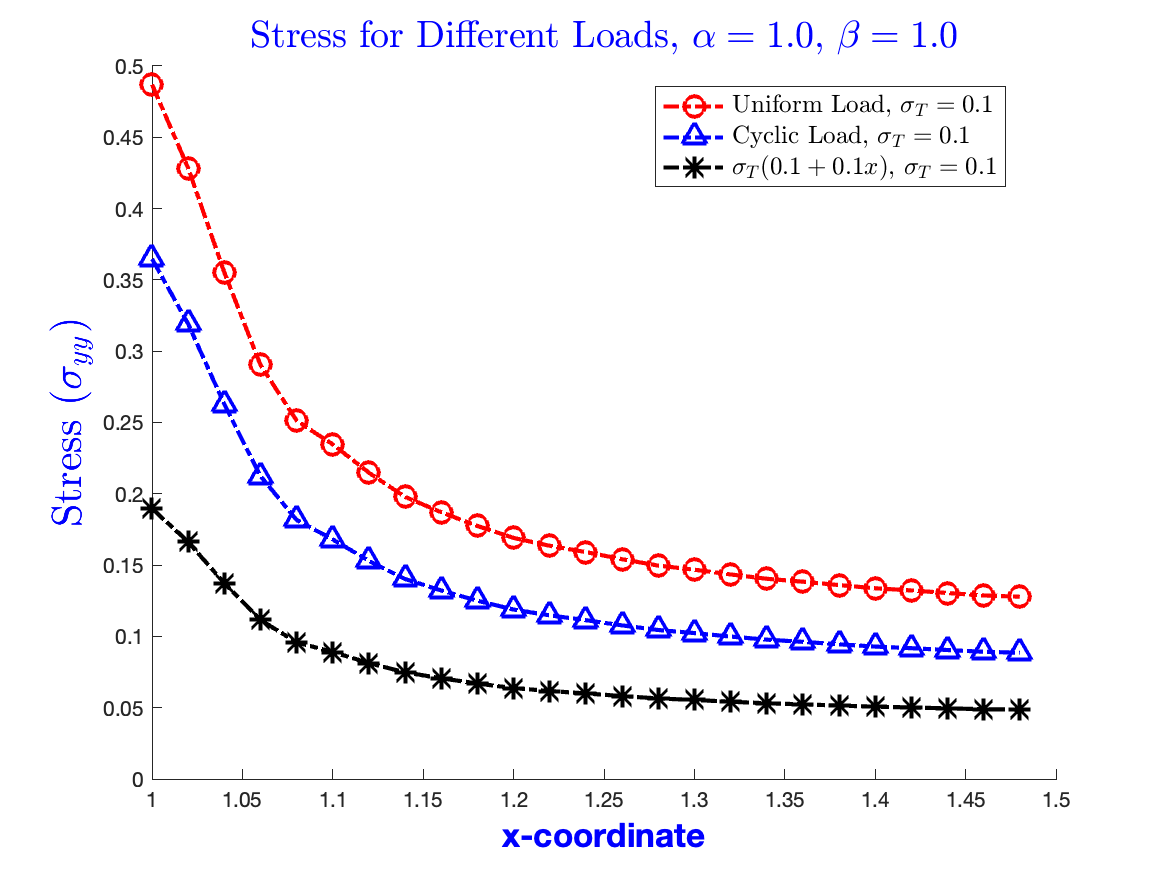}
        \caption{Load = 0.1}
        \label{fig:stress_lp1}
    \end{subfigure}
    \caption{{Comparison of crack-tips stress for various loading scenarios with the material fibers aligned with $x$-axis. }}
    \label{stress_xcase3}
\end{figure}

\begin{figure}[H]
    \centering
    \begin{subfigure}[b]{0.32\linewidth}
        \centering
        \includegraphics[width=\linewidth]{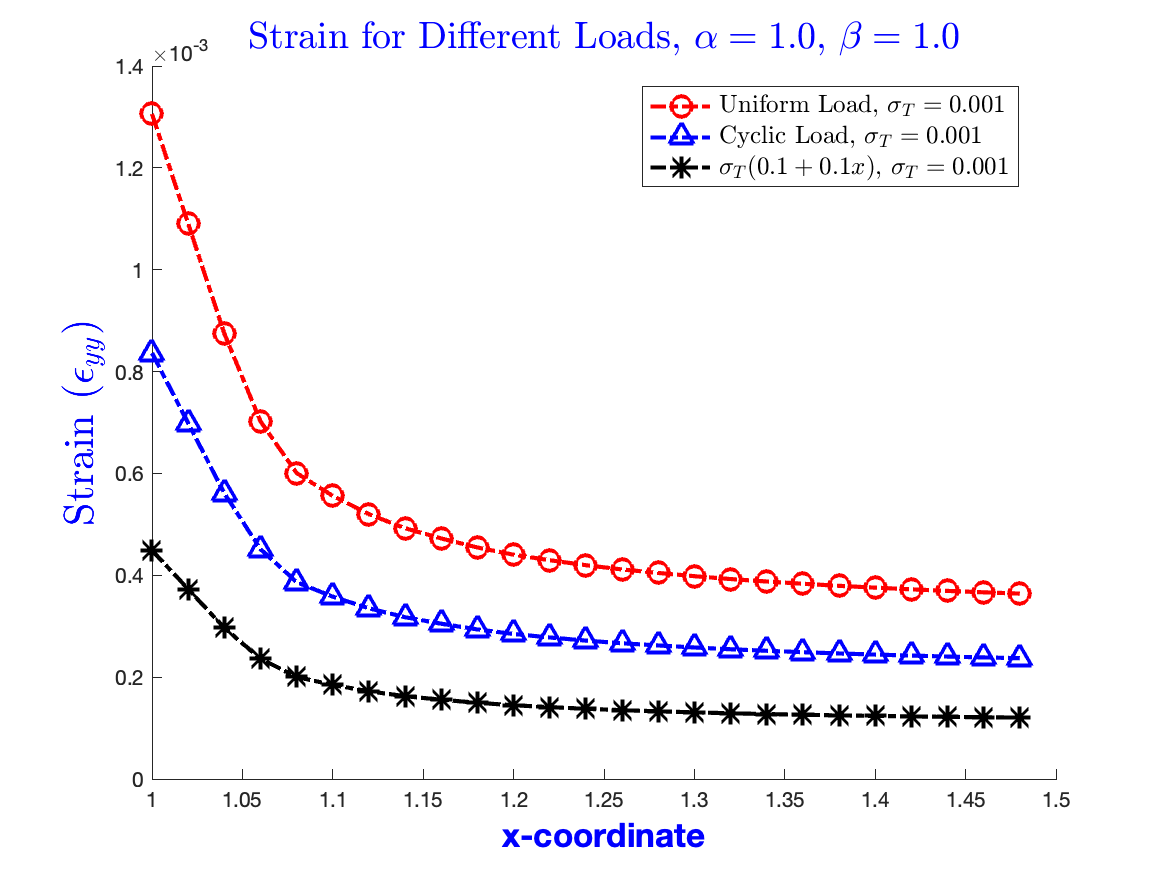}
        \caption{Load = 0.001}
        \label{fig:strain_lp001}
    \end{subfigure}
    \hfill
    \begin{subfigure}[b]{0.32\linewidth}
        \centering
        \includegraphics[width=\linewidth]{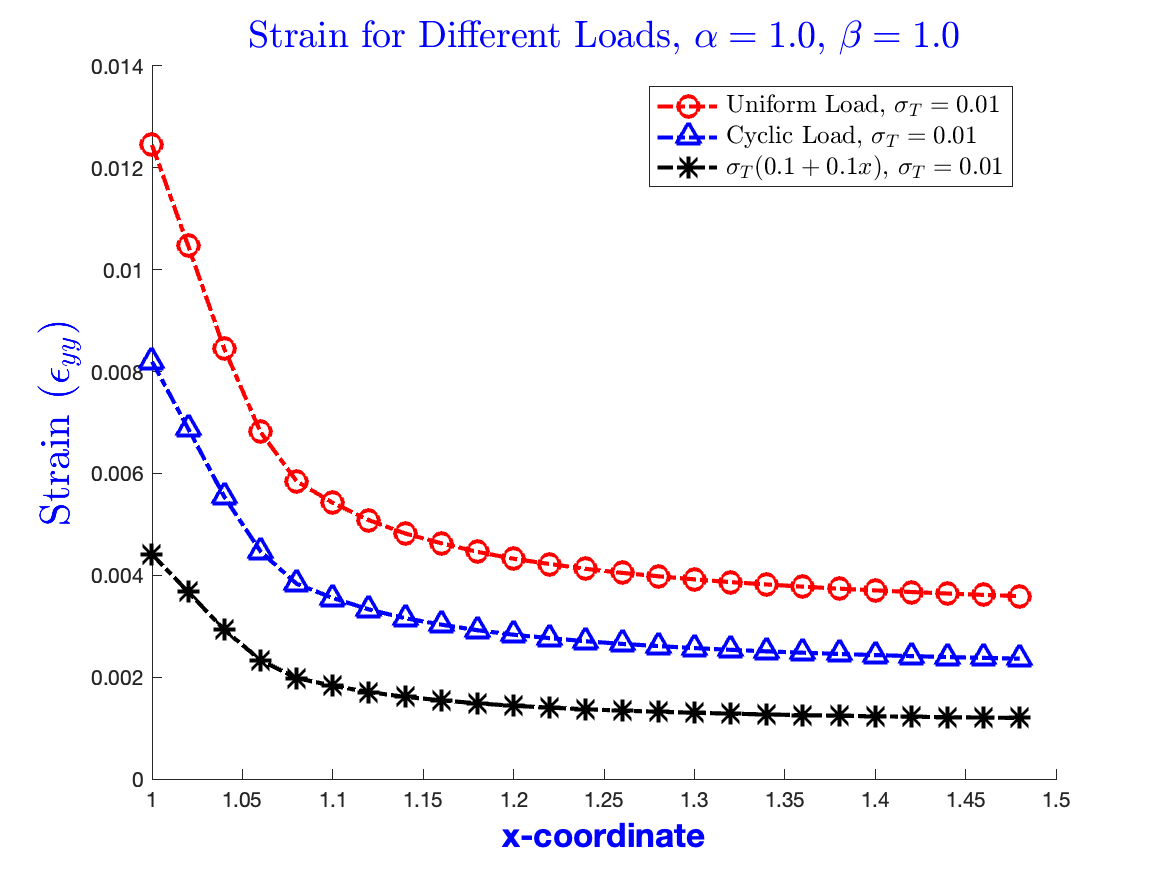}
        \caption{Load = 0.01}
        \label{fig:strain_lp01}
    \end{subfigure}
    \hfill
    \begin{subfigure}[b]{0.32\linewidth}
        \centering
        \includegraphics[width=\linewidth]{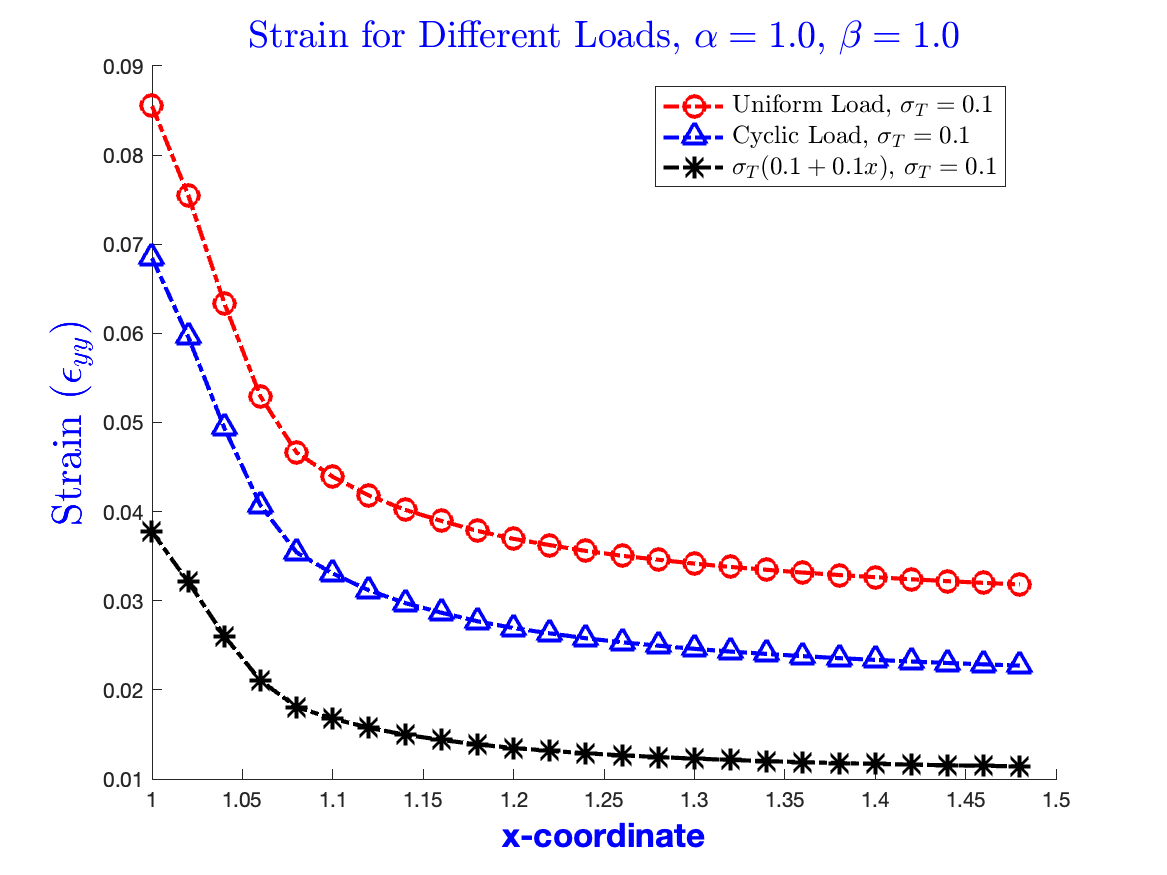}
        \caption{Load = 0.1}
        \label{fig:strain_lp1}
    \end{subfigure}
    \caption{{Comparison of crack-tips strain for various loading scenarios with the material fibers aligned with $x$-axis.}}
    \label{strain_xcase3}
\end{figure}

\begin{figure}[H]
    \centering
    \begin{subfigure}[b]{0.32\linewidth}
        \centering
        \includegraphics[width=\linewidth]{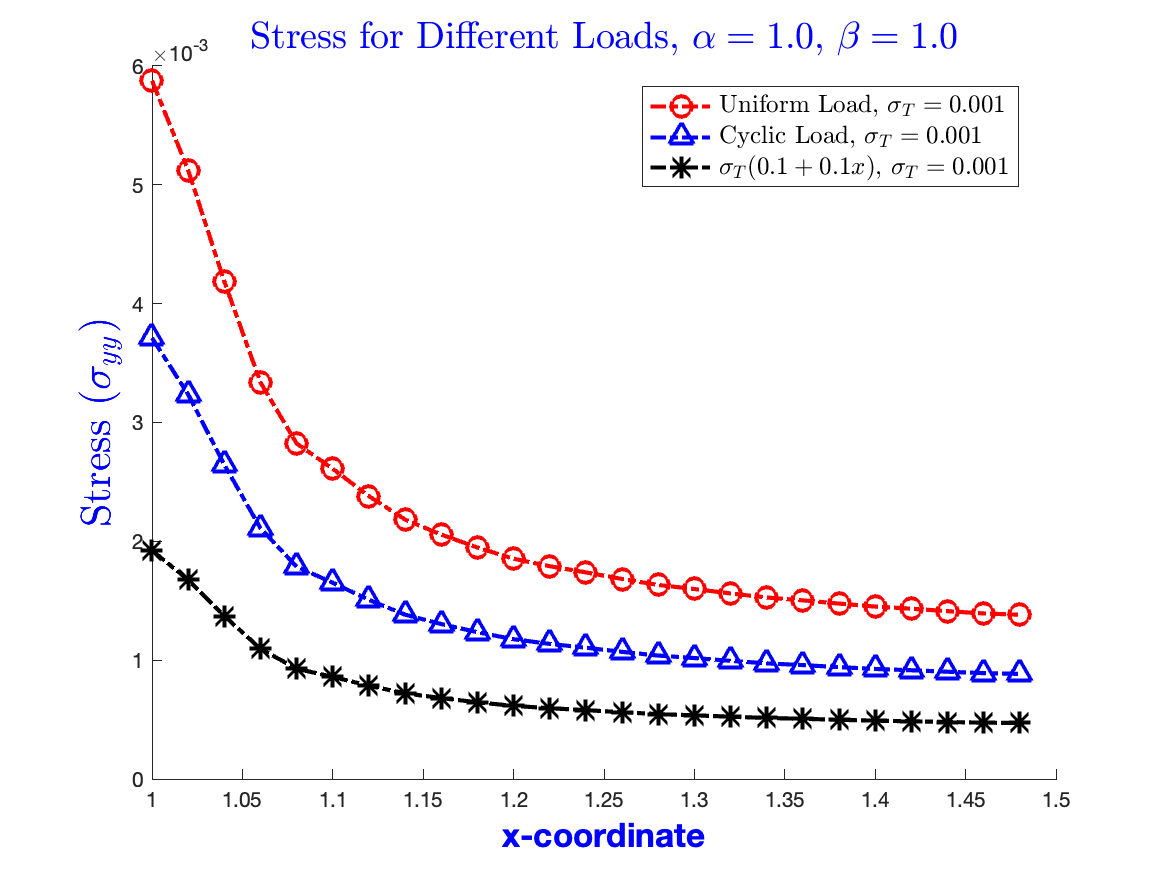}
        \caption{Load = 0.001}
        \label{fig:stress_M2_lp001}
    \end{subfigure}
    \hfill
    \begin{subfigure}[b]{0.32\linewidth}
        \centering
        \includegraphics[width=\linewidth]{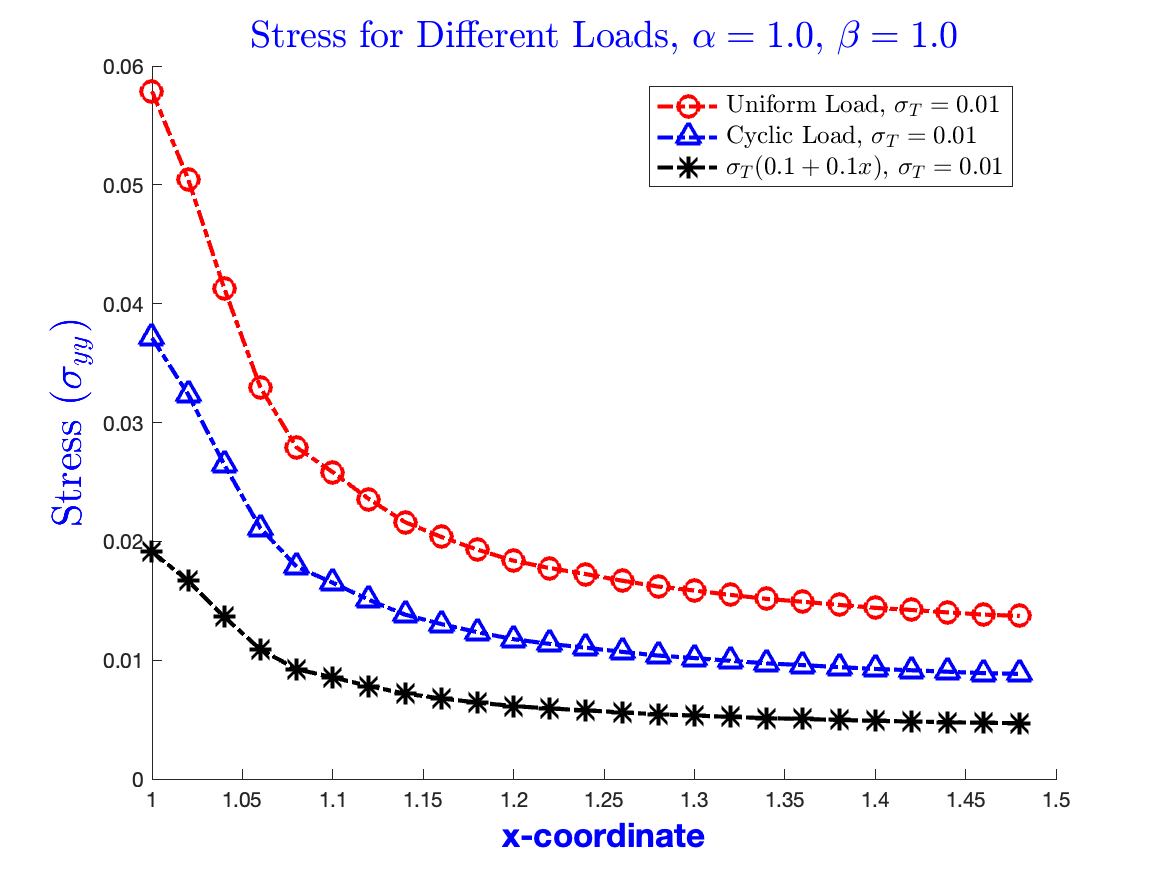}
        \caption{Load = 0.01}
        \label{fig:stress_M2_lp01}
    \end{subfigure}
    \hfill
    \begin{subfigure}[b]{0.32\linewidth}
        \centering
        \includegraphics[width=\linewidth]{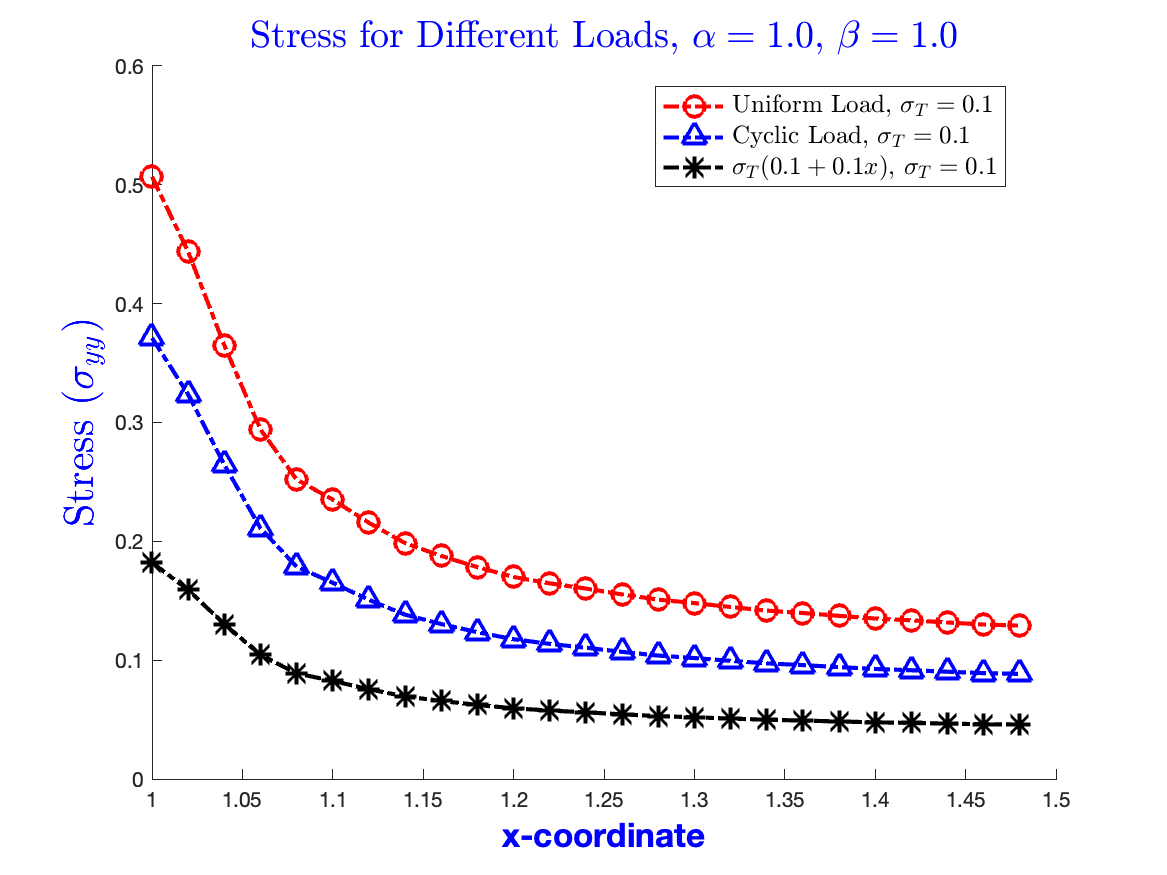}
        \caption{Load = 0.1}
        \label{fig:stress_M2_lp1}
    \end{subfigure}
    \caption{{Comparison of crack-tips stress for various loading scenarios with the material fibers aligned with $y$-axis.}}
    \label{stress_case4}
\end{figure}

\begin{figure}[H]
    \centering
    \begin{subfigure}[b]{0.32\linewidth}
        \centering
        \includegraphics[width=\linewidth]{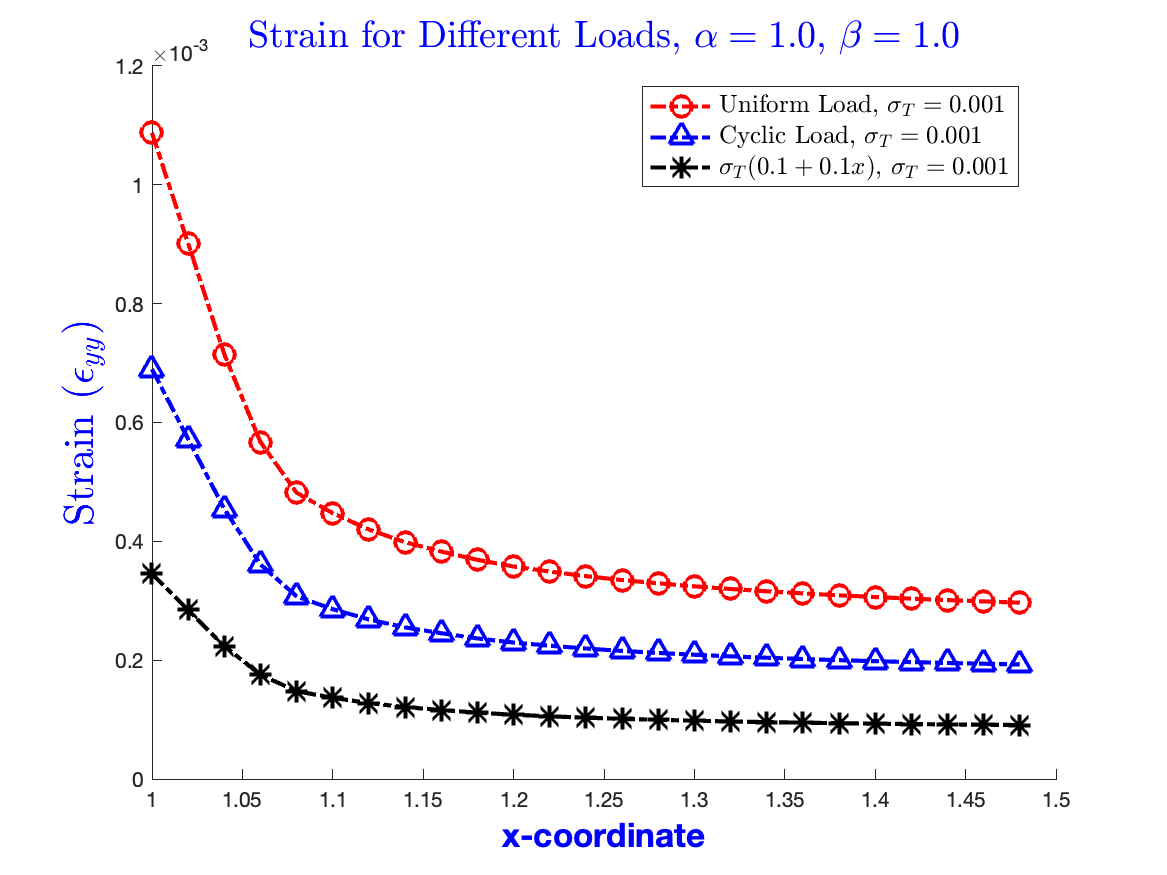}
        \caption{Load = 0.001}
        \label{fig:strain_M2_lp001}
    \end{subfigure}
    \hfill
    \begin{subfigure}[b]{0.32\linewidth}
        \centering
        \includegraphics[width=\linewidth]{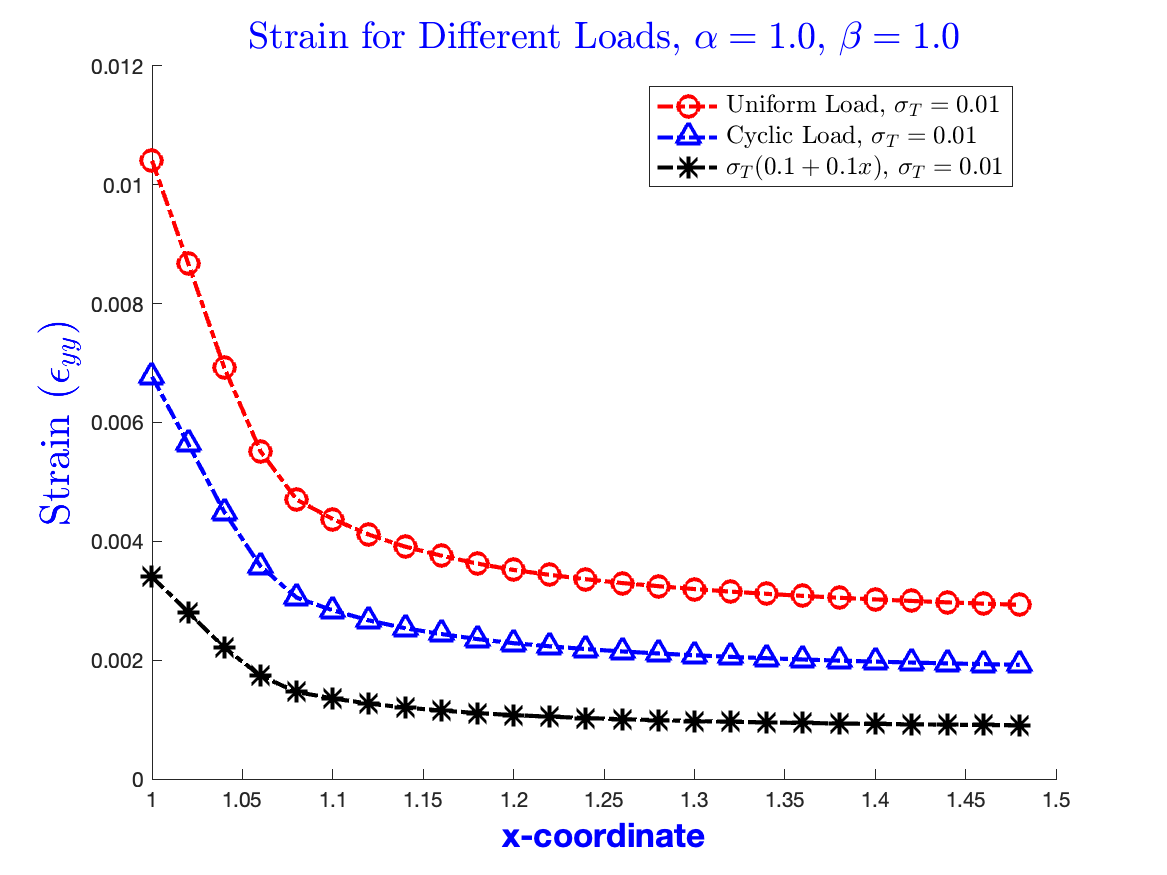}
        \caption{Load = 0.01}
        \label{fig:strain_M2_lp01}
    \end{subfigure}
    \hfill
    \begin{subfigure}[b]{0.32\linewidth}
        \centering
        \includegraphics[width=\linewidth]{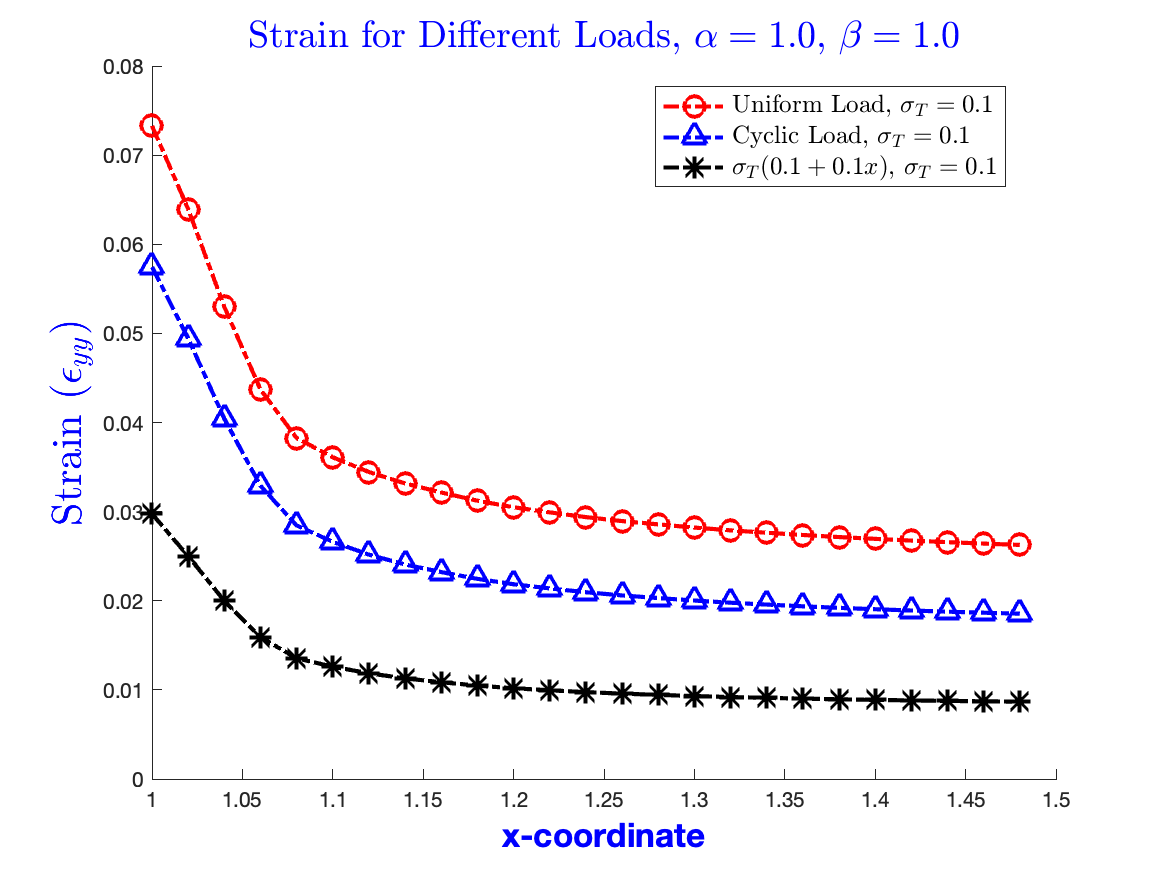}
        \caption{Load = 0.1}
        \label{fig:strain_M2_lp1}
    \end{subfigure}
    \caption{{Comparison of crack-tips strain for various loading scenarios with the material fibers aligned with $y$-axis.}}
    \label{strain_case4}
\end{figure}

Figures~\ref{stress_xcase3}, \ref{strain_xcase3}, \ref{stress_case4}, and \ref{strain_case4} illustrate the critical crack-tip stress and strain distributions under various loading conditions. These figures clearly demonstrate that a uniform tensile load consistently induces significantly higher stress and strain concentrations at the crack tip compared to the other two loading types investigated. This pronounced concentration under uniform tension highlights its substantial impact on material behavior near crack initiation points.

\section{Conclusion}
\label{sec:conclusion}

This study successfully introduced a novel mathematical model that meticulously characterizes the mechanical response of orthotropic elastic materials through a unique, algebraically nonlinear constitutive relationship. Our core objective was to deepen the understanding of crack-tip fields in strain-limiting bodies, specifically examining the influence of two distinct fiber orientations. The model's robustness is underpinned by the rigorous assumption of a monotonic and Lipschitz continuous constitutive response, which was crucial for ensuring the well-posedness of our continuous Galerkin formulation and, consequently, the existence and uniqueness of the weak solution. By synergistically combining these nonlinear constitutive laws with the fundamental balance of linear momentum, we meticulously formulated the complete boundary value problem. To effectively solve this complex, vector-valued, quasi-linear elliptic boundary value problem, we employed a sophisticated computational approach. This involved integrating Picard's iterative algorithm with the continuous conforming Galerkin finite element method, enabling us to achieve accurate and efficient numerical solutions for the intricate material interactions. Furthermore, our comprehensive sensitivity analysis revealed the significant impact of various model parameters, diverse loading conditions, and different fiber orientations on the material's behavior.
\begin{itemize}
\item A pivotal finding from this analysis is that, for both fiber orientations considered, an increase in the parameter $\beta$ consistently led to a moderate decrease in the peak values of stress, strain, and strain energy density at the crack tip. This inverse relationship unequivocally demonstrates that $\beta$ acts as a toughening or crack-mitigating parameter. This suggests that mechanisms governed by $\beta$, such as fiber bridging or localized plasticity, effectively shield the crack tip from applied loads. Ultimately, this means that increasing $\beta$ directly enhances the material's resilience to fracture by substantially reducing stress concentration at this critical point.
\item  In stark contrast to the beneficial effects of $\beta$, the other key modeling parameter, $\alpha$, has been found to exert a detrimental influence on the crack-tip fields, directly impacting the material's susceptibility to fracture and ultimate failure. An increase in the value of $\alpha$ consistently leads to a marked intensification of stress and strain concentrations at the crack tip. This critical observation implies a significant reduction in the material’s inherent ability to resist crack propagation as $\alpha$ increases. From a fracture mechanics standpoint, these elevated stress and strain values signify a concerning reduction in the external force required to reach the critical threshold for catastrophic failure. In essence, as $\alpha$ rises, the material becomes more brittle and prone to premature failure under lower applied loads. This heightened localized intensity at the crack tip accelerates the onset of ductile tearing or brittle fracture, depending on the material's specific failure mechanisms. Consequently, from a vital design and safety perspective, high values of $\alpha$ are not only undesirable but also inherently dangerous, as they dramatically increase the risk of unexpected and potentially catastrophic structural failure. This highlights $\alpha$ as a crucial parameter to minimize in the development of fracture-resistant materials and structures.
\item The strain energy density is consistently found to be highest in the immediate neighborhood of the crack-tip for both the parameters $\beta$ and $\alpha$. This localized concentration of strain energy provides a significant and robust local fracture criterion, crucial for thoroughly studying the evolution and propagation of crack tips under various loading conditions. Specifically, this criterion allows for the prediction of crack initiation and growth, as the material accumulates energy beyond its critical threshold. Its applicability extends across diverse scenarios, including pure mechanical loading \cite{lee2022finite,yoon2021quasi,manohar2025convergence} and complex coupled thermo-mechanical loads \cite{yoon2022CNSNS}. Furthermore, the magnitude and distribution of strain energy density can offer insights into the material's ductility or brittleness at the failure point, aiding in the development of more accurate predictive models for material integrity and lifespan.
\end{itemize}

The current research provides a robust foundation that can be significantly extended in several promising directions further to advance our understanding of crack behavior in complex materials. One crucial avenue involves investigating crack-tip fields within porous elastic bodies, where the material's moduli are intrinsically dependent on its density \cite{gou2025computational,yoon2024finite}. This extension would allow for more accurate modeling of materials with internal voids or cellular structures, which are prevalent in many engineering applications. Another crucial future undertaking is the challenge of resolving crack-tip fields in three-dimensional bodies \cite{gou2023finite,gou2023computational}. Moving beyond two-dimensional simplifications will enable a more realistic representation of crack propagation and stress distribution in actual components, providing insights vital for designing safer and more durable structures. Finally, a thorough numerical analysis comparing the performance of both continuous and discontinuous Galerkin finite element methods for these types of problems represents another important future research topic \cite{manohar2024hp}. Such an analysis would not only optimize computational efficiency but also enhance the accuracy and reliability of numerical predictions, contributing significantly to the field of computational fracture mechanics.

\section*{Acknoledgement}
The first author, SG, would like to thank the University of Texas Rio Grande Valley for providing a Presidential Research Fellowship during his PhD studies. SMM's work is supported by the National Science Foundation under Grant No. 2316905.
 
\bibliographystyle{plain}
\bibliography{ref}

\end{document}